\documentclass[1p, final]{elsarticle}

\usepackage{subfig}
\usepackage{amsmath,bm}
\usepackage{color}
\usepackage{enumerate}
\usepackage{multirow}
\usepackage{amssymb}
\usepackage{epstopdf}
\usepackage{epsfig}
\usepackage[table,xcdraw]{xcolor}
\usepackage{algpseudocode,algorithm,algorithmicx}

\usepackage{lineno,hyperref}
\modulolinenumbers[5]

\journal{Computational Geosciences}

\begin{document}

\begin{frontmatter}

\title{\LARGE Inexact Methods for Sequential Fully Implicit (SFI) Reservoir Simulation}

\author{Jiamin Jiang}
\cortext[mycorrespondingauthor]{Corresponding author}
\ead{jiamin.jiang@chevron.com}

\author{Pavel Tomin}

\author{Yifan Zhou}

\address{Chevron Energy Technology Co.}
\address{1500 Louisiana St., Houston, TX 77002, USA}

\date{November 9, 2020}

\begin{abstract}

The sequential fully implicit (SFI) scheme was introduced~(Jenny et al.~2006) for solving coupled flow and transport problems. Each time step for SFI consists of an outer loop, in which there are inner Newton loops to implicitly and sequentially solve the pressure and transport sub-problems. In standard SFI, the sub-problems are usually solved with tight tolerances at every outer iteration. This can result in wasted computations that contribute little progress towards the coupled solution. The issue is known as `over-solving'.~Our objective is to minimize the cost of inner solvers while maintaining the convergence rate of SFI.

We first extended a nonlinear-acceleration (NA) framework (Jiang and Tchelepi 2019) to multi-component compositional models, for ensuring robust outer-loop convergence. We then developed inexact-type methods that alleviate `over-solving'. It is found that there is no need for one sub-problem to strive for perfection, while the coupled (outer) residual remains high due to the other sub-problem.

The new SFI solver was tested using several complex cases. The problems involve multi-phase and EoS-based compositional fluid systems. We compared different strategies such as fixed relaxations on absolute and relative tolerances for the inner solvers, as well as an adaptive approach. The results show that the basic SFI method is quite inefficient. Away from a coupled solution, additional accuracy achieved in inner solvers is wasted, contributing to little or no reduction of the overall outer residual. By comparison, the adaptive inexact method provides relative tolerances adequate for the current convergence state of the sub-problems. We show across a wide range of flow conditions that the new solver can effectively resolve the over-solving issue, and thus greatly improve the overall efficiency.

\end{abstract}

\end{frontmatter}

\section{Introduction}

Numerical reservoir simulation is an essential tool for improving our understanding of underground resources, including oil and gas recovery, groundwater remediation, and $\mathrm{CO_2}$ subsurface sequestration. Detailed geological models with heterogeneous coefficients are usually constructed as input to reservoir simulators. In addition, predicting fluid dynamics evolution requires solving partial differential equations (PDEs) that represent multi-phase flow and transport in natural porous media. These PDEs are highly nonlinear and exhibit an intricate mixture of elliptic and hyperbolic characteristics. The above features make the development of robust, efficient, and accurate discretization and solution schemes quite challenging.

Several temporal discretization schemes are available to solve the mass-conservation equations that describe coupled flow and transport (Aziz and Settari 1979). The use of explicit schemes poses severe restrictions on time-step sizes and is usually considered impractical for heterogeneous three-dimensional problems, in which the Courant-Friedrichs-Lewy (CFL) numbers can vary by several orders of magnitude throughout the domain~(Coats 1980; Jenny et al. 2009). Therefore, implicit schemes such as fully implicit (FI) and sequential implicit (SI) methods are preferred in practice. The resulting nonlinear system is usually cast in residual form and solved using a Newton-based solver. For a target time-step, a sequence of nonlinear iterations is performed until the converged solution is achieved. Each iteration involves the construction of the Jacobian matrix and solution of the corresponding linear algebraic equations (Younis et al. 2010). 

Compared with the FI method, sequential methods handle flow and transport separately and differently. This adds flexibility in the choice of discretization scheme, solution strategy, and time-stepping for the sub-problems (M{\o}yner and Lie 2016; Kozlova et al.~2016; Lie et al. 2017). Specialized solvers that are well suited to the specific systems of equations can be designed to optimize performances. Acs et al. (1985), Watts~(1986), Trangenstein and Bell (1989a, 1989b) developed sequential implicit (SI) formulations for compositional flow problems. Among SI, the sequential fully implicit (SFI) method which converges to the same solution as FI, has recently received an increasing attention. The SFI method was introduced along with the development of the multiscale finite volume (MSFV) framework (Jenny et al. 2006), to simulate immiscible multi-phase problems. Each time step for SFI consists of an outer loop, in which there is one inner Newton loop to implicitly solve the pressure equation and another loop to implicitly solve the transport equations. SFI has been extended to the black-oil model (Lee et al. 2008; Watanabe et al. 2016; Kozlova et al.~2016; M{\o}yner and Lie 2016), and recently to general compositional models by Hajibeygi and Tchelepi~(2014), Moncorgé et al. (2018), and M{\o}yner and Tchelepi (2018).

For tightly coupled problems, the basic SFI method suffers from convergence difficulties, resulting in oscillations or divergence of the outer iterative process (Kuttler and Wall 2008; Jiang and Tchelepi 2019; Li et al. 2019; Moncorgé et al. 2019). Moncorgé et al. (2019) proposed a consistent scheme based on implicit hybrid upwinding (IHU) for both the flow and transport systems. Their results showed that IHU leads to large reductions of outer iterations on the problems with strong buoyancy or capillarity. Jiang and Tchelepi (2019) carried out a detailed analysis of the coupling mechanisms between the two sub-problems to better understand the iterative behaviors of SFI. They further developed a solution framework based on nonlinear acceleration (NA) techniques which greatly improve the convergence performance of outer loops.

In the standard SFI method, the sub-problems are usually solved to high precision at every outer iteration. It is found that, however, there is no need for one sub-problem to strive for perfection (`over-solving'), while the coupled (outer) residual remains high due to the other sub-problem. Over-solving may result in wasted computations that contribute little or no progress towards the coupled solution (Dembo et al. 1982; Eisenstat and Walker 1996; Klie 1996; Dawson et al. 1997; Klie and Wheeler 2005; Senecal and Ji 2017). Our objective here is to minimize the cost of inner solvers while not degrading the convergence rate of SFI.

In this work, we first extended the NA framework to compositional models, for ensuring robust outer-loop convergence.~We then developed inexact-type methods that alleviate `over-solving'. The motivation is similar to the inexact Newton method (Dembo et al. 1982; Eisenstat and Walker 1996; Klie 1996; Dawson et al.~1997), where the linear (i.e. inner) iterations are controlled in a way that the Newton (i.e.~outer) convergence is not degraded, but overall computational efforts are largely decreased. In particular, we proposed an~adaptive strategy that provides relative tolerances based on the convergence rates of coupled problems.

The new inexact SFI solver was tested using several complex cases.~The problems involve multi-phase and EoS-based compositional fluid systems. We compared different strategies starting from fixed relaxations on absolute and relative tolerances for the inner solvers. From the results we observe that the feedback from one inner solver can cause the residual of the other to rebound to a much higher level during outer iterations. When couplings are strong, the outer convergence is mainly restricted by the initial residuals of the sub-problems. The studies demonstrate that the basic SFI method is quite inefficient. Away from a coupled solution, additional accuracy achieved in inner solvers is wasted, contributing to little or no reduction of the overall residual. By comparison, the proposed adaptive inexact method provides relative tolerances adequate for the sub-problems. We show across a wide range of flow conditions that the new approach can effectively resolve the over-solving issue, and thus greatly improve the overall efficiency.

\section{Immiscible multi-phase flow}

\subsection{Governing equations}

We consider compressible and immiscible flow and transport with $n_p$ fluid phases. Pressure-dependent functions are incorporated to relate fluid volumes at reservoir and surface conditions. The conservation equation for phase $l$ is
\begin{equation} 
\label{eq:mass_con}
\frac{\partial \left ( \phi b_{l} s_{l} \right )}{\partial t } + \nabla \cdot \left (b_{l} \textbf{u}_{l} \right ) = b_{l} q_{l},
\quad l \in \left \{ 1,...,n_p \right \},
\end{equation}
where $\phi$ is the rock porosity, $t$ is the time, $b_{l}$ is the inverse of the phase formation volume factor (FVF), $s_{l}$ is the phase saturation, with the constraint that the sum of saturations is equal to one $\sum_{l} s_{l} = 1 $, $q_{l}$ is the well flow rate (source and sink terms). Without loss of generality, hereafter we ignore $q_{l}$ because it is zero everywhere except at a cell with well. $\textbf{u}_l$ is the phase velocity, which is expressed as a function of phase potential gradient $\nabla \Phi_l $ using the extended Darcy's law
\begin{equation} 
\label{eq:phase_vel}
\textbf{u}_l = -K \lambda_l \nabla \Phi_l = -K\lambda_l\left ( \nabla p - \rho_l g \nabla h \right ),
\end{equation}
where $K$ is the rock permeability. $p$ is the pressure (capillary forces are assumed to be negligible so that there is only a single pressure), $g$ is the gravitational acceleration and $h$ is the height. Phase mobility is given as $\lambda_{l} = k_{rl}/\mu_l$, where $k_{rl}$ and $\mu_l$ are the relative permeability and the viscosity, respectively. Phase density is evaluated through $\rho_l = b_l \rho_l^S $, where $\rho_l^S$ is the surface density.

The phase velocity can also be written using a fractional-flow formulation, which involves the total velocity $\textbf{u}_T $, defined as the sum of the phase velocities
\begin{equation} 
\label{eq:tol_vel}
\textbf{u}_T = \sum_l \textbf{u}_l = -K \lambda_T \nabla p + K \sum_l \lambda_l \rho_l g \nabla h .
\end{equation}
Eq. (\ref{eq:tol_vel}) is used to express the pressure gradient as a function of $\textbf{u}_T$ in order to eliminate the pressure variable from Eq. (\ref{eq:phase_vel}), obtaining
\begin{equation} 
\mathbf{u}_l = \frac{\lambda_{l}}{\lambda_T} \mathbf{u}_T + K g \nabla h \sum_{m} \frac{\lambda_{m}\lambda_{l}}{\lambda_T}  \left ( \rho_{l}-\rho_{m} \right).
\label{eq:phase_vel_frac}
\end{equation}

\subsection{Discretized equations}

The coupled multi-phase system describes the interplay between viscous and gravitational forces. A standard finite-volume scheme (Cao 2002) is employed as the spatial discretization for the conservation equations. A two-point flux approximation (TFPA) is used to approximate the flux at a cell interface.

The fully-implicit discretization of a cell is
\begin{equation} 
\label{eq:dis_mass}
\left ( \phi b_{l} s_{l} \right )^{n+1}_i - \left ( \phi b_{l} s_{l} \right )^{n}_i + \frac{\Delta t}{V_i} \sum_{j\in adj(i)} \left ( b_{l} F_{l} \right )^{n+1}_{ij} = 0 ,
\end{equation}
where subscript $i$ denotes quantities associated with cell $i$ and $ij$ denotes quantities associated with the interface between cells $i$ and $j$. $adj(i)$ is the set of neighbors of cell $i$. Superscripts denote the time level. $\Delta t$ is time-step size, $V_i$ is the volume of cell $i$.

It is often advantageous to reformulate Eq. (\ref{eq:dis_mass}) into one elliptic or parabolic equation for the pressure and several hyperbolic equations for the saturations. This reformulation allows us to employ discretization approaches specially developed and suited for the corresponding type of equations. Furthermore, a sequential solution procedure can be applied to reduce the number of unknowns in each solution step. Both the pressure and transport equations are nonlinear and need to be solved iteratively. 

To derive a discrete pressure equation, we first multiply Eq. (\ref{eq:dis_mass}) by $\alpha_l = 1/b_{l}^{n+1} $. Then with $\sum_l s_l = 1$, the summation of the resulting equations gives 
\begin{equation} 
\label{eq:press_eq}
\phi_i^{n+1} - \phi_i^{n} \sum_{l} \left ( \frac{b_{l}^{n}}{b_{l}^{n+1}}s_{l}^{n} \right )_i + \frac{\Delta t}{V_i} \sum_{j\in adj(i)} \sum_{l} \frac{b_{l,ij}^{n+1}}{b_{l,i}^{n+1}} F_{l,{ij}}^{n+1} = 0 ,
\end{equation}
where the saturation dependency at the current time level, $n+1$, is eliminated in the accumulation term. 

The numerical flux $F_{l,ij}$ can be written as 
\begin{equation} 
\label{eq:Fl_tol_v}
F_{l,ij} = \frac{\lambda_{l,ij}}{\lambda_{T,ij}} u_{T,ij} + \Upsilon_{ij} \sum_{m} \frac{\lambda_{l,ij} \lambda_{m,ij}}{\lambda_{T,ij}} \left ( g_{m,ij} - g_{l,ij} \right ),
\end{equation}
where $\lambda_{T,ij} = \sum_{l} \lambda_{l,ij}$. $\Upsilon_{ij}$ is the interface transmissibility. The total velocity discretization is given by
\begin{equation} 
\label{eq:disc_tol_v}
u_{T,ij} = \sum_{l} \Upsilon_{ij} \lambda_{l,ij} \Delta \Phi_{l,ij} ,
\end{equation}
where $\Delta \Phi_{l,ij} = \Delta p_{ij} - g_{l,ij}$ is the phase potential difference with the discrete weights $g_{l,ij} = \rho_{l,ij} \, g \Delta h_{ij}$.

The mobility terms in $F_{l,ij}$ are usually evaluated using upwinding schemes. The~phase-potential upwinding (PPU) (Aziz and Settari 1979) and the phase upwinding~(PU) (Brenier and Jaffré 1991) are popular in reservoir simulation practice. As revealed by Li and Tchelepi (2015), Lee et al. (2015), Jiang and Younis (2017), these schemes can produce switches of upstream directions, thus causing nonlinear convergence difficulties in the presence of buoyancy. To address this flow reversal issue, hybrid upwinding (HU) scheme was proposed and extended for obtaining a smooth numerical flux (Lee et al. 2015; Lee and Efendiev 2016). In addition, Jiang and Younis~(2017) devised an alternative scheme called C1-PPU, which improves smoothness with respect to saturations as well as phase potentials. 

We recently demonstrated that the discontinuous behavior of PU also largely degrades the outer-loop convergence of SFI. By comparison, the HU scheme alleviates the convergence difficulty because it reduces couplings between the flow and transport sub-problems (Jiang and Tchelepi 2019; Moncorgé et al. 2019). In this paper we employ HU for the transport problem.

\subsection{Sequential formulation}

The Sequential Fully Implicit (SFI) method (Jenny et al. 2006) has received increasing attention in recent years. Each time step for SFI consists of an outer loop, in which there are inner Newton loops to implicitly and sequentially solve the pressure and transport sub-problems. For each iteration of outer loop, computations proceed as follows: compute the pressure field iteratively to a certain tolerance, and update the total flux, then compute the saturation iteratively. Note that the total flux is always evaluated using the PPU scheme. The updated saturation defines a new mobility field for the subsequent pressure problem. These steps can be iterated until convergence of all variables at the current time level.

The equation system under the fractional-flow formulation will consist of the pressure equation, Eq.~(\ref{eq:press_eq}), and $n_p-1$ transport equations, Eq.~(\ref{eq:dis_mass}). Correspondingly, the primary variables are the pressure and $n_p-1$ phase saturations. 

The solution procedure of SFI for a single time-step is demonstrated in Algorithm~\ref{alg:SFI}. The counter $\nu$ denotes an outer iteration, which should not be confused with the inner iterations over the individual solvers. $\epsilon_p^{out}$~and $\epsilon_t^{out}$ are the increment tolerances of pressure and saturations between outer iterations, respectively. $r_p$ and $r_t$ are the residuals of the pressure and transport equations. $\mathcal{R}_{p}$ and $\mathcal{R}_{t}$ are the normalized residuals. $\epsilon_p$~and $\epsilon_t$ are absolute tolerances of residual norm for the pressure and transport solvers. $J_p = \left [ \frac{\partial r_p}{\partial p} \right ] $ and $J_s = \left [ \frac{\partial r_s}{\partial s} \right ] $ are the Jacobian matrices.



\begin{algorithm}
	\caption{Sequential fully implicit method} \label{alg:SFI}
	\begin{algorithmic}[1]
		\State $\nu = 0$, initialize $p^{\nu} = p^n$, $s^{\nu} = s^n$ 
		\While{$\left \| p^{\nu} - p^{\nu-1} \right \|_{\infty} > \epsilon_p^{out} \ , \ \left \| s^{\nu} - s^{\nu-1} \right \|_{\infty} > \epsilon_t^{out}$}
		\Comment{Outer coupling loop}
		\smallskip
		\State $k = 0$, $p_{k} = p^{\nu}$ 
		\While{$\left \| \mathcal{R}_{p} \right \|_{\infty} > \epsilon_p$ } \Comment{Pressure loop}
		\State Solve linearized pressure equation:
		\State $J_p \delta p = -r_p$
		\State $p_{k + 1} = p_{k} + \delta p$
		\State $k \leftarrow k + 1$
		\EndWhile 
		\State $p^{\nu + 1} = p_{k}$
		
		\medskip
		\State Compute total flux by summing phase fluxes
		\medskip
		
		\State $k = 0$, $s_{k} = s^{\nu}$ 
		\While{$\left \| \mathcal{R}_{t} \right \|_{\infty} > \epsilon_t$ } \Comment{Transport loop}
		\State Solve linearized transport equations:
		\State $J_s \delta s = -r_s$
		\State $s_{k + 1} = s_{k} + \delta s$
		\State $k \leftarrow k + 1$
		\EndWhile 
		
		\State $s^{\nu + 1} = s_{k}$
		\State $\nu \leftarrow \nu + 1$
		
		\EndWhile 
	\end{algorithmic}
\end{algorithm}

\section{Compositional flow}

\subsection{Governing equations}

An important aspect of any compositional modeling is the choice of formulation. Two commonly used formulations are natural variables (Coats 1980; Cao 2002) and overall-composition variables (Acs et al. 1985; Collins et al. 1992; Voskov and Tchelepi 2012). In this work, we rely on the overall-composition formulation which avoids variable substitution. The equations and unknowns are the same for any phase state, and the variables set is

(1) $p$ $-$ pressure, 

(2) $z_c$ $-$ overall composition of each component.

The overall-composition formulation greatly simplifies the application of the nonlinear acceleration described in Section~\ref{NA} which can be readily used to the formulation with a~consistent variables set.

We consider compressible gas-oil flow without capillarity. We ignore water that does not exchange mass with the hydrocarbon phases. The conservation equations for the isothermal compositional problem containing $n_c$ components are written as 
\begin{equation} 
\label{eq:mass_con_comp}
\frac{\partial}{\partial t } \left ( \phi z_c \rho_T \right ) + \nabla \cdot \left ( x_c \rho_{o} \textbf{u}_{o} + y_c \rho_{g} \textbf{u}_{g} \right ) = q_{c}, \quad c \in \left \{ 1,...,n_c \right \},
\end{equation}
where $\rho_T=\sum_l \rho_l s_l$ is the total density, $\rho_l$ is the phase molar density, $x_{c}$ and $y_{c}$ are the mole fractions of component $c$ in the oil and gas phases, respectively, $q_{c}$ is the well flow rate. Phase velocities are defined by Eq.~(\ref{eq:phase_vel}) or Eq.~(\ref{eq:phase_vel_frac}).

For a Newton iteration, a phase stability test is performed to determine if the system can be split into two phases. If a cell is determined to be in a single-phase state, no additional computations are necessary. If a two-phase state is detected, flash calculations are performed to obtain the phase compositions $x_{c}$, $y_{c}$ and molar phase fractions $\nu_o$, $\nu_g$ by solving a local nonlinear system of equations.
For a mixture of $n_c$ components and two phases, the mathematical model expressing the thermodynamic equilibrium is (Voskov and Tchelepi 2012) 

\begin{equation} 
f_{c,o}(p,x_1,..,x_{n_c}) - f_{c,g}(p,y_1,..,y_{n_c}) = 0 , \quad c \in \left \{ 1,...,n_c \right \},
\label{eq:vle_fu}
\end{equation}

\begin{equation} 
z_c - \nu_o x_{c} - \nu_g y_{c} = 0 , \quad c \in \left \{ 1,...,n_c \right \},
\end{equation}

\begin{equation} 
\sum_{c=1}^{n_c} \left ( x_{c} - y_{c} \right ) = 0.
\label{eq:vle_sum}
\end{equation}
with additional constraints used to eliminate a part of the unknowns. These include the phase constraints
\begin{equation} 
	\sum_{c=1}^{n_c} x_{c} = 1, \qquad \sum_{c=1}^{n_c} y_{c} = 1,
\end{equation}
and the phase fraction constraint
\begin{equation} 
	\nu_o + \nu_g = 1.
\end{equation}

We assume that $p$ and $\textbf{z}$ are known, and that $f_{c,l}$ are governed by an equation of state (EoS) model. The objective is to find all the $x_{c}$, $y_{c}$, and $\nu_l$. We employ a two-stage procedure to compute the equilibrium phase behavior. The standard stability analysis is performed first to determine if a single-phase mixture is likely to split into two phases. Then flash calculations are performed to obtain the compositions of the existing phases. After Eq.~(\ref{eq:vle_fu}-\ref{eq:vle_sum}) are solved, the derivatives of $x_{c}$, $y_{c}$, and $\nu_l$ with respect to $p$ and $\textbf{z}$ are obtained using the inverse theorem approach (Collins et al. 1992; Voskov and Tchelepi 2012). The saturations are then computed as $s_l = \frac{\nu_l / \rho_l}{\sum_l \nu_l / \rho_l}$.

\subsection{Sequential formulation}

Acs et al. (1985), Watts (1986), Trangenstein and Bell (1989a, 1989b) developed sequential implicit (SI) formulations for compositional flow problems. The sequential method involves solving a pressure equation and a set of transport equations.

The pressure equation that is essentially a volume balance is derived from a weighted sum of the component conservation equations
\begin{equation} 
r_p = \sum_{c} r_c w_c,
\end{equation}
where $r_c$ and $w_c$ are the conservation equation and weighting factor of component $c$, respectively. The purpose of the weighting factors is to eliminate the dependency of pressure accumulation term with respect to non-pressure variables. They are found from the solution of a local linear system for each cell at each pressure iteration. We follow the algebraic approach proposed by Coats (2000), Moncorgé et al. (2018), M{\o}yner and Tchelepi (2018) to construct the pressure equation. The overall compositions, $z_c$, are kept fixed during the pressure solution. 

For the transport, either one conservation equation or the volume closure equation must be dropped to avoid an over-determined system. Here we remove the volume closure, because it is more practical for the cases with significant compressibility or phase-change effects (M{\o}yner and Tchelepi 2018). The total saturation $s_T$ is used as an additional variable, $s_T$ appears in the accumulation terms and as a phase flux multiplier (Watts 1986, M{\o}yner and Tchelepi 2018). The resulting transport system contains $n_c$ conservation equations, and the variable set to be solved is $\left ( z_1 , ... , z_{n_c-1}, s_T \right )$. With the $s_T$ formulation in transport, a volume error is introduced such that $s_T \neq 1$, which can be interpreted as the splitting error: $s_T$ deviates from unity in cells with large coupling errors (M{\o}yner and Tchelepi 2018).

The pressure and total-velocity are fixed during the transport solution. The spatial and temporal discretizations for the sub-problems follow the same methods as in the immiscible multi-phase model.

\section{Nonlinear acceleration}
\label{NA}

Previously we have shown that direct applications of SFI iterations may encounter severe convergence difficulties (Li et al. 2019). Therefore, nonlinear acceleration techniques are necessary to improve the convergence of SFI. We apply the quasi-Newton~(QN) method for accelerating the convergence of SFI (Jiang and Tchelepi 2019). The objective is to reduce outer iteration counts.

It is worth noting that QN-type methods have been applied to improve the nonlinear convergence of Newton-Krylov methods (Klie 1996; Dawson et al.~1997; Klie and Wheeler 2005).~But their focus is on solving the fully-implicit systems for coupled flow and transport problems.

\subsection{Block Gauss-Seidel formulation}

The iterative form of SFI is equivalent to the block Gauss-Seidel (BGS) process
\begin{equation} 
\label{eq:n_GS}
\begin{cases}
\textbf{x}_p^{\nu+1} = \mathcal{P} \left ( \textbf{x}_p^{\nu},\textbf{x}_t^{\nu} \right ) , \\ 
\textbf{x}_t^{\nu+1} = \mathcal{T} \left ( \textbf{x}_p^{\nu+1},\textbf{x}_t^{\nu} \right ) .
\end{cases}
\end{equation}
where the operators $\mathcal{P}$ and $\mathcal{T}$ represent respectively the pressure and transport solvers. $\textbf{x}_p$ and $\textbf{x}_t$ represent vectors of the pressure and transport unknowns. The transport solver works on the results from the pressure solver. Note that the coupling is also subject to the transmission condition with fixed total flux. The BGS coupling involves repeated applications of the update given by Eq.~(\ref{eq:n_GS}). The global problem is converged when the solution to Eq.~(\ref{eq:n_GS}) is consistent and both the sub-problems are converged.

The BGS iteration can be written in compact form as
\begin{equation} 
\label{eq:c_n_GS}
\widetilde{\textbf{x}_t}^{\nu+1} = \mathcal{T} \left ( \textbf{x}_p^{\nu+1}, \textbf{x}_t^{\nu} \right ) = \mathcal{T} \Big ( \mathcal{P} \left ( \textbf{x}_t^{\nu} \right ), \textbf{x}_t^{\nu} \Big ) = \mathcal{T} \circ \mathcal{P} \left ( \textbf{x}_t^{\nu} \right ).
\end{equation}
where a tilde $\widetilde{(\cdot )}$ denotes the current unmodified solution during the iteration. The nonlinear operator $\mathcal{T} \circ \mathcal{P} $ takes an input vector $\textbf{x}_t^{\nu}$ and generates an output vector $\widetilde{\textbf{x}_t}^{\nu+1}$ of the same size.

To enable an implicit treatment, a residual form is introduced
\begin{equation} 
\label{eq:resi_n_GS}
r^{\nu+1} = \widetilde{\textbf{x}_t}^{\nu+1} - \textbf{x}_t^{\nu} = \mathcal{T} \circ \mathcal{P} \left ( \textbf{x}_t^{\nu} \right ) - \textbf{x}_t^{\nu} .
\end{equation}
Iterative correction is required to bring the pressure and saturation responses into balance so that the residual vanishes. 

From Eq. (\ref{eq:resi_n_GS}) we see that convergence criteria of SFI can be readily determined based on solution increments between outer iterations. The increments of primary variables are checked to ensure the convergence of the coupled problem.

\subsection{Quasi-Newton nonlinear acceleration}

Once the system is reformulated in the BGS form, it can be treated as a root-finding problem, Eq.~(\ref{eq:resi_n_GS}), to be tackled by the Newton method. In order to determine the new increment $\Delta \textbf{x}_t^{\nu}$, a linear equation system has to be solved
\begin{equation} 
J^{\nu} \Delta \textbf{x}_t^{\nu} = -r^{\nu+1} ,
\end{equation}
where $J^{\nu}$ denotes the Jacobian matrix of the residual operator. $\Delta \textbf{x}_t^{\nu}$ is the difference between the current and previous solutions. Then the new modified solution is 
\begin{equation} 
\label{eq:1_Aitken}
\textbf{x}_t^{\nu+1} = \textbf{x}_t^{\nu} + \Delta \textbf{x}_t^{\nu} .
\end{equation}

The QN method approximates the Jacobian directly from generated nonlinear vector sequence (Degroote et al. 2010). The QN update is given as 
\begin{equation} 
\Delta \textbf{x}_t^{\nu} = \left ( \widehat{\left. \frac{\partial r}{\partial \textbf{x}_t} \right|_{\textbf{x}_t^{\nu}}} \right )^{-1} \left ( -r^{\nu+1} \right ) .
\end{equation}

The inverse of the Jacobian does not have to be constructed explicitly; we only need the product of it with the right-hand side vector. $\Delta \textbf{x}_t^{\nu}$ is directly approximated through solving an unconstrained form of the least-squares problem 
\begin{equation} 
\label{eq:unc_ls}
\underset{\gamma}{\textrm{min}} \left \| r^{\nu+1} - \Delta R^{\nu} \gamma \right \|_2 ,
\end{equation}
where the solution vector is denoted by $\gamma =\left ( \gamma_0,...,\gamma_{m_{\nu}-1} \right )^T $, $\left \| \ \cdot \ \right \|_2 $ is the Euclidean norm on $\mathbb{R}^n$, matrix $\Delta R^{\nu} = \left ( \Delta r^{\nu-m_{\nu}+1}, ..., \Delta r^{\nu} \right ) $ is constructed by stacking vectors $\Delta r^i = r^{i+1} - r^{i} $.

The QN method is adapted to SFI for coupled flow and transport.~The algorithmic description is given in Algorithm \ref{alg:QN_SFI}. Note that $\widetilde{\textbf{x}}_t^{\nu+1} = \mathcal{T} \circ \mathcal{P} \left ( \textbf{x}_t^{\nu} \right )$ is one SFI iteration, comprising the inner loops described in Algorithm \ref{alg:SFI}.~Input and output vectors from maximum $m$ previous iterations are combined to provide a better prediction for the next iteration.~Parameter $m$ defines the number of previous iterations used to obtain the secant information, and $m_{\nu} = \textrm{min} \left \{ m, \nu \right \}$.~If $m$ is small, then the secant information may be too limited to provide desirably fast convergence.~However, if $m$ is large, the least-squares problem may be badly conditioned.~Moreover, outdated secant information from previous iterations may be kept, leading to convergence degradation (Walker and Ni 2011). Therefore a proper choice of $m$ is likely to be problem-dependent.

The least-squares problem with the unconstrained form (\ref{eq:unc_ls}) can be solved using QR decomposition (Walker and Ni 2011). Only the economy (thin) QR decomposition is necessary. As the algorithm proceeds, the QR decomposition of $\Delta R^{\nu}$ can be efficiently achieved through updating that of $\Delta R^{\nu-1}$: $\Delta R^{\nu}$ is obtained from $\Delta R^{\nu-1}$ by appending a new column on the right and possibly dropping one column from the left. If the QR decomposition is $\Delta R^{\nu} = \mathcal{Q}^{\nu} \times \mathcal{R}^{\nu}$, then the solution of Eq.~(\ref{eq:unc_ls}) can be expressed as
\begin{equation} 
\gamma = \textrm{arg} \ \underset{\gamma}{\textrm{min}} \left \| \left ( \mathcal{Q}^{\nu} \right )^T r^{\nu+1} - \mathcal{R}^{\nu} \gamma \right \|_2 ,
\end{equation}
which is obtained by solving the $m_{\nu} \times m_{\nu}$ triangular system $\mathcal{R}^{\nu} \gamma = \left ( \mathcal{Q}^{\nu} \right )^T r^{\nu+1} $. After $\gamma$ is found, $\Delta \textbf{x}_t^{\nu}$ can be computed as
\begin{equation}
\Delta \textbf{x}_t^{\nu} = r^{\nu+1} - \left ( \Delta X^{\nu} + \Delta R^{\nu} \right ) \gamma^{\nu},
\end{equation}
where $\Delta X^{\nu} = \left ( \Delta \textbf{x}_t^{\nu-m_{\nu}}, ..., \Delta \textbf{x}_t^{\nu-1} \right ) $, $\Delta \textbf{x}_t^{\nu} =\textbf{x}_t^{\nu+1} - \textbf{x}_t^{\nu}$.

A damping parameter $\omega$ is used to reduce step lengths when iterates are not near a~solution. The update with damping is given as
\begin{equation}
\textbf{x}_t^{\nu+1} = \textbf{x}_t^{\nu} + \omega \, r^{\nu+1} - \left ( \Delta X^{\nu} + \omega \Delta R^{\nu} \right ) \gamma
\end{equation}
We found that $\omega \in [0.5, 0.8]$ is quite effective for the cases with strong couplings between flow and transport. The value of $\omega$ is fixed throughout a whole simulation.

\begin{algorithm}
	\caption{Quasi-Newton for the SFI method} \label{alg:QN_SFI}
	\begin{algorithmic}[1]
		\State Set $\omega^{0} $, $\omega$ and $m \geq 1$
		\State Initialize $\textbf{x}_p^{0} = \textbf{x}_p^n$, $\textbf{x}_t^{0} = \textbf{x}_t^n$ 
		\State $\widetilde{\textbf{x}_t}^{1} = \mathcal{T} \circ \mathcal{P} \left ( \textbf{x}_t^{0} \right ) $
		\State $r^{1} = \widetilde{\textbf{x}_t}^{1} - \textbf{x}_t^{0} $
		\State $\textbf{x}_t^{1} = \textbf{x}_t^{0} + \omega^{0} r^{1} $
		\State $\nu = 1$
		\While{$\left \| \textbf{x}_p^{\nu} - \textbf{x}_p^{\nu-1} \right \|_{\infty} > \epsilon_p^{out} \ , \ \left \|  \textbf{x}_t^{\nu} - \textbf{x}_t^{\nu-1} \right \|_{\infty} > \epsilon_t^{out}$}
		\Comment{Outer coupling loop}
		\smallskip
		\State Set \ $m_{\nu} = \textrm{min} \left \{ m, \nu \right \}$ 
		\State $\widetilde{\textbf{x}}_t^{\nu+1} = \mathcal{T} \circ \mathcal{P} \left ( \textbf{x}_t^{\nu} \right ) $
		\Comment{One SFI iteration}
		\State $r^{\nu+1} = \widetilde{\textbf{x}}_t^{\nu+1} - \textbf{x}_t^{\nu} $
		\State $\Delta X^{\nu} = \left ( \Delta \textbf{x}_t^{\nu-m_{\nu}}, ..., \Delta \textbf{x}_t^{\nu-1} \right ) $, where $\Delta \textbf{x}_t^i =\textbf{x}_t^{i+1} - \textbf{x}_t^{i}$ 
		\State $\Delta R^{\nu} = \left ( \Delta r^{\nu-m_{\nu}+1}, ..., \Delta r^{\nu} \right ) $, where $\Delta r^i = r^{i+1} - r^{i} $ 
		\State $\gamma = \left ( \left ( \Delta X^{\nu} \right )^T \Delta R^{\nu} \right )^{-1} \left ( \Delta X^{\nu} \right )^T r^{\nu+1} $
		\State $\textbf{x}_t^{\nu+1} = \textbf{x}_t^{\nu} + \omega \, r^{\nu+1} - \left ( \Delta X^{\nu} + \omega \Delta R^{\nu} \right ) \gamma$ 
		\State $\nu \leftarrow \nu + 1$
		\smallskip
		\EndWhile 
	\end{algorithmic}
\end{algorithm}

The nonlinear acceleration provides a general approach both for immiscible and compositional models. For immiscible case the variable set is $\textbf{x}_t = \left ( s_1 , ... , s_{n_p-1}\right )$. When applying QN to the compositional model, $\textbf{x}_t = \left ( z_1 , ... , z_{n_c-1}, s_T \right )$ is taken as the solution vector for the SFI iteration.

\section{Inexact methods for SFI}



\begin{figure}[!htb]
\centering
\includegraphics[scale=0.5]{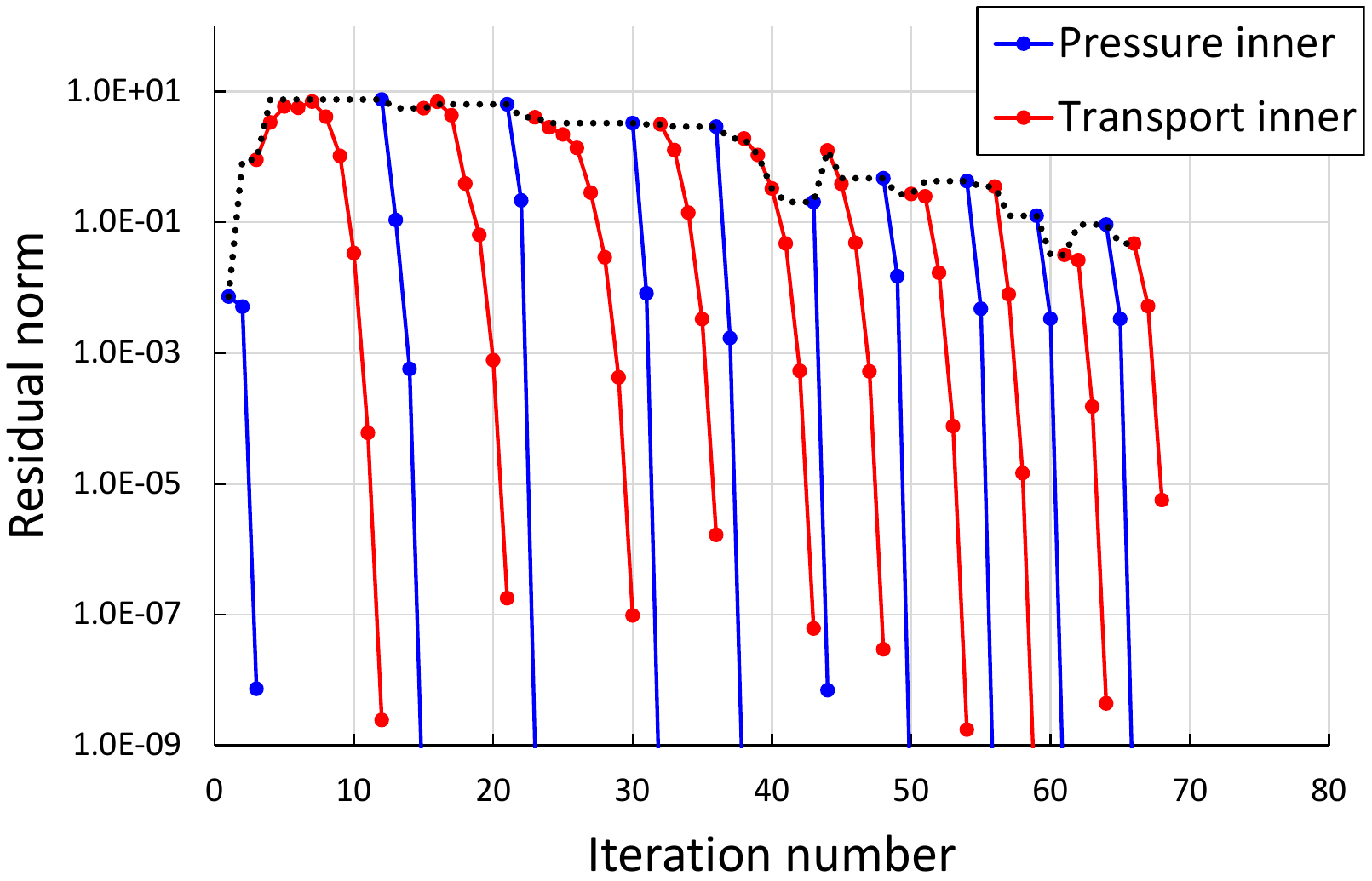}
\caption{Residual history of the inner solvers with tight tolerances for a 2D gravity-driven two-phase problem. Dotted line shows the outer loop convergence envelop.}
\label{fig:in_re_plot}
\end{figure}

The main objective of this work is to minimize the cost of inner solvers while not degrading the global convergence rate of SFI. In standard SFI, sub-problems are usually solved to high precision at every outer iteration. That over-solving may result in wasted computations contributing little progress towards the coupled solution (Dembo~et~al.~1982; Eisenstat and Walker 1996; Senecal and Ji 2017). To demonstrate this issue, we examine a test case described in Section 6.1, which is a two-dimensional gravity-driven two-phase problem. The tight convergence tolerances $\epsilon_p = 10^{-7}$ and $\epsilon_t = 10^{-5}$ are used for the inner solvers. The residual history of the inner solvers is plotted in \textbf{Fig.~\ref{fig:in_re_plot}}. 

A residual rebound is observed at the beginning of each inner loop, due to the feedback from the other solver. The magnitude of the rebound is one heuristic measure of the coupling strength. We can see that the outer convergence is determined by the initial residuals rather than the final residuals, that is, the top rather than the bottom envelope of the curves. Consequently, there is no need for one solver to strive for perfection, while the coupled residual remains high due to the other solver.

To mitigate the over-solving issue, relaxations on the convergence tolerances of the inner solvers are necessary. The motivation here is similar to the inexact Newton method (Dembo et al. 1982; Eisenstat and Walker 1996; Klie 1996; Dawson et al. 1997), where the linear (i.e. inner) iterations are controlled in a way that the Newton (i.e. outer) convergence is not degraded, but overall computational efforts are largely decreased (Senecal and Ji 2017). Next we present three options of inexact methods for SFI.

\subsection{Absolute relaxation}

Perhaps the simplest method is to perform fixed relaxations on absolute tolerances $\epsilon_p$ and $\epsilon_t$. The disadvantage of the method is its lack of generality. Tuning of parameters is required for different cases to provide a robust performance.

\subsection{Relative relaxation}

The classical inexact Newton (IN) method is a generalization of Newton's method for solving a nonlinear problem, $r(x) = 0$. At the $k$-th iteration of IN, the step $\delta x_{k}$ from the current solution $x_{k}$ has to meet the condition (Dembo et al. 1982)
\begin{equation} 
\label{eq:ada_rel_1}
\left \| r(x_{k}) + J(x_{k}) \delta x_{k} \right \| \leq \xi_{k} \left \| r(x_{k}) \right \|
\end{equation}
for a so-called forcing term $\xi_{k} \in [0, 1)$. Away from a solution, choosing $\xi_{k}$ too small may lead to over-solving the Newton equation. Therefore a less accurate approximation of the Newton step could be more effective. The proper choice of the forcing term in Eq.~(\ref{eq:ada_rel_1}) is critical to the efficiency and robustness of the IN method.

As in the IN method, each inner solver of SFI can be assigned a relaxed relative tolerance. The relative termination criterion is given as (Birken 2015)
\begin{equation} 
\label{eq:rel_rel_1}
\frac{\left \| \mathcal{R}_{k} \right \|_{\infty}}{\left \| \mathcal{R}_{0} \right \|_{\infty}} < \xi 
\end{equation}
or
\begin{equation} 
\label{eq:rel_rel_2}
\left \| \mathcal{R}_{k} \right \|_{\infty} < \epsilon \ , \quad \textrm{with} \ \ \epsilon = \xi \left \| \mathcal{R}_{0} \right \|_{\infty} ,
\end{equation}
where $\mathcal{R}_{k}$ is the normalized residual at inner iteration $k$. Here the value of $\xi$ is fixed during the simulation. Based on our experience, reducing the residual of an inner solver by only one order of magnitude per outer iteration provides acceptable performance.

\subsection{Adaptive relative relaxation}

The relative tolerance $\xi$ in SFI can be adaptively computed following a strategy similar to the one proposed by Eisenstat and Walker (1996). A more general form of the termination criterion for the inner solvers is given by 
\begin{equation} 
\label{eq:crr_p_1}
\left \| \mathcal{R}_{k}^{\nu} \right \|_{\infty} < \epsilon \ , \quad \textrm{with} \ \ \epsilon = \xi^{\nu} \left \| \mathcal{R}_{0}^{\nu} \right \|_{\infty} ,
\end{equation}
where the indices $\nu$ and $k$ refer to the outer and inner iterations, respectively. Thus the value of $\xi^{\nu}$ in not constant and is computed by
\begin{equation} 
\label{eq:crr_p_2}
\xi^{\nu} = \beta \, \frac{\left \| \mathcal{R}_{0}^{\nu} \right \|_{\infty}}{\left \| \mathcal{R}^{\nu-1}_{0} \right \|_{\infty}}  \ , \ \quad \nu>0 .
\end{equation}

As we can see, $\xi^{\nu}$ is chosen based on the coupled problem's convergence rate, which is estimated by the ratio of the current and previous initial residual norm. We set $\xi^{0} = \beta$. After the first outer iteration, the estimate for the convergence rate is available. The relative tolerances are restricted to the range $\xi \in \left [0.01 , \beta \right ]$. 

Eqs. (\ref{eq:crr_p_1}) and (\ref{eq:crr_p_2}) are separately applied to each inner solver, with possibly different values of $\beta$. For each outer iteration, the initial residual norms $\left \| \mathcal{R}_{0}^{\nu} \right \|_{\infty}$ are recorded, and $\xi^{\nu}$ are updated at the beginning of each inner solver. The algorithmic description of the adaptive relative relaxation method is given in Algorithm \ref{alg:SFI_in}. 

It should be noted that solution qualities of the inner solvers, especially the pressure solutions, may require higher accuracy. This is because the total-flux field obtained from pressure could have a large impact on the nonlinear convergence of the transport solver. In practice, it may be necessary to impose additional safeguards that restrict the absolute tolerance $\epsilon$ to a value range for robustness.

\begin{algorithm}
\caption{Adaptive relative relaxation method} \label{alg:SFI_in}
\begin{algorithmic}[1]
\State $\nu = 0$, initialize $p^{\nu} = p^n$, $\textbf{x}_t^{\nu} = \textbf{x}_t^n$ 
\smallskip
\State $\xi_p^{0} = \beta_p \ , \ \xi_s^{0} = \beta_s$
\smallskip
\While{$\left \| p^{\nu} - p^{\nu-1} \right \|_{\infty} > \epsilon_p^{out} \ , \ \left \| \textbf{x}_t^{\nu} - \textbf{x}_t^{\nu-1} \right \|_{\infty} > \epsilon_t^{out}$}
\Comment{Outer coupling loop}
\medskip
\State $k = 0$, $p_{k} = p^{\nu}$ 
\smallskip
\State $\epsilon_p = \xi_p^{\nu} \left \| \mathcal{R}_{p,k}^{\nu} \right \|_{\infty} $
\smallskip
\While{$\left \| \mathcal{R}_{p} \right \|_{\infty} > \epsilon_p$ } \Comment{Pressure loop}
\State Solve linearized pressure equation:
\State $J_p \delta p = -r_p$
\State $p_{k + 1} = p_{k} + \delta p$
\State $k \leftarrow k + 1$
\EndWhile 
\State $p^{\nu + 1} = p_{k}$

\medskip
\State Compute total flux by summing phase fluxes
\medskip

\State $k = 0$, $\textbf{x}_{t,k} = \textbf{x}_{t}^{\nu}$ 
\smallskip
\State $\epsilon_t = \xi_s^{\nu} \left \| \mathcal{R}_{s,k}^{\nu} \right \|_{\infty} $
\smallskip
\While{$\left \| \mathcal{R}_{t} \right \|_{\infty} > \epsilon_t$ } \Comment{Transport loop}
\State Solve linearized transport equations:
\State $J_t \delta \textbf{x}_{t} = -r_t$
\State $\textbf{x}_{t, k + 1} = \textbf{x}_{t, k} + \delta \textbf{x}_{t}$
\State $k \leftarrow k + 1$
\EndWhile 
\State $\textbf{x}_t^{\nu + 1} = \textbf{x}_{t,k}$
\smallskip
\State $\nu \leftarrow \nu + 1$
\smallskip
\State $\xi_p^{\nu} = \beta_p \frac{\left \| \mathcal{R}_{p,0}^{\nu} \right \|_{\infty} }{ \left \| \mathcal{R}^{\nu-1}_{p,0} \right \|_{\infty} }$, \, $\xi_t^{\nu} = \beta_t \frac{ \left \| \mathcal{R}_{t,0}^{\nu} \right \|_{\infty} }{ \left \| \mathcal{R}^{\nu-1}_{t,0} \right \|_{\infty} } $
\smallskip
\EndWhile 
\end{algorithmic}
\end{algorithm}

\section{Results: immiscible multi-phase flow}

We validated the effectiveness of the inexact SFI framework using the immiscible multi-phase problems. The results of the outer-loop convergence for the immiscible problems are not reported; comprehensive comparisons between basic SFI and NA techniques were provided in the previous work (Jiang and Tchelepi 2019).

All the numerical studies in this section were conducted within the open-source Matlab Reservoir Simulation Toolbox (MRST) software (Lie et al. 2012; M{\o}yner and Lie 2016). The rock is slightly compressible and the fluid phases are assumed to have constant compressibility so that
\begin{equation} 
b_l = b_{l}(p_0) \, e^{\left ( p - p_0 \right ) c_l} ,
\end{equation}
where $p_0$ is a reference pressure and $c_l$ the compressibility factor for the phase $l$.

We consider simple relative permeabilities for the three phases
\begin{equation} 
k_{rw}\left ( s_w \right ) = s_{w}^{n}, \quad k_{rg}\left ( s_g \right ) = s_{g}^{n}, \quad k_{row}\left ( s_o \right ) = k_{rog}\left ( s_o \right ) = s_{o}^{n} ,
\end{equation}
where $n = 2$ or $3$ in this paper. For oil relative permeability, we employ the weighted interpolation model by Baker (1988), which is the default setting of the commercial simulator Eclipse (Schlumberger 2013). This choice ensures positive values and continuous derivatives provided that $k_{row}$, $k_{rog}$ fulfill the same criteria (Baker 1988; Lee and Efendiev 2016).

The hybrid upwinding scheme for numerical flux and the QN method for nonlinear acceleration are applied to ensure the convergence of outer loop. For QN $m = 3 $ and $\omega = 0.5 $. The convergence tolerances of the inner solvers for the different methods are summarized in Table \ref{tab:conv_tol_1}. 

\begin{table}[!htb]
\centering
\caption{Convergence parameters of inner solvers for different solution strategies}
\label{tab:conv_tol_1}
\begin{tabular}{|c|c|c|c|c|c|c|c|}
\hline
\multicolumn{2}{|c|}{Tight tolerances} & \multicolumn{2}{|c|}{Absolute relaxation} &
\multicolumn{2}{|c|}{Relative relaxation} &
\multicolumn{2}{|c|}{Adaptive relaxation}
\\ \hline
$\epsilon_p$ & $\epsilon_t$ & $\epsilon_p$ & $\epsilon_t$ & $\xi_p$ & $\xi_s$ & $\beta_p$ & $\beta_s$           \\ \hline
$10^{-7}$ & $10^{-5}$ & 1.0 & 0.1 & 0.1 & 0.1 & 0.5 & 0.5
\\ \hline
\end{tabular}
\end{table}

\subsection{Case 1: gravity-driven two-phase flow}

We consider a gravity-driven two-phase flow scenario. The specification of the base model is shown in Table \ref{tab:specification}. Quadratic relative permeability functions are used. The simulation control parameters are summarized in Table \ref{tab:control}.

\begin{table}[!htb]
\centering
\caption{Case 1. Specification of the base model}
\label{tab:specification}
\begin{tabular}{|c|c|c|}
\hline
Parameter                  &  Value           & Unit      \\ \hline
NX / NZ                    &  60 / 60         &           \\ \hline
LX / LY / LZ               &  600 / 10 / 600  & ft        \\ \hline
Initial pressure           &  2000            & psi       \\ \hline
Rock permeability          &  100             & mD        \\ \hline
Rock porosity              &  0.1             &           \\ \hline
Water density              &  1000            & $\textrm{kg}/\textrm{m}^3$  \\ \hline
Oil   density              &  500             & $\textrm{kg}/\textrm{m}^3$  \\ \hline
Water viscosity            &  1               & cP        \\ \hline
Oil   viscosity            &  4               & cP        \\ \hline
Oil   compressibility      &  $6.9\cdot 10^{-6}$        & 1/psi     \\ \hline
Rock  compressibility      &  $10^{-6}$            & 1/psi     \\ \hline
\end{tabular}
\end{table}

\begin{table}[!htb]
\centering
\caption{Case 1. Simulation control parameters}
\label{tab:control}
\begin{tabular}{|c|c|c|}
\hline
Parameter                 & Value          & Unit   \\ \hline
Maximum time-step size    & 50             & day    \\ \hline
Total simulation time     & 400            & day    \\ \hline
Maximum number of outer iterations      & 30         &     \\ \hline
Convergence tolerance of outer loop     & 0.001      &     \\ \hline
\end{tabular}
\end{table}

\subsubsection{Case 1a: homogeneous lock-exchange problem}

\begin{figure}[!htb]
\centering
\includegraphics[scale=0.5]{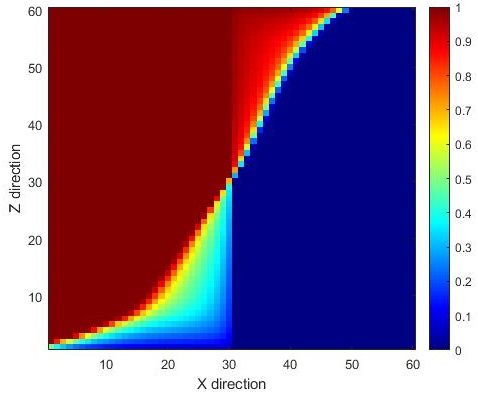}
\caption{Oil saturation profile for Case 1a: homogeneous lock-exchange problem.}
\label{fig:s_h_le_pg}
\end{figure}

We first tested a lock-exchange problem on a homogeneous square model. Oil initially occupies the left half of the domain, while water fills the right half. The lock-exchange problem is challenging for the SFI method, because gravity significantly contributes to the total velocity and induces a global re-circulation flow pattern. The oil saturation profile at the end of simulation is plotted in \textbf{Fig.~\ref{fig:s_h_le_pg}}.

\subsubsection{Case 1b: heterogeneous lock-exchange problem}

\begin{figure}[!htb]
\centering
\subfloat[Permeability]{
\includegraphics[scale=0.5]{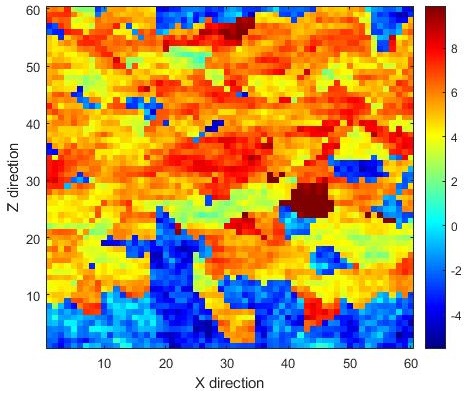}}
\subfloat[Porosity]{
\includegraphics[scale=0.5]{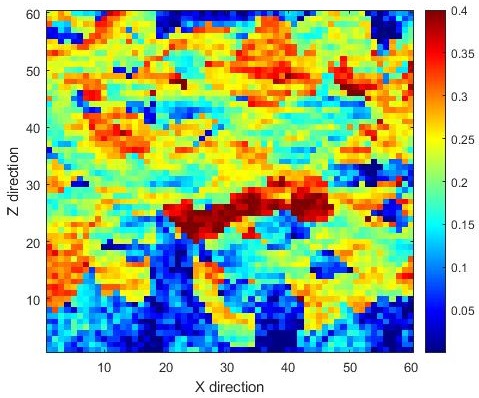}}
\caption{Random permeability (log) and porosity fields for the heterogeneous square model.}
\label{fig:perm_poro_sq}
\end{figure}

The lock-exchange problem was also tested on a heterogeneous square model.~The random rock properties of the model are shown in \textbf{Fig.~\ref{fig:perm_poro_sq}}. Oil density is changed to 800~$\textrm{kg}/\textrm{m}^3$ and maximum time-step size becomes 20~days, with total simulation time as 200~days for the following two-phase cases. The simulation parameters are modified because the case is very difficult for the basic SFI method to converge.

\begin{figure}[!htb]
\centering
\subfloat[Tight tolerances]{
\includegraphics[scale=0.4]{FIGURE/residual/2p/2p_no.pdf}}
\\
\subfloat[Absolute relaxation]{
\includegraphics[scale=0.4]{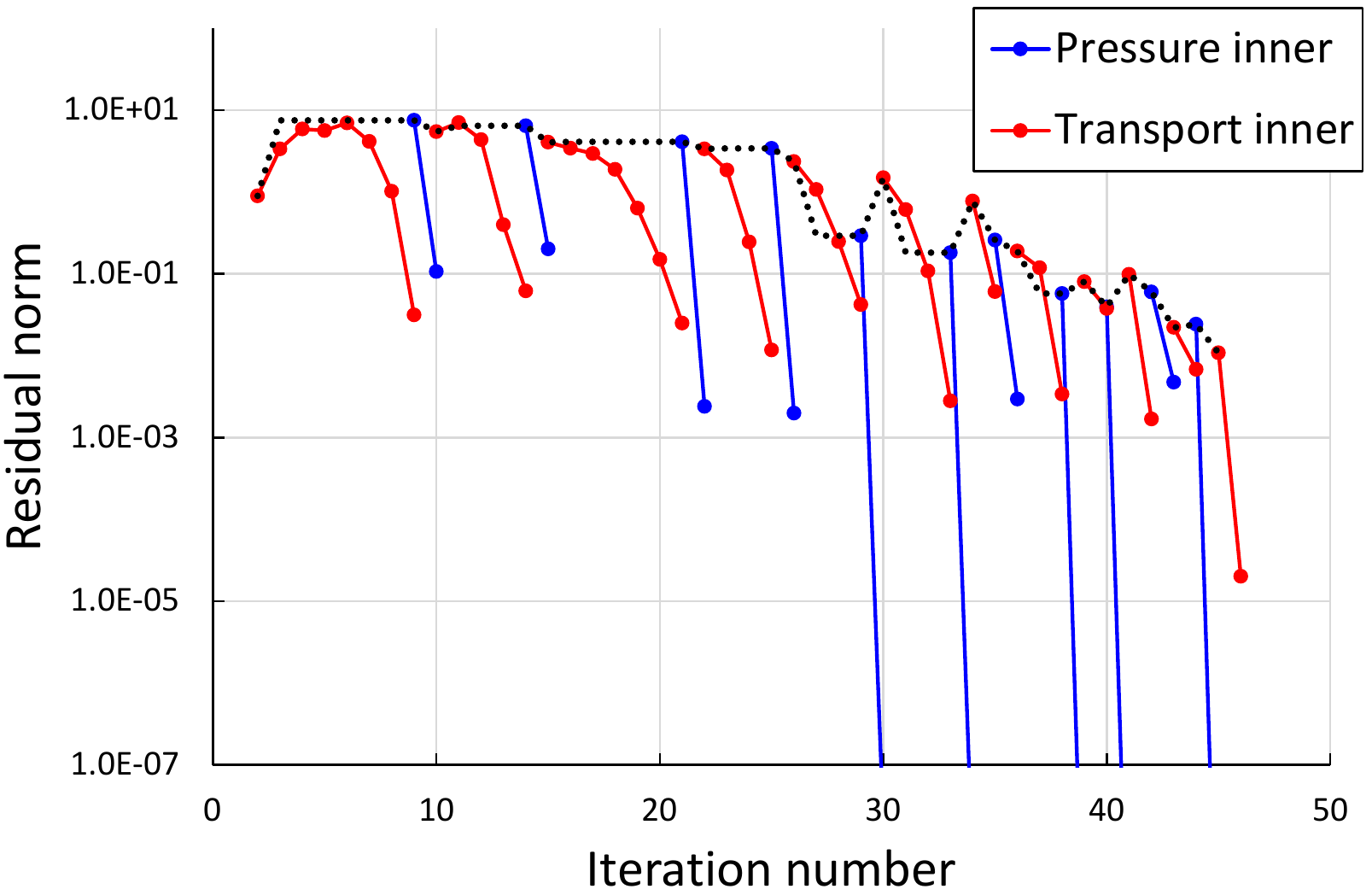}}
\\
\subfloat[Adaptive relative relaxation]{
\includegraphics[scale=0.4]{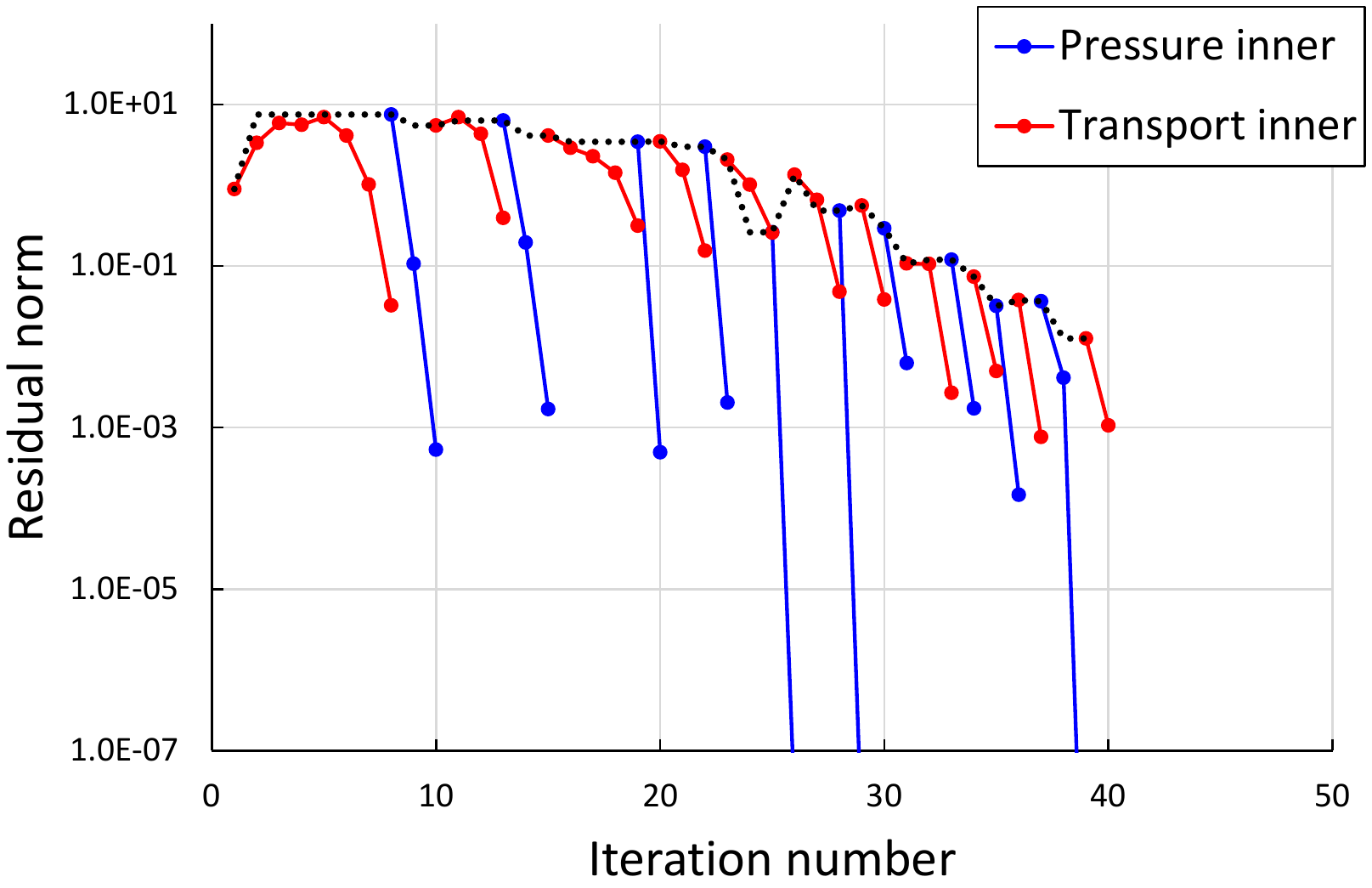}}
\caption{Residual history of the inner solvers for Case 1b: heterogeneous lock-exchange problem.}
\label{fig:resi_1b}
\end{figure}

\textbf{Fig.~\ref{fig:resi_1b}} shows the residual history of the inner solvers for the second timestep. As~discussed above, the outer convergence is mainly restricted by the initial residuals of the sub-problems. The residual rebound of one inner solver is due to the feedback from the other solver. In standard SFI (tight tolerances), the two sub-problems are solved to high precision at every outer iteration. However, additional accuracy achieved in inner solvers is wasted, contributing to little reduction of the overall residual. By comparison, the inexact methods apply relaxed tolerances and thus alleviate the over-solving issue. This enables reaching the same coupled solutions with fewer inner iterations.

\begin{figure}[!htb]
\centering
\includegraphics[scale=0.6]{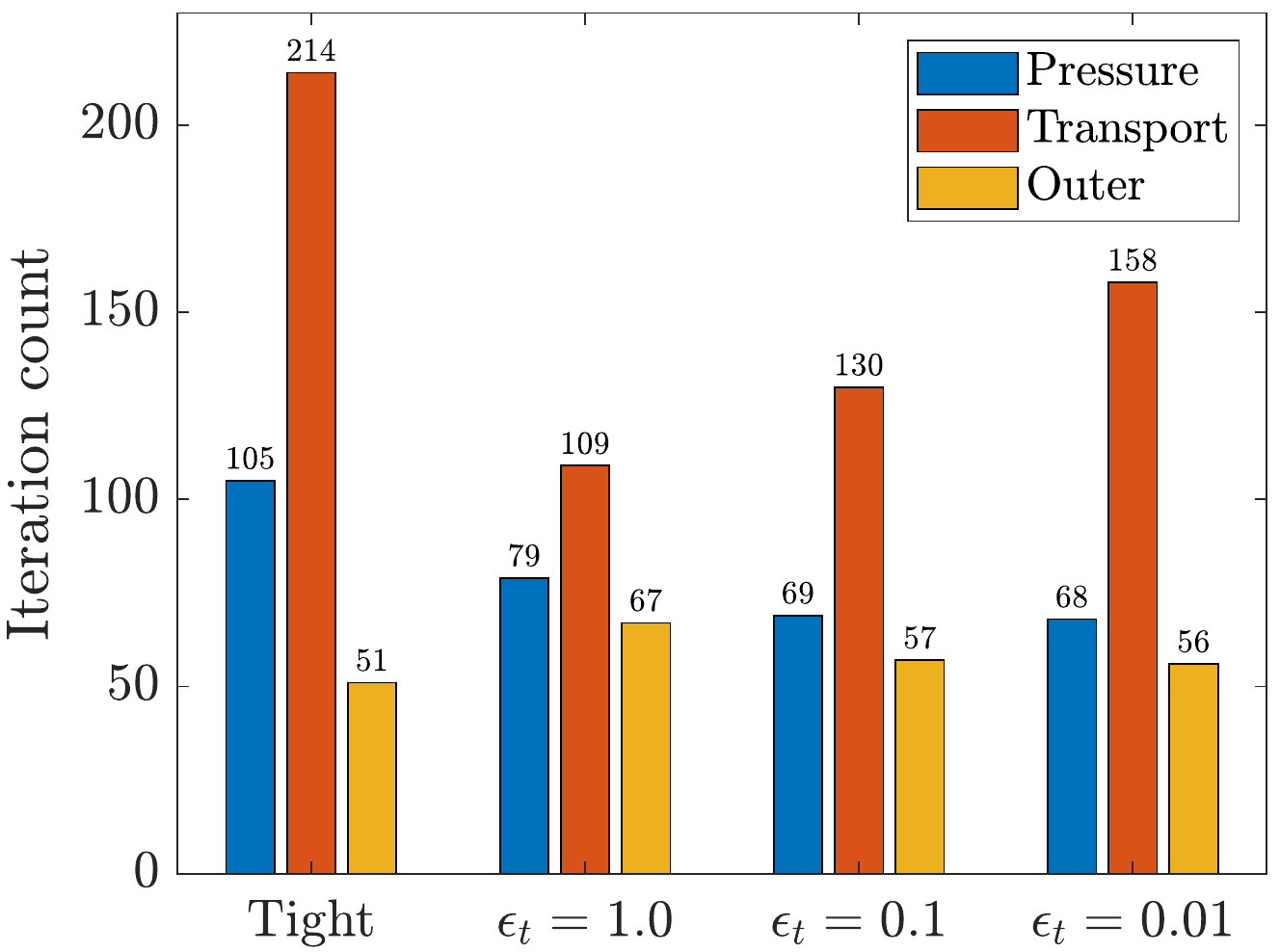}
\caption{Iteration performance of the absolute relaxation strategy applied to the transport solver for Case 1b: heterogeneous lock-exchange problem.}
\label{fig:iter_1b}
\end{figure}

We tested the fixed absolute relaxations for the transport solver. The iteration performance of the cases with different values of $\epsilon_t$ is summarized in \textbf{Fig.~\ref{fig:iter_1b}}. From the results we can see that the overall performance is sensitive to the relaxation level. A~tight tolerance leads to large number of transport iterations. On the other hand, a much relaxed tolerance worsens the performance of outer loop as well as pressure solver, even though the transport cost is reduced.


\subsubsection{Case 1c: heterogeneous gravity segregation problem}

\begin{figure}[!htb]
\centering
\includegraphics[scale=0.5]{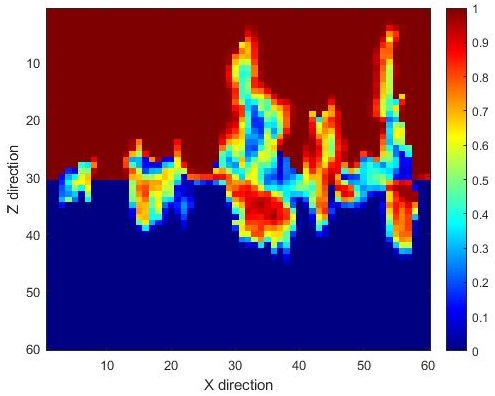}
\caption{Oil saturation profile of Case 1c: heterogeneous gravity segregation problem.}
\label{fig:s_he_se_pg}
\end{figure}

We also tested a counter-current flow problem. In this case, a complex fluid dynamics quickly develops due to gravity segregation during the simulation. The oil saturation profile is plotted in \textbf{Fig.~\ref{fig:s_he_se_pg}}.

\subsubsection{Summary of iteration performance}

\begin{figure}[!htb]
\centering
\subfloat[Case 1a]{
\includegraphics[scale=0.5]{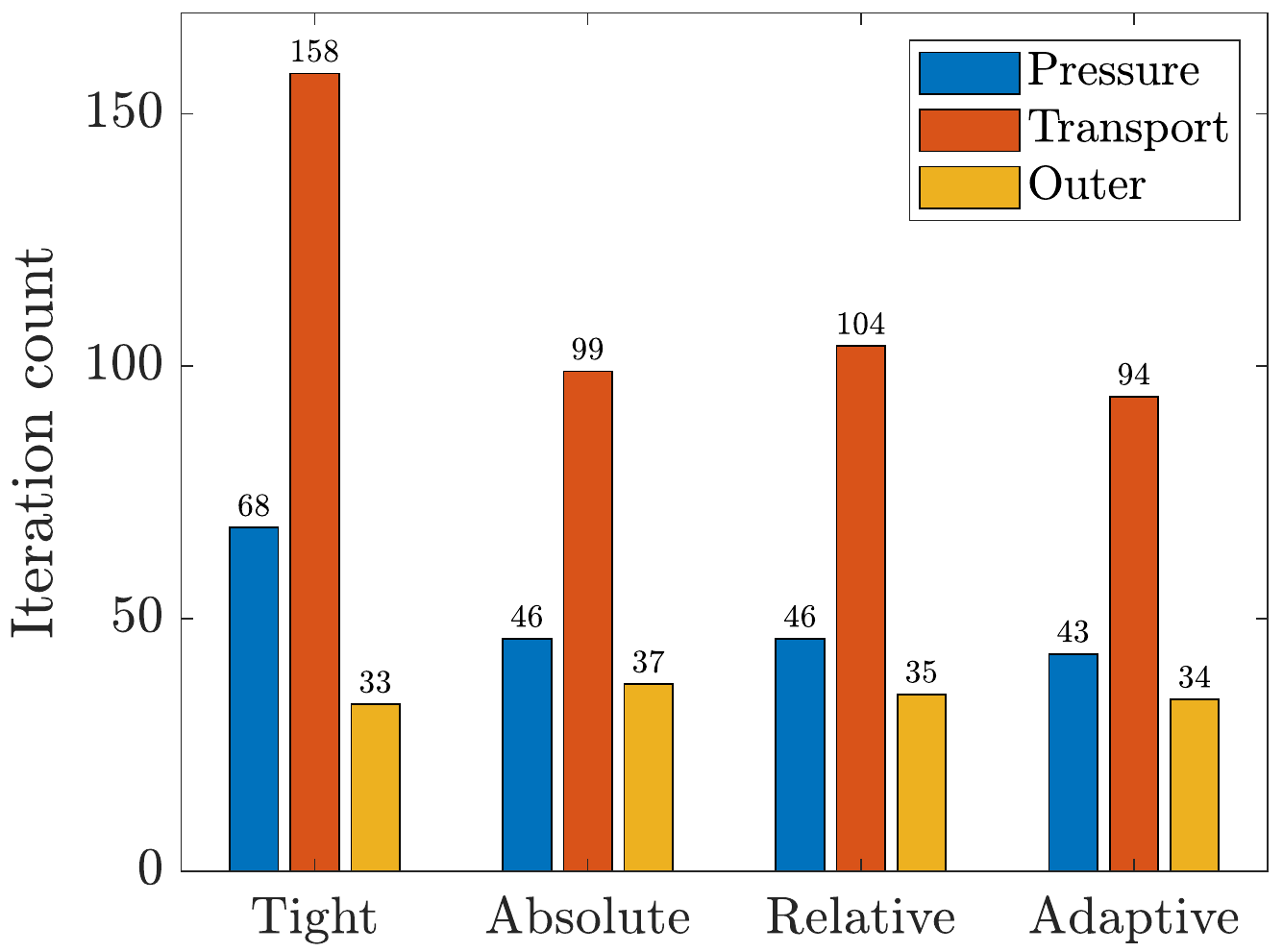}}
\\
\subfloat[Case 1b]{
\includegraphics[scale=0.5]{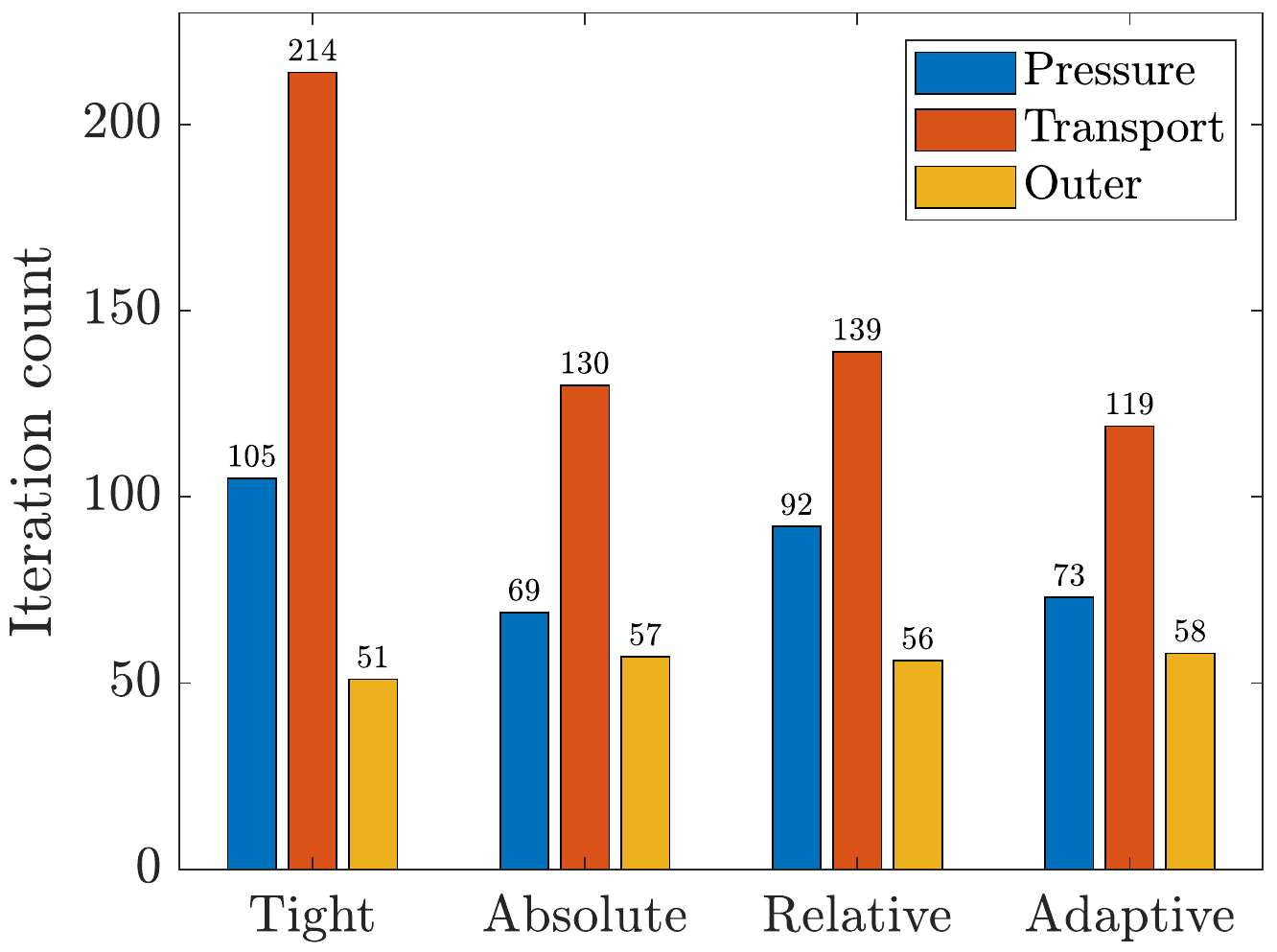}}
\\
\subfloat[Case 1c]{
\includegraphics[scale=0.5]{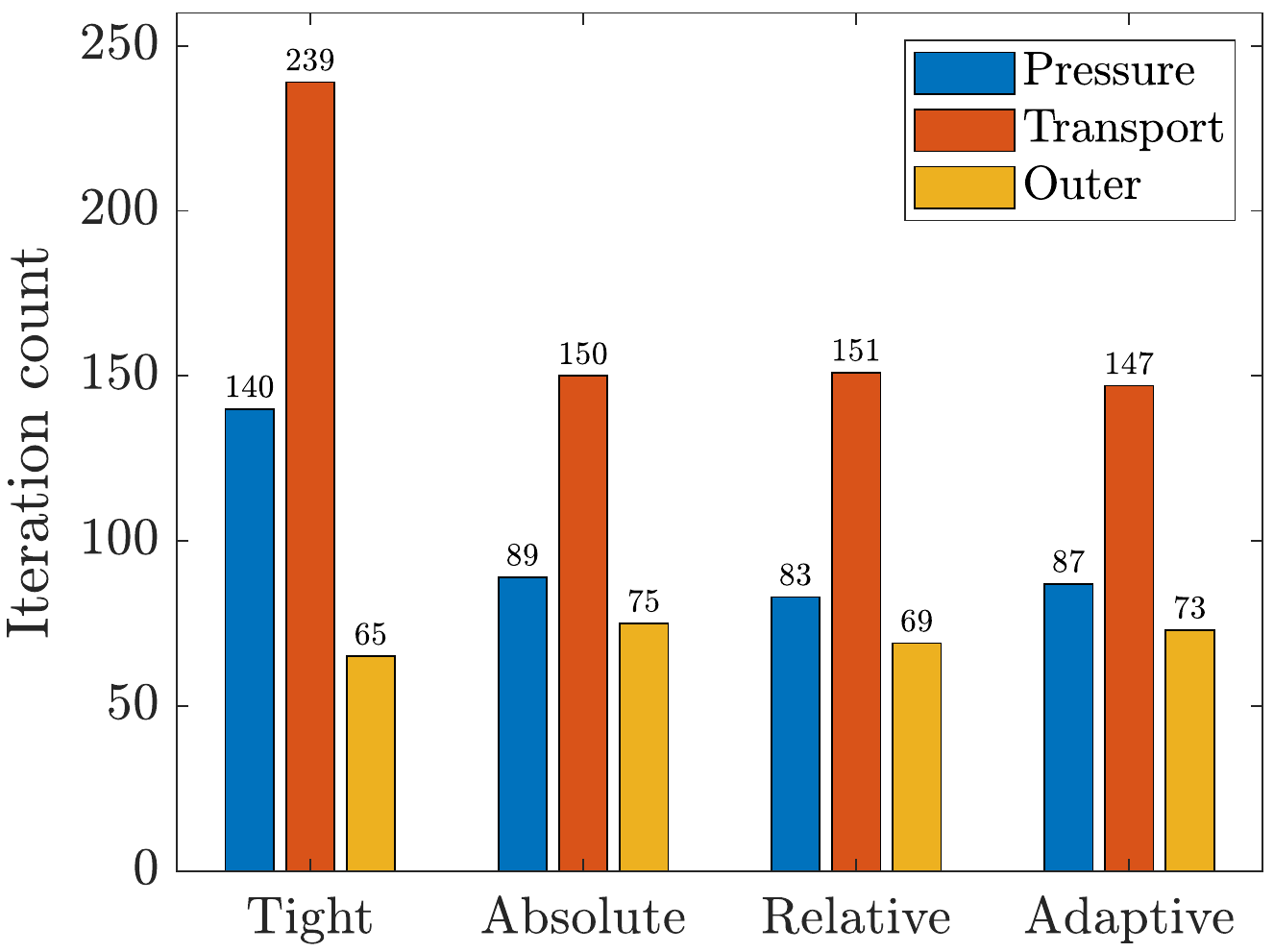}}
\caption{Case 1 iteration performance.}
\label{fig:iter_1}
\end{figure}

The iteration performance for the three cases is summarized in \textbf{Fig.~\ref{fig:iter_1}}. The results demonstrate that the inexact methods greatly improve the overall efficiency by reducing the number of inner iterations, while preserving the outer-loop convergence.~Among the three inexact methods, the adaptive strategy which provides the relative tolerances adequate for the sub-problems enjoys the least iteration counts. For this particular set of cases, the absolute relaxation strategy also shows a good performance; however, as will be demonstrated for other problems below, this is not generally the case.

\subsection{Case 2: two-phase gravity segregation, SPE 10 model}

\begin{figure}[!htb]
\centering
\subfloat[Permeability]{
\includegraphics[scale=0.5]{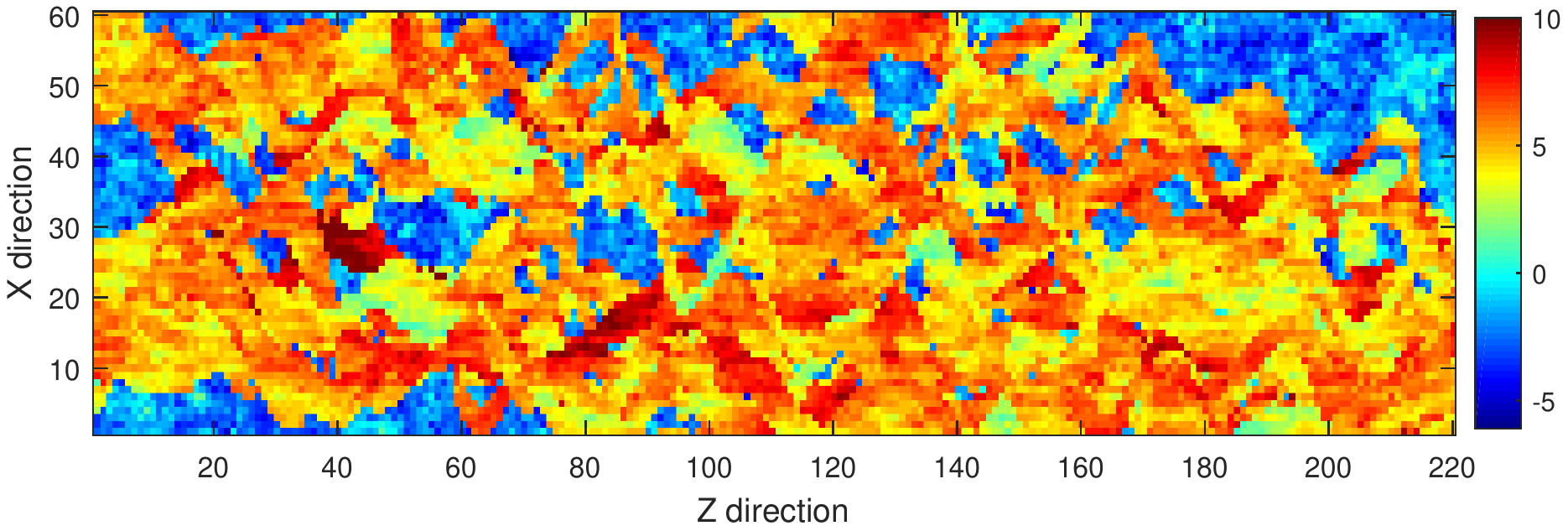}}
\\
\subfloat[Porosity]{
\includegraphics[scale=0.5]{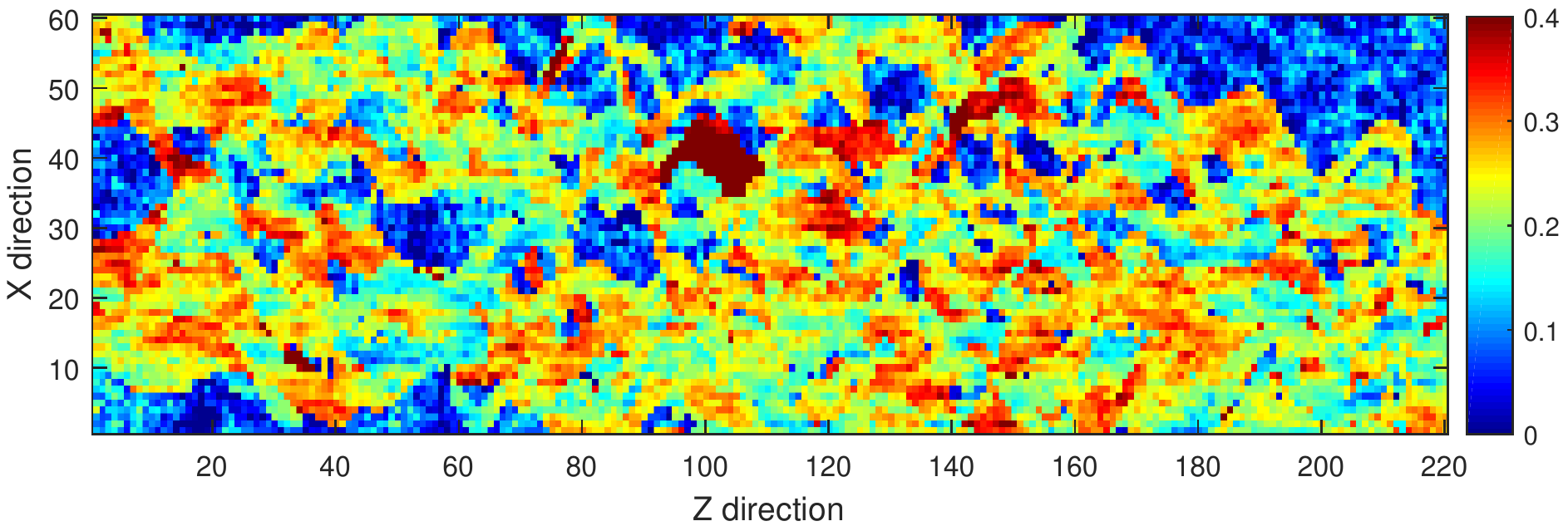}}
\caption{Permeability (log) and porosity fields of the bottom layer of the SPE 10 model.}
\label{fig:perm_poro}
\end{figure}
\begin{figure}[!htb]
\centering
\includegraphics[scale=0.5]{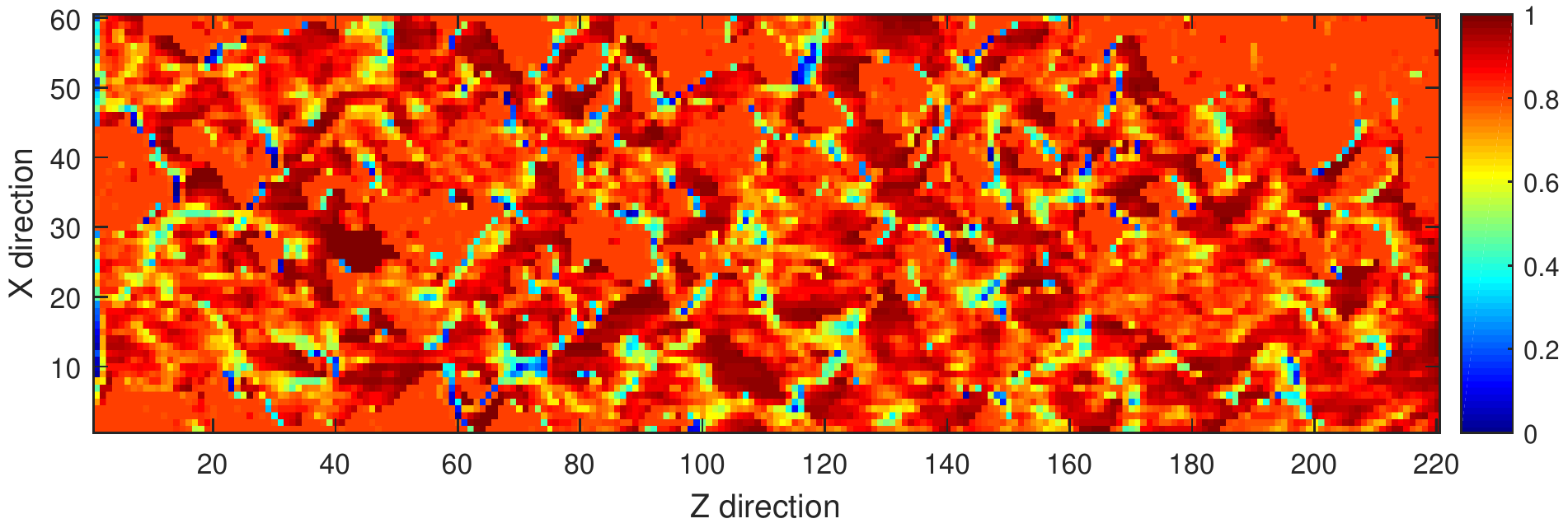}
\caption{Oil saturation profile for Case 2: two-phase gravity segregation, SPE 10 model.}
\label{fig:s_ne_spe_pg}
\end{figure}
\begin{figure}[!htb]
\centering
\includegraphics[scale=0.6]{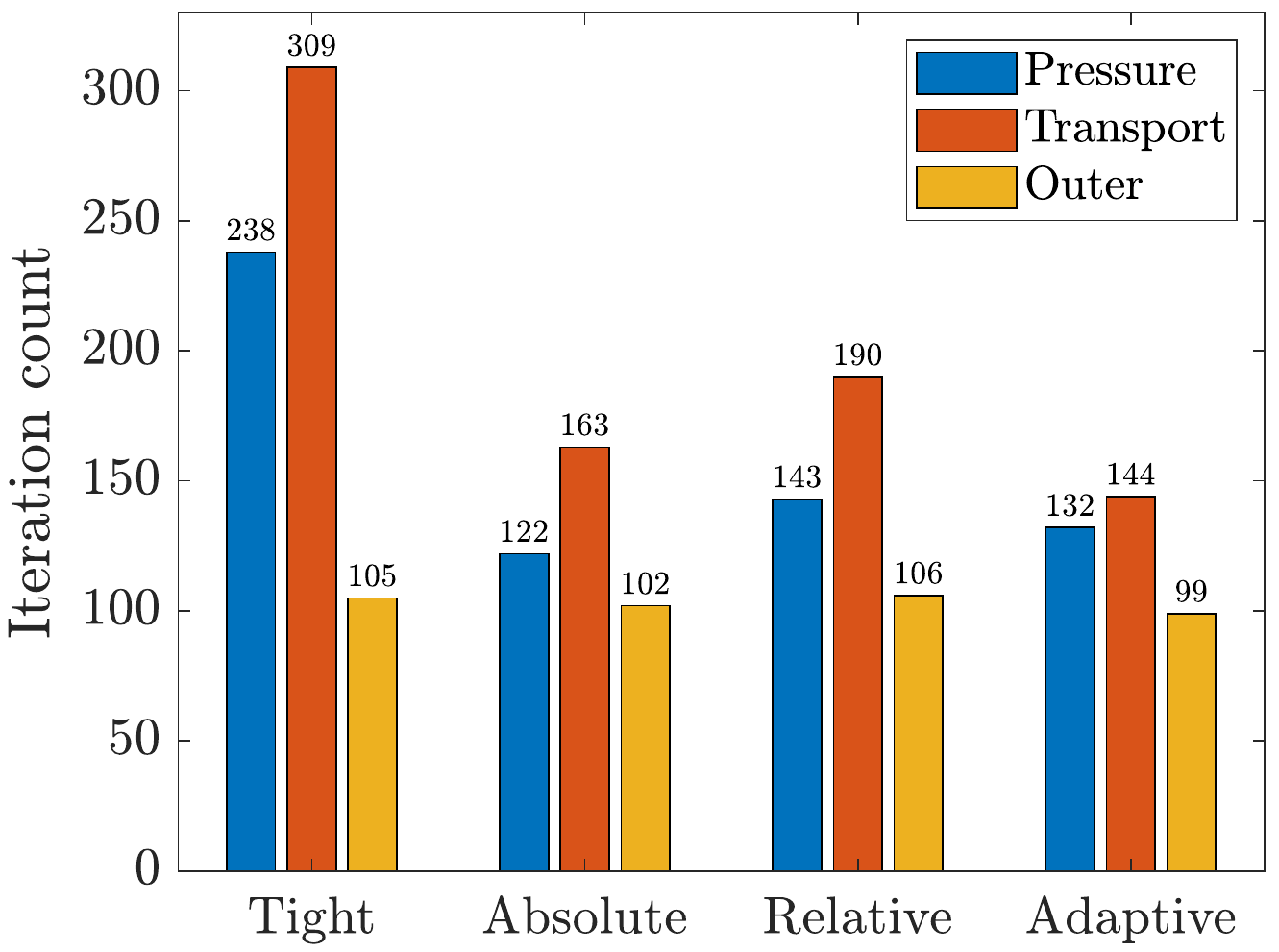}
\caption{Iteration performance of Case 2: two-phase gravity segregation, SPE 10 model.}
\label{fig:iter_2}
\end{figure}

The model setting is identical to the previous case as specified in Table~\ref{tab:specification}. The bottom layer of the SPE~10 model (Christie and Blunt 2001) is used. The rock properties are shown in \textbf{Fig.~\ref{fig:perm_poro}}. Maximum time-step size is 10~days, with the total simulation time of 200~days. Initial water saturation is uniformly set to 0.2. Due to the non-equilibrium initial condition, gravity segregation starts to take place in all cells from the beginning. The oil saturation profile is plotted in \textbf{Fig.~\ref{fig:s_ne_spe_pg}}. The iteration performance of Case 2 is summarized in \textbf{Fig.~\ref{fig:iter_2}}. Compared to the standard method, the inexact methods lead to substantial reductions in the total iterations.

\subsection{Case 3: two-phase flow with viscous and gravitational forces, SPE 10 model}

We consider a scenario with combined viscous and gravitational forces. A quarter-five spot well pattern is applied: water is injected at the middle of the domain and producers are placed at the four corners. The injection rate is $2.5\times 10^{-4}$ pore volume (664 $\textrm{ft}^3$) per day. The modified model parameters of Case 3 are summarized in Table~\ref{tab:specification_VG}. The other parameters specified in Case 1a remain unchanged.

\begin{table}[!htb]
\centering
\caption{Modified model parameters of Case 3}
\label{tab:specification_VG}
\begin{tabular}{|c|c|c|}
\hline
Parameter                  &  Value            & Unit      \\ \hline
NX / NZ                    &  60 / 220         &           \\ \hline
LX / LY / LZ               &  600 / 10 / 2200  & ft        \\ \hline
Initial water saturation   &  0.01             &           \\ \hline
Injection rate             &  664              & $\textrm{ft}^3/\textrm{day}$  \\ \hline
Production BHP             &  500              & psi          \\ \hline
Maximum time-step size     &  20            & day    \\ \hline
Total simulation time      &  200           & day    \\ \hline


\end{tabular}
\end{table}

\begin{figure}[!htb]
\centering
\includegraphics[scale=0.5]{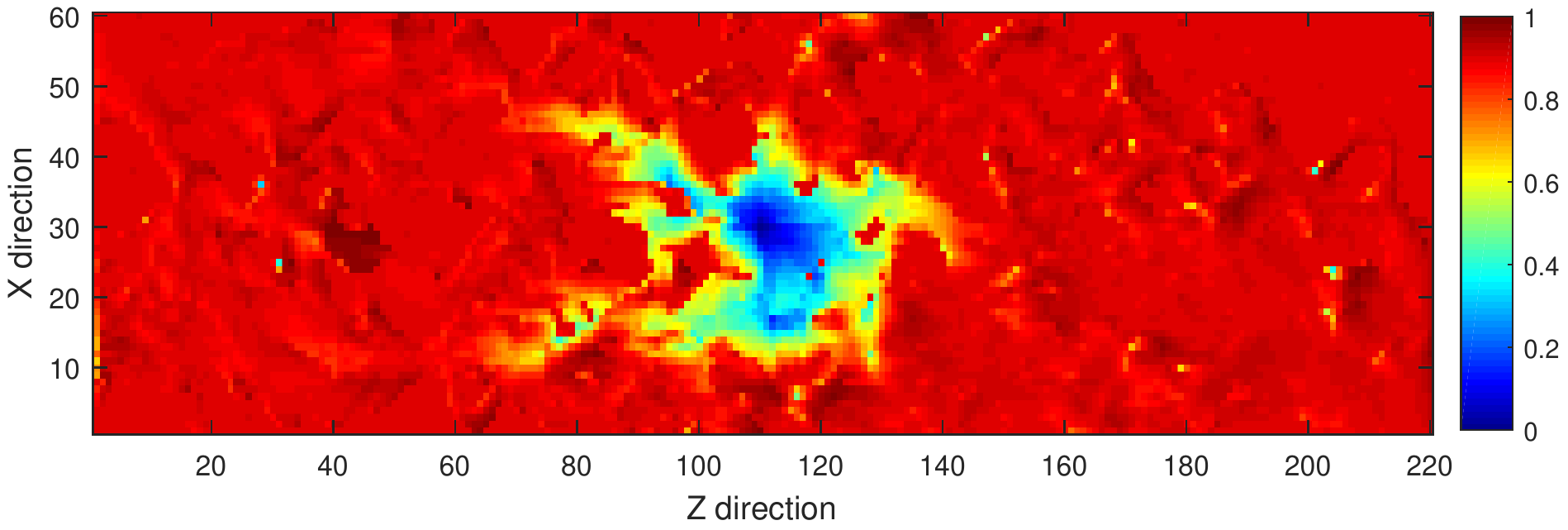}
\caption{Oil saturation profile of Case 3.}
\label{fig:s_VG_2_spe}
\end{figure}  

\begin{figure}[!htb]
\centering
\includegraphics[scale=0.6]{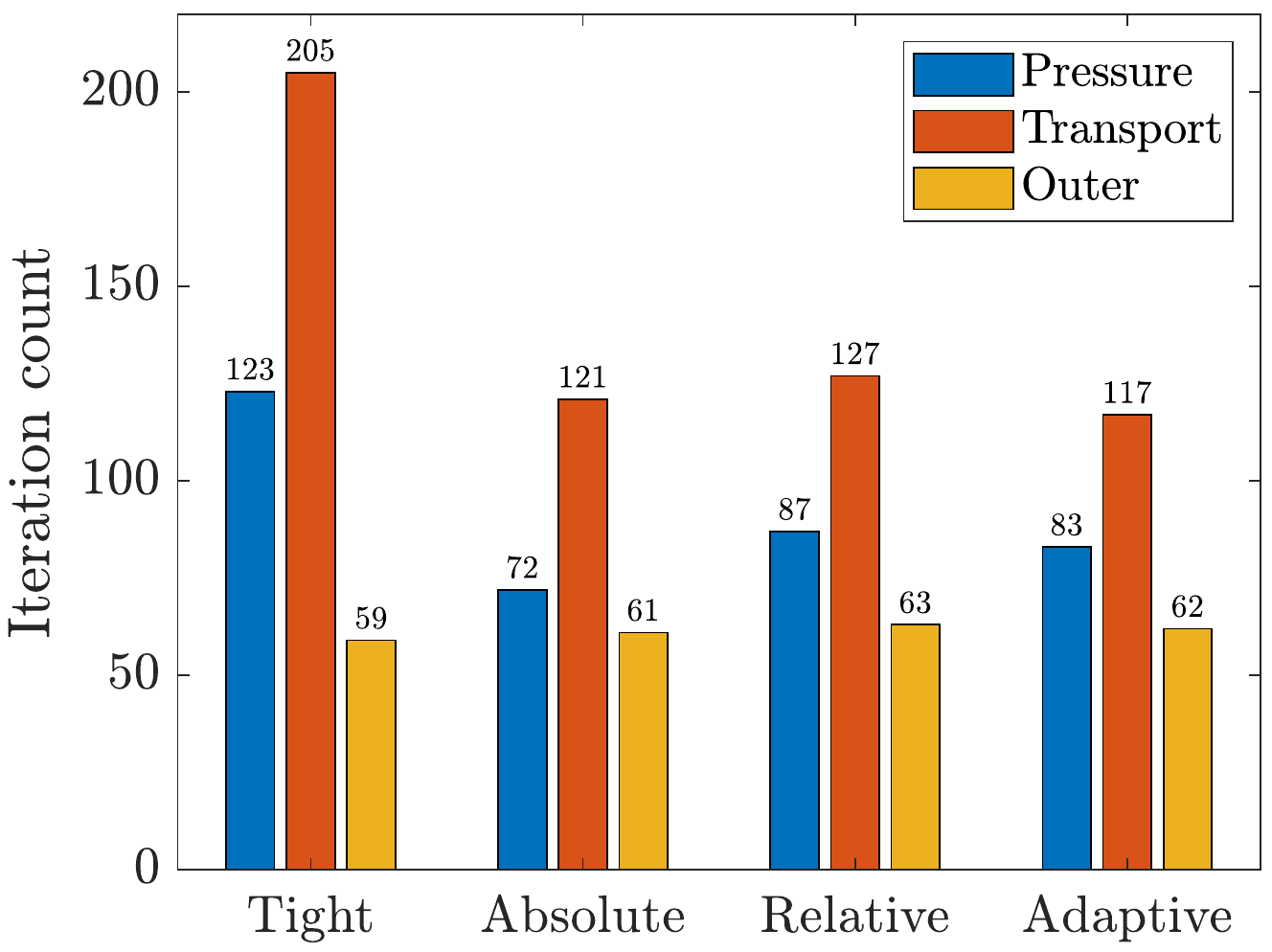}
\caption{Case 3 iteration performance.}
\label{fig:iter_3}
\end{figure}

The oil saturation profile is plotted in \textbf{Fig.~\ref{fig:s_VG_2_spe}}. The iteration performance of Case 3 is summarized in \textbf{Fig.~\ref{fig:iter_3}}. We can see that the adaptive strategy improves the overall efficiency, though the outer iteration numbers slightly increase. One should always keep in mind that most of simulation cost comes from inner solvers, thus reducing their iteration numbers is the primary goal.

\subsection{Case 4: pure viscous flow, SPE 10 model}

\subsubsection{Case 4a: water injection into oil reservoir}

We consider a scenario with pure viscous force and tested a model with water injecting into an oil reservoir. The model parameters of Case 3 are still used. The injection rate is changed to 1328 $\textrm{ft}^3/\textrm{day}$.

\subsubsection{Case 4b: water injection into gas reservoir}

\begin{figure}[!htb]
\centering
\includegraphics[scale=0.59]{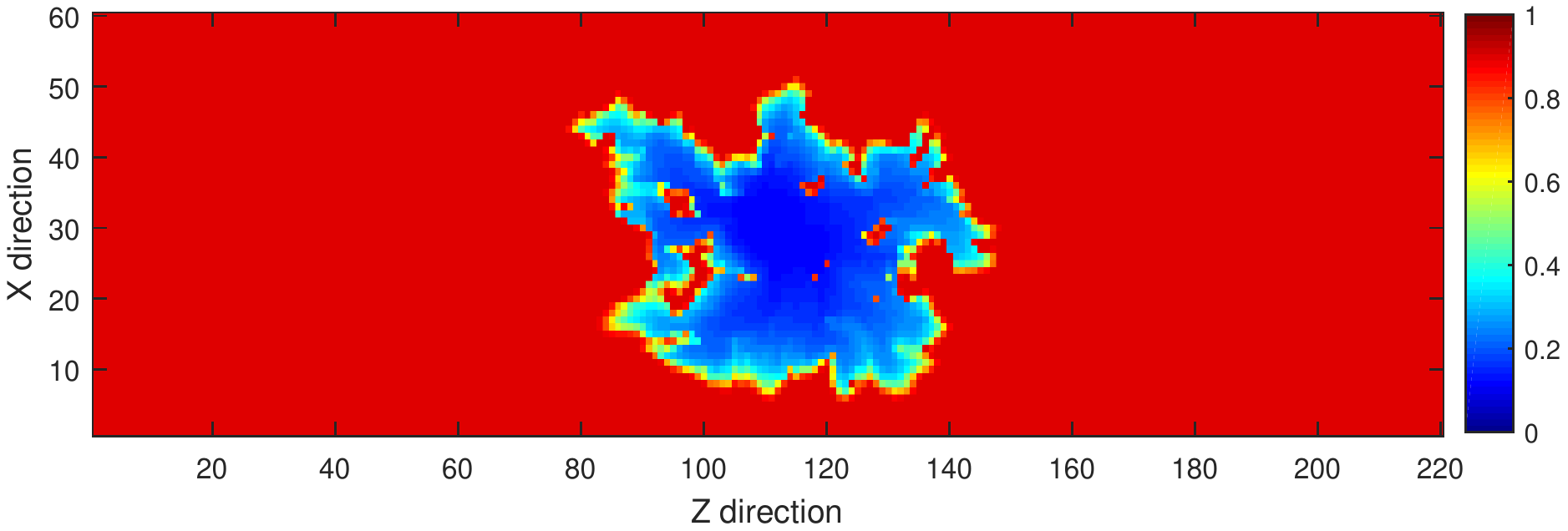}
\caption{Gas saturation profile of Case 4b: water injection into gas reservoir.}
\label{fig:s_WiG_spe}
\end{figure}

We tested a model with water injecting into a gas reservoir. The specification of the model is shown in Table \ref{tab:specification_WiG}. Cubic relative-permeability functions are used, and the viscosity ratio is $M = \mu_w / \mu_g = 4$. The gas saturation profile is plotted in \textbf{Fig.~\ref{fig:s_WiG_spe}}. From the plot we observe that the saturation front is sharp, due to the property of the fractional flow function. The propagation of the sharp front produces large mobility changes and results in a tight coupling between the sub-problems.

\begin{table}[!htb]
\centering
\caption{Specification of the model with water injecting into a gas reservoir}
\label{tab:specification_WiG}
\begin{tabular}{|c|c|c|}
\hline
Parameter                  &  Value            & Unit      \\ \hline
Initial water saturation   &  0.1              &           \\ \hline
Initial pressure           &  2000             & psi       \\ \hline
Water viscosity            &  1                & cP        \\ \hline
Gas   viscosity            &  0.25             & cP        \\ \hline
Gas   compressibility      &  $6.9\cdot 10^{-5}$         & 1/psi     \\ \hline
Rock  compressibility      &  $10^{-6}$             & 1/psi     \\ \hline
Injection rate             &  1328             & $\textrm{ft}^3/\textrm{day}$  \\ \hline
Production BHP             &  500              & psi          \\ \hline
Maximum time-step size     &  20               & day    \\ \hline
Total simulation time      &  200              & day    \\ \hline
\end{tabular}
\end{table}

\begin{figure}[!htb]
\centering
\subfloat[Case 4a]{
\includegraphics[scale=0.46]{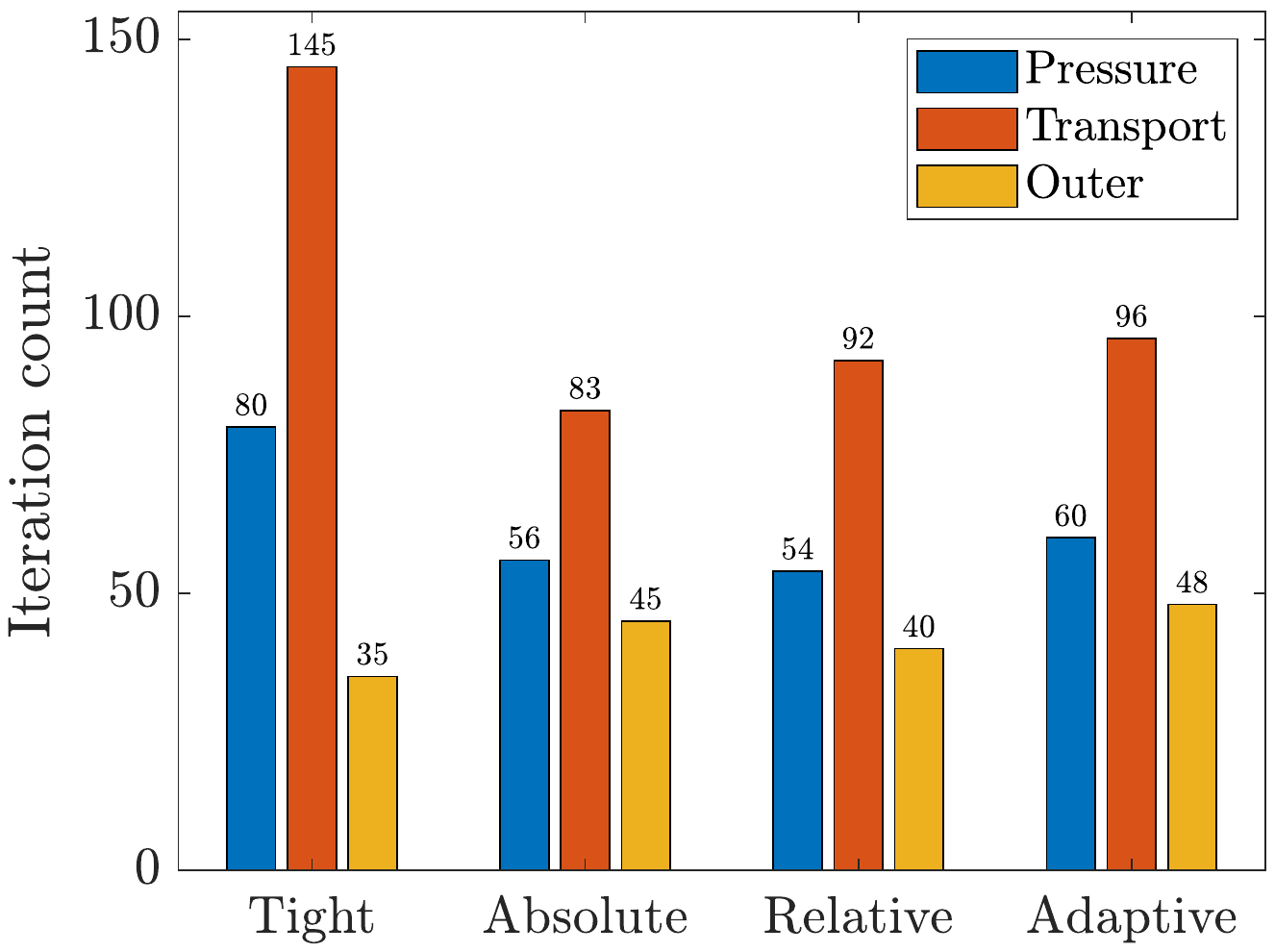}}
\
\subfloat[Case 4b]{
\includegraphics[scale=0.46]{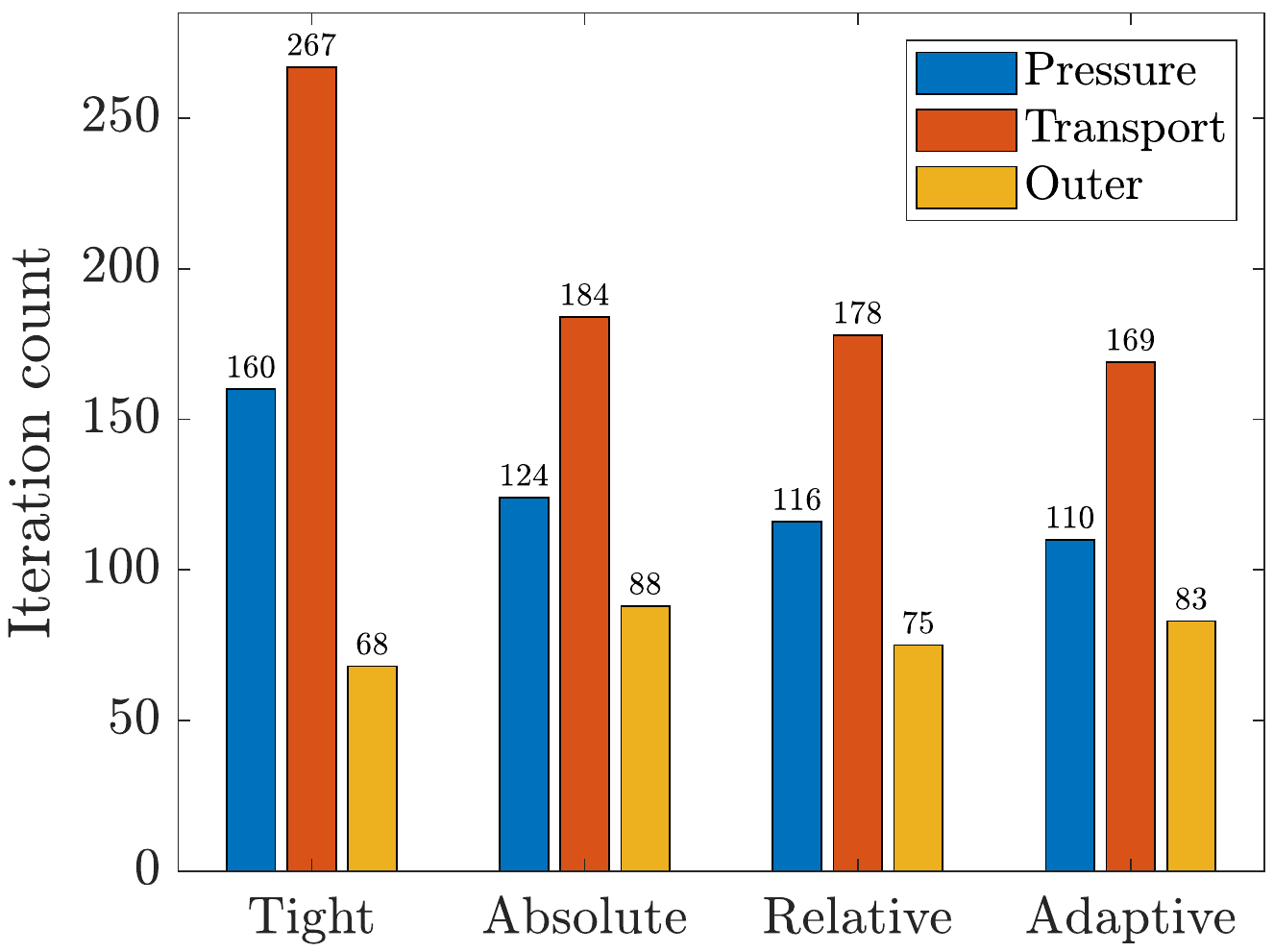}}
\caption{Iteration performance of Case 4: pure viscous flow, SPE 10 model.}
\label{fig:iter_4}
\end{figure}

The iteration performance for Case 4 is summarized in \textbf{Fig.~\ref{fig:iter_4}}. The inexact methods exhibit better performances than the standard method with tight tolerances.

\subsection{Case 5: three-phase pure viscous flow, SPE 10 model}

The specification of the three-phase model is shown in Table~\ref{tab:specification_WAG}. In the presence of the three phases, a challenging scenario is generated with strong coupling terms between the flow and transport.

\subsubsection{Case 5a: water-alternating-gas injection}

\begin{figure}[!htb]
\centering
\subfloat[Oil saturation profile]{
\includegraphics[scale=0.5]{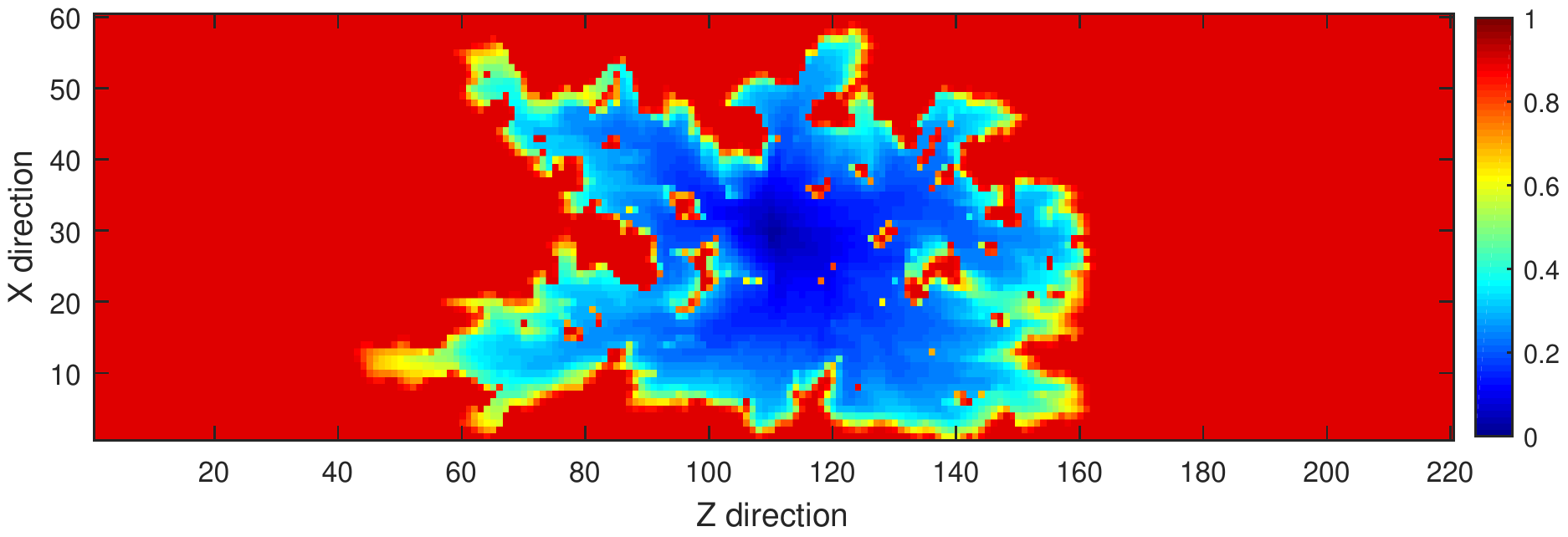}}
\\
\subfloat[Gas saturation profile]{
\includegraphics[scale=0.5]{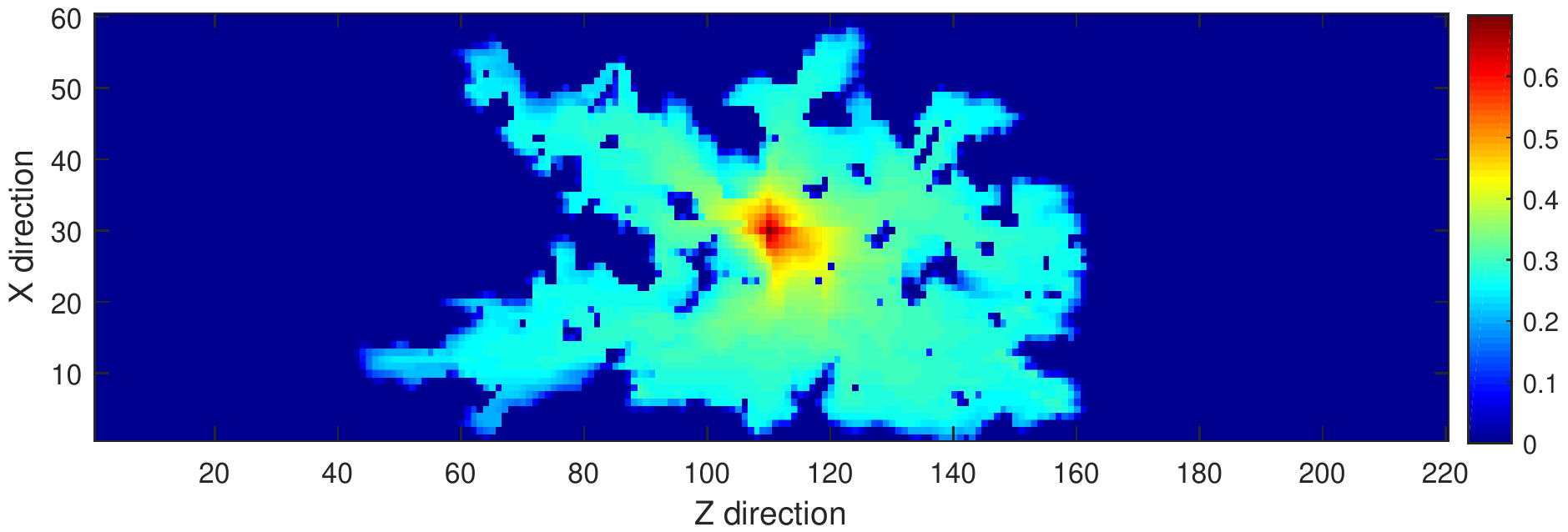}}
\caption{Saturation profiles of Case 5a: water-alternating-gas injection.}
\label{fig:s_WAG_spe}
\end{figure}

We first consider a case with water-alternating-gas (WAG) injection. Water is injected first, then followed by gas injection, and the process alternates for each time interval. Cubic relative-permeability functions are used. The oil and gas saturation profiles are plotted in \textbf{Fig.~\ref{fig:s_WAG_spe}}. The simulation scenario is quite challenging because of the frequently changing well schedule.

\begin{table}[!htb]
\centering
\caption{Specification of the three-phase model}
\label{tab:specification_WAG}
\begin{tabular}{|c|c|c|}
\hline
Parameter                  &  Value            & Unit      \\ \hline
Initial water saturation   &  0.1              &           \\ \hline
Initial oil saturation     &  0.9              &           \\ \hline
Initial pressure           &  2000             & psi       \\ \hline
Water viscosity            &  1                & cP        \\ \hline
Oil   viscosity            &  1                & cP        \\ \hline
Gas   viscosity            &  0.25             & cP        \\ \hline
Oil   compressibility      &  $6.9\cdot 10^{-6}$         & 1/psi     \\ \hline
Gas   compressibility      &  $6.9\cdot 10^{-5}$         & 1/psi     \\ \hline
Rock  compressibility      &  $10^{-6}$             & 1/psi     \\ \hline
Injection rate             &  1328             & $\textrm{ft}^3/\textrm{day}$  \\ \hline
Production BHP             &  500              & psi          \\ \hline
Time-interval size         &  20               & day    \\ \hline
Total simulation time      &  400              & day    \\ \hline
\end{tabular}
\end{table}

\subsubsection{Case 5b: water injection}

We also consider a simpler water injection case. The modified model parameters are summarized in Table~\ref{tab:specification_WI}.

\begin{figure}[!htb]
\centering
\subfloat[Case 5a]{
\includegraphics[scale=0.46]{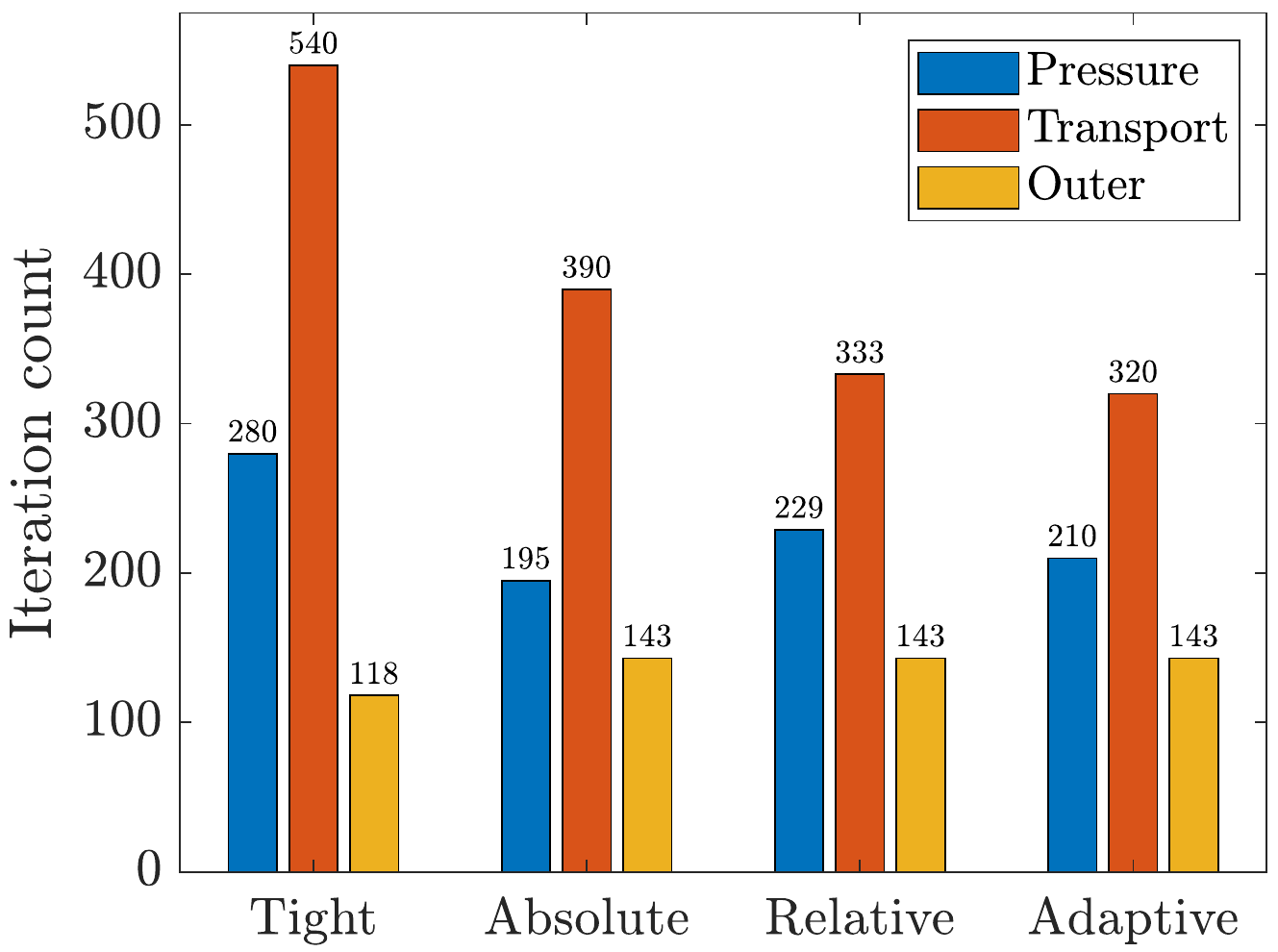}}
\
\subfloat[Case 5b]{
\includegraphics[scale=0.46]{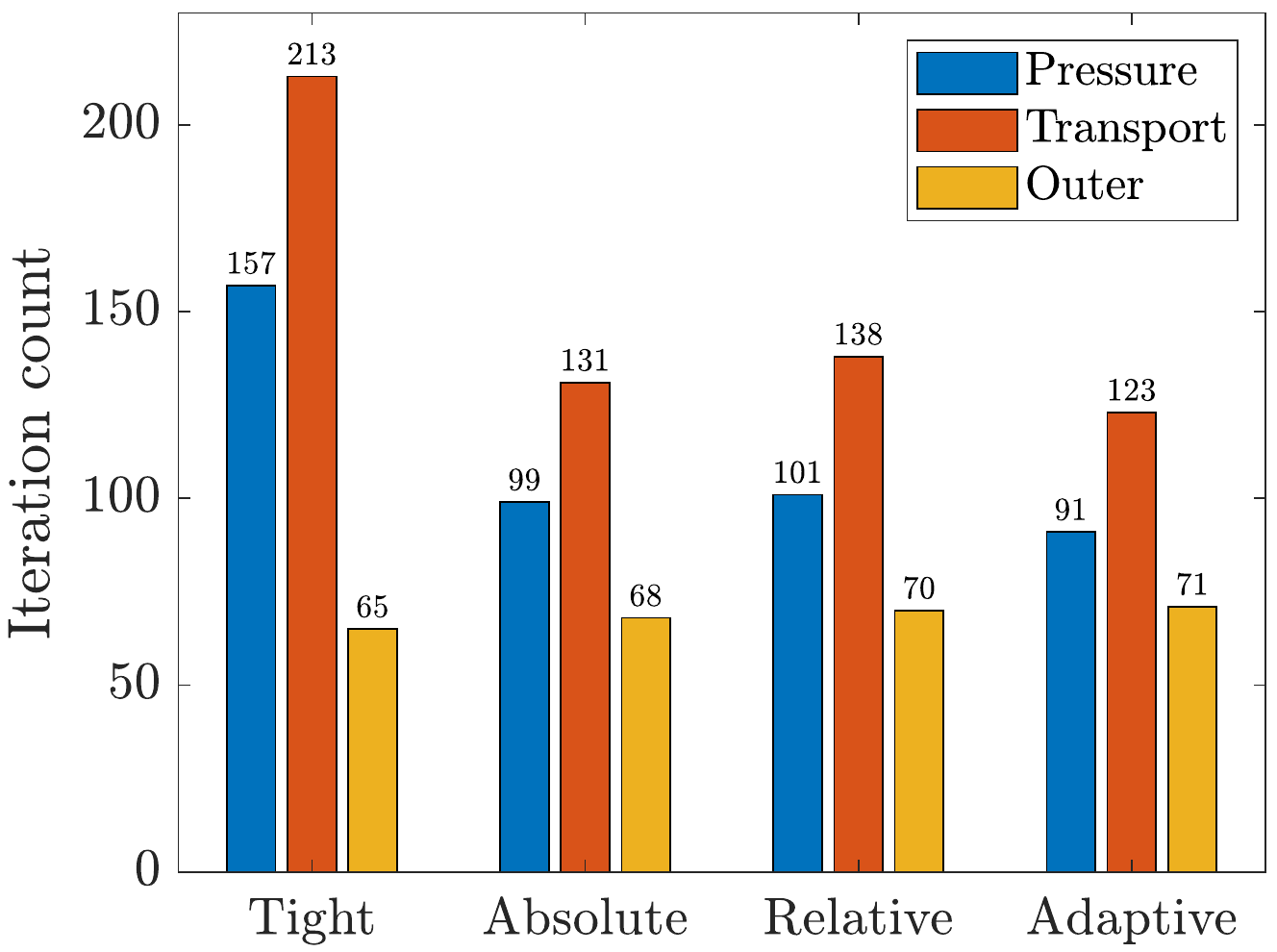}}
\caption{Iteration performance of Case 5: three-phase pure viscous flow, SPE 10 model.}
\label{fig:iter_5}
\end{figure}

The iteration performance for Case 5 is summarized in \textbf{Fig.~\ref{fig:iter_5}}. The results show that the inexact methods lead to more outer iterations, but lower overall cost. In particular, compared with tight tolerances, the adaptive strategy achieves 50\% reductions in the~transport iterations.

\begin{table}[!htb]
\centering
\caption{Modified model parameters of Case 5b}
\label{tab:specification_WI}
\begin{tabular}{|c|c|c|}
\hline
Parameter                  &  Value            & Unit      \\ \hline
Initial water saturation   &  0.1              &           \\ \hline
Initial oil saturation     &  0.2              &           \\ \hline
Oil   viscosity            &  4                & cP        \\ \hline
Injection rate             &  2656             & $\textrm{ft}^3/\textrm{day}$  \\ \hline
Total simulation time      &  200              & day    \\ \hline
\end{tabular}
\end{table}

\section{Results: compositional flow}

The inexact-SFI framework developed in this work was integrated into the Automatic Differentiation General Purpose Research Simulator (AD-GPRS) (Younis and Aziz 2007; Zhou et al. 2011; Voskov and Tchelepi 2012; Rin et al. 2017; Garipov et al. 2018). We validated the effectiveness of the inexact methods using several problems of increasing complexity. Three different fluids, and three different reservoir models are considered. The models include: (a)~homogeneous 1D model; (b)~gravity-driven lock exchange; (c)~bottom layer of the SPE 10 problem. In the following cases, the relative permeabilities are quadratic for both phases, unless otherwise indicated.

The linear systems are solved using the Intel MKL PARDISO solver. Here we only study the iterative performance of the outer and nonlinear inner loops. During nonlinear iterations, the standard Newton method can lead to unphysical values for some variables and must be explicitly corrected. Given the results from a~linear solution, the variable-set is updated cell-wise in two stages. In the first stage, all fractions are kept within the physical interval [0, 1]. The second stage is a local Appleyard chopping of the update (Appleyard and Cheshire 1982), such that it is smaller than a predefined value (e.g., 0.1 for overall compositions). The value of $10^{-4}$ is used for both $\epsilon_p$ and $\epsilon_t$.

\subsection{Case 6: one-dimensional gas injection}

The specification of 1D compositional model is given in Table \ref{tab:specification_comp}. Pressure is kept constant at the both injection and production ends.

\begin{table}[!htb]
\centering
\caption{Specification of 1D compositional model}
\label{tab:specification_comp}
\begin{tabular}{|c|c|c|}
\hline
Parameter                  &  Value            & Unit   \\ \hline
NX                         &  500              &        \\ \hline
DX / DY / DZ               &  32.8 / 32.8 / 32.8     & ft      \\ \hline
Permeability               &  1000             & mD     \\ \hline
Porosity                   &  0.2              &        \\ \hline
Rock compressibility       &  $6.9\cdot 10^{-7}$ & 1/psi  \\ \hline
Total simulation time      &  400              & day    \\ \hline
\end{tabular}
\end{table}

\subsection{Case 6a: three-component fluid}

\begin{figure}[!htb]
\centering
\subfloat[Gas saturation]{
\includegraphics[scale=0.46]{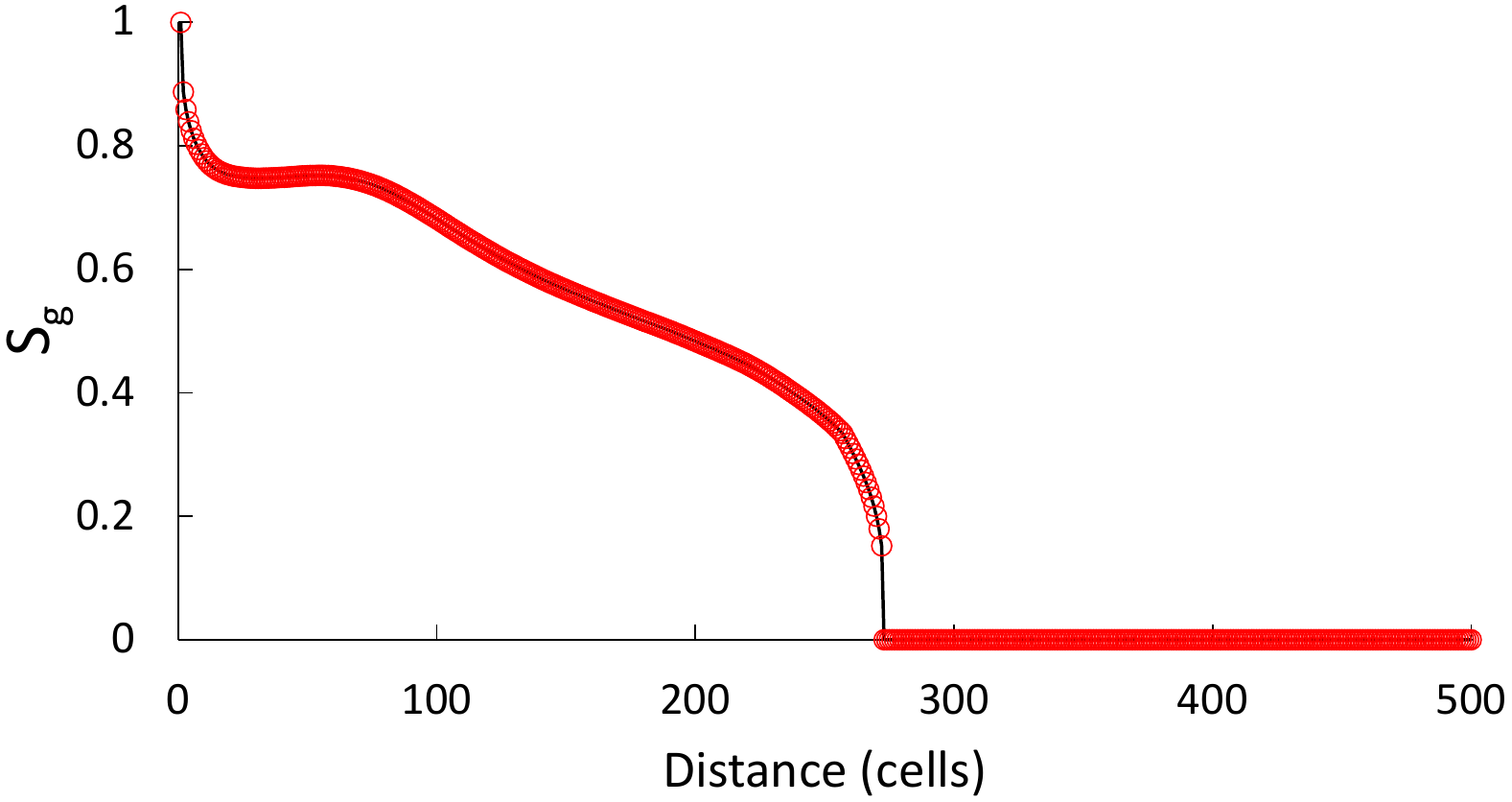}}
\\
\subfloat[Overall composition of $\textrm{C}_1$]{
\includegraphics[scale=0.46]{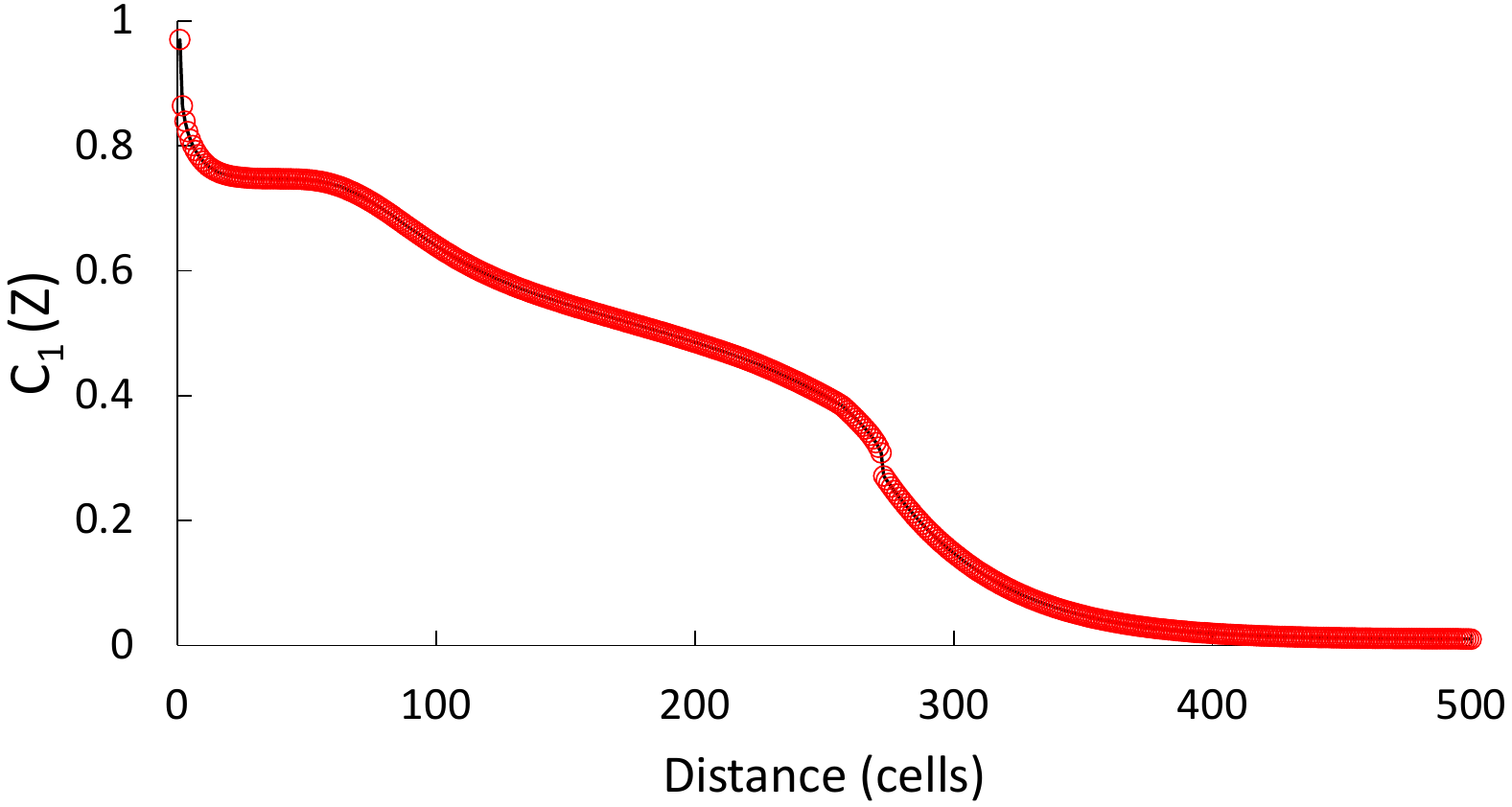}}
\caption{Gas saturation and overall composition profiles of Case 6a: three-component fluid.}
\label{fig:sg_z_6a}
\end{figure}

The example is for a three-component fluid system where the initial oil is made of~$\left \{ \textrm{C}_1 (1 \%), \textrm{C}_4 (50 \%), \textrm{C}_{10} (49 \%) \right \}$, at an initial pressure of 725.2 psi and a temperature of~373~$\textrm{K}$. Production pressure is 580.1 psi, injection pressure is 1740.4 psi, and injection gas is a mixture of $\left \{ \textrm{C}_1 (97 \%), \textrm{C}_4 (2 \%), \textrm{C}_{10} (1 \%) \right \}$. Cubic relative permeabilities are employed. Gas saturation and overall composition profiles are plotted in \textbf{Fig.~\ref{fig:sg_z_6a}}. 

\begin{figure}[!htb]
\centering
\subfloat[Tight tolerances]{
\includegraphics[scale=0.4]{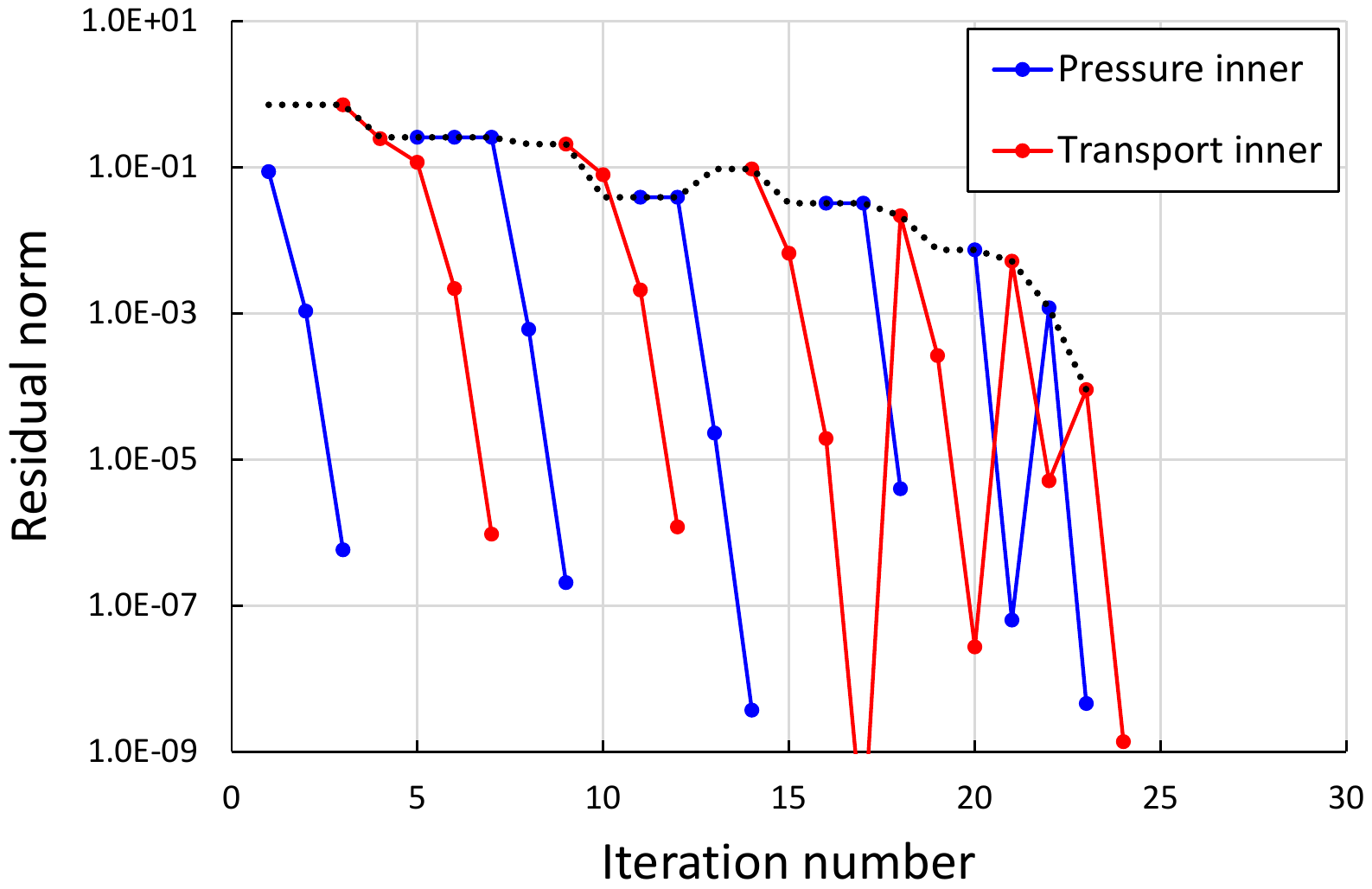}}
\
\subfloat[Adaptive relative relaxation]{
\includegraphics[scale=0.4]{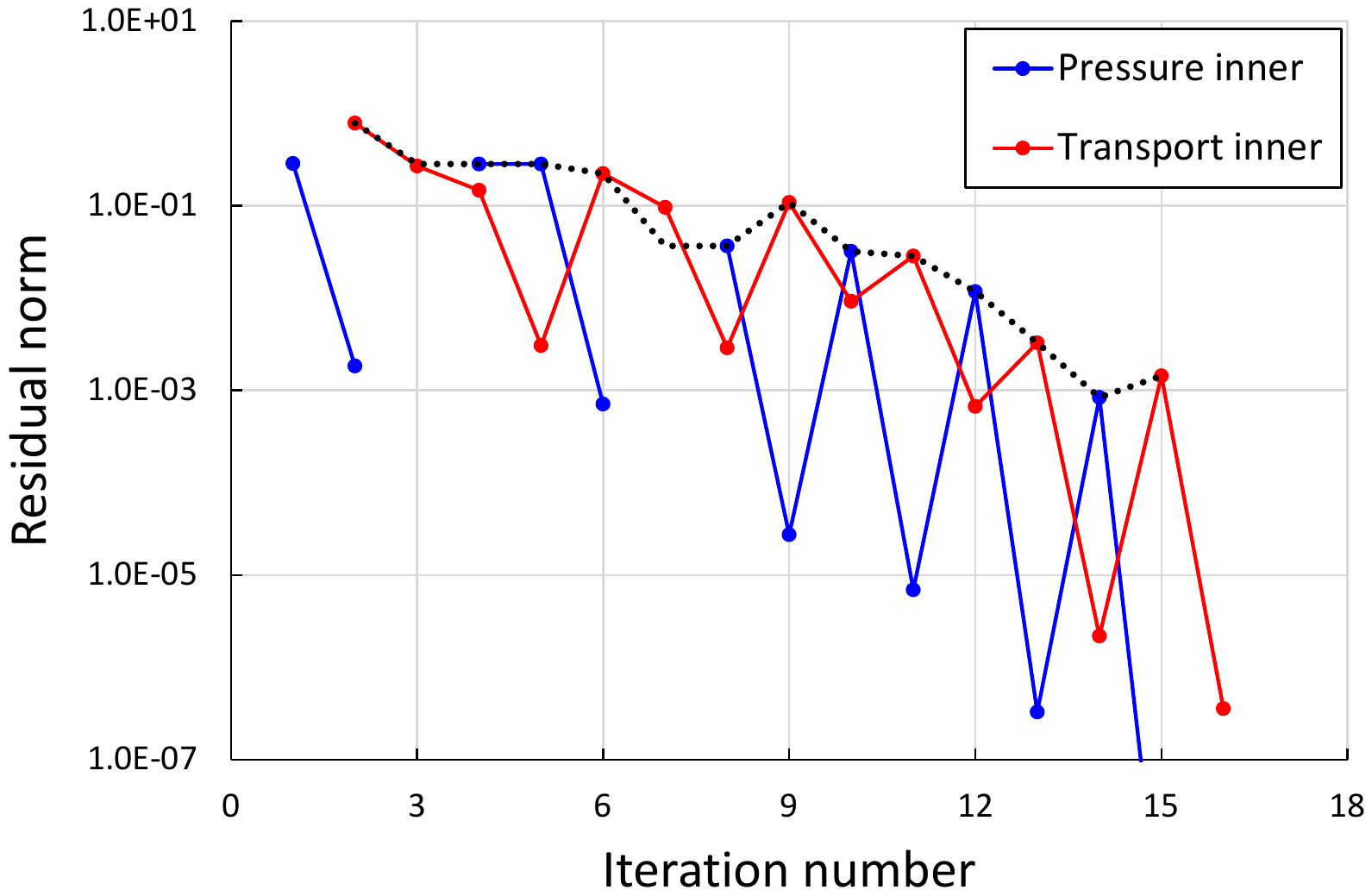}}
\caption{Residual history of the inner solvers for Case 6a: three-component fluid.}
\label{fig:resi_6a}
\end{figure}

\textbf{Fig.~\ref{fig:resi_6a}} plots the residual history of the inner solvers for the second timestep. We can see that less accurate inner solutions are more effective in terms of total iterations. Compared with tight tolerances, the adaptive strategy provides relative tolerances based on the problem's convergence rate.

\subsection{Case 6b: four-component fluid}

\begin{figure}[!htb]
\centering
\subfloat[Gas saturation]{
\includegraphics[scale=0.46]{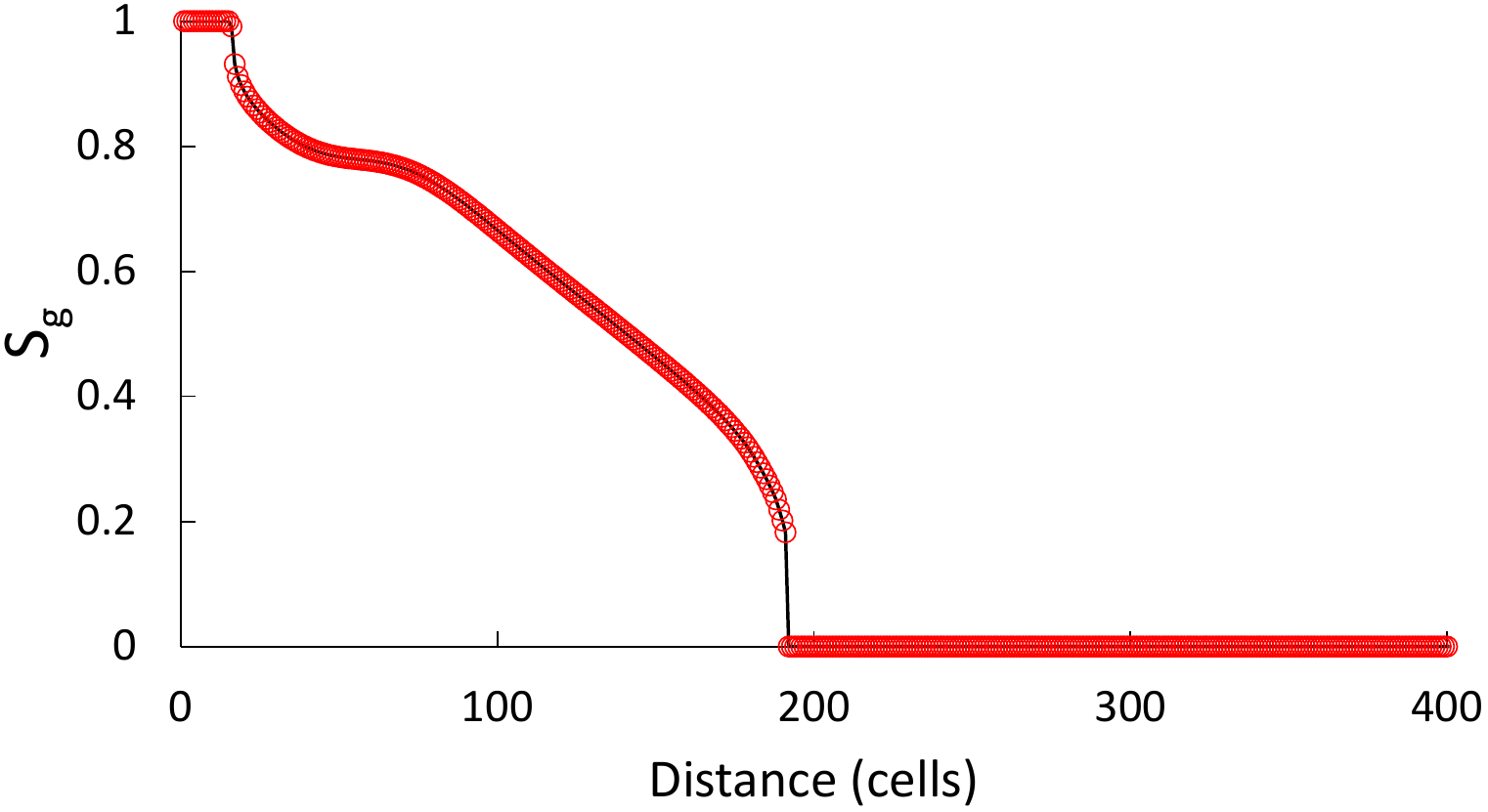}}
\\
\subfloat[Overall composition of $\textrm{C}_1$]{
\includegraphics[scale=0.45]{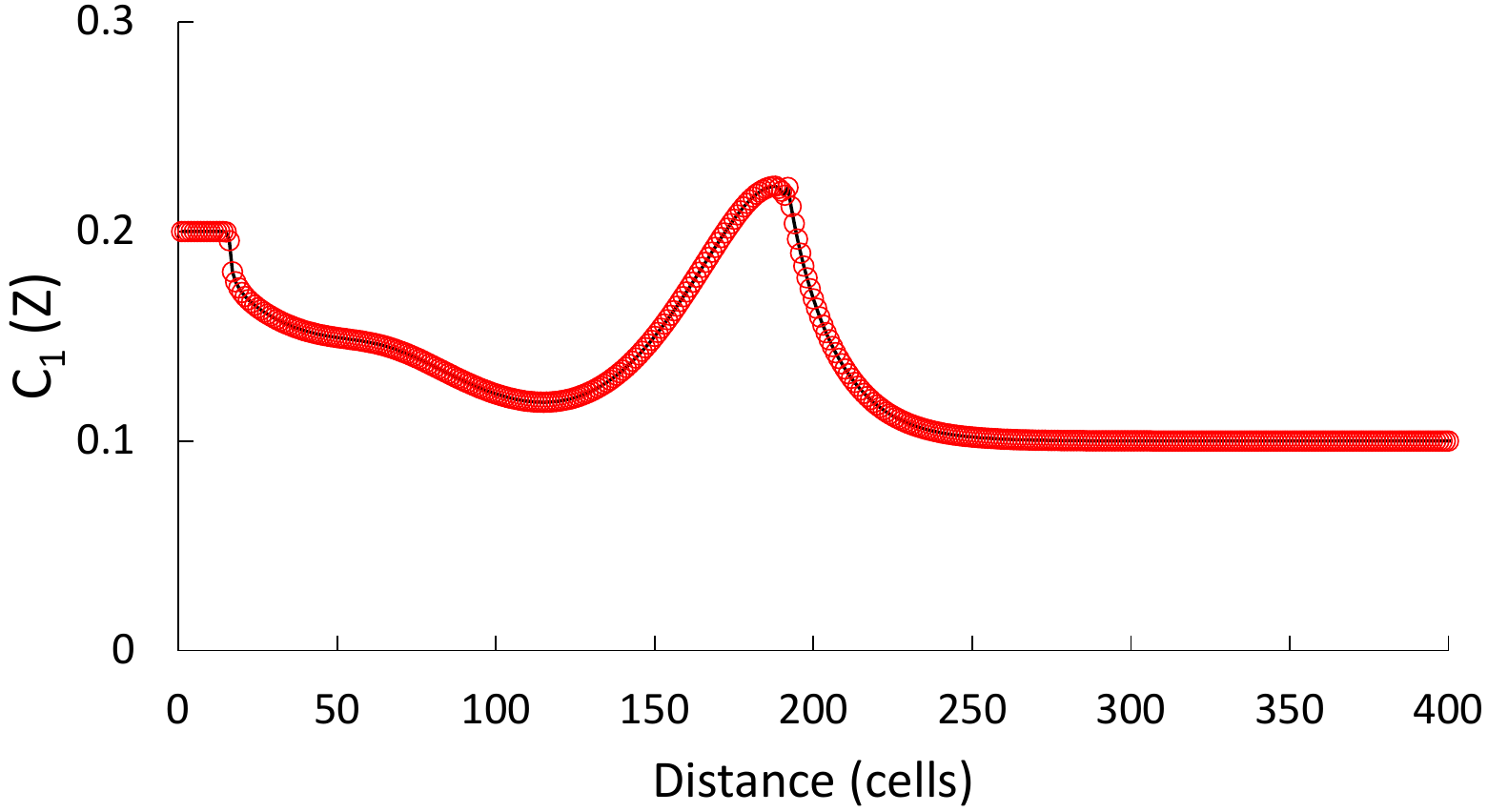}}
\\
\subfloat[Overall composition of $\textrm{CO}_2$]{
\includegraphics[scale=0.46]{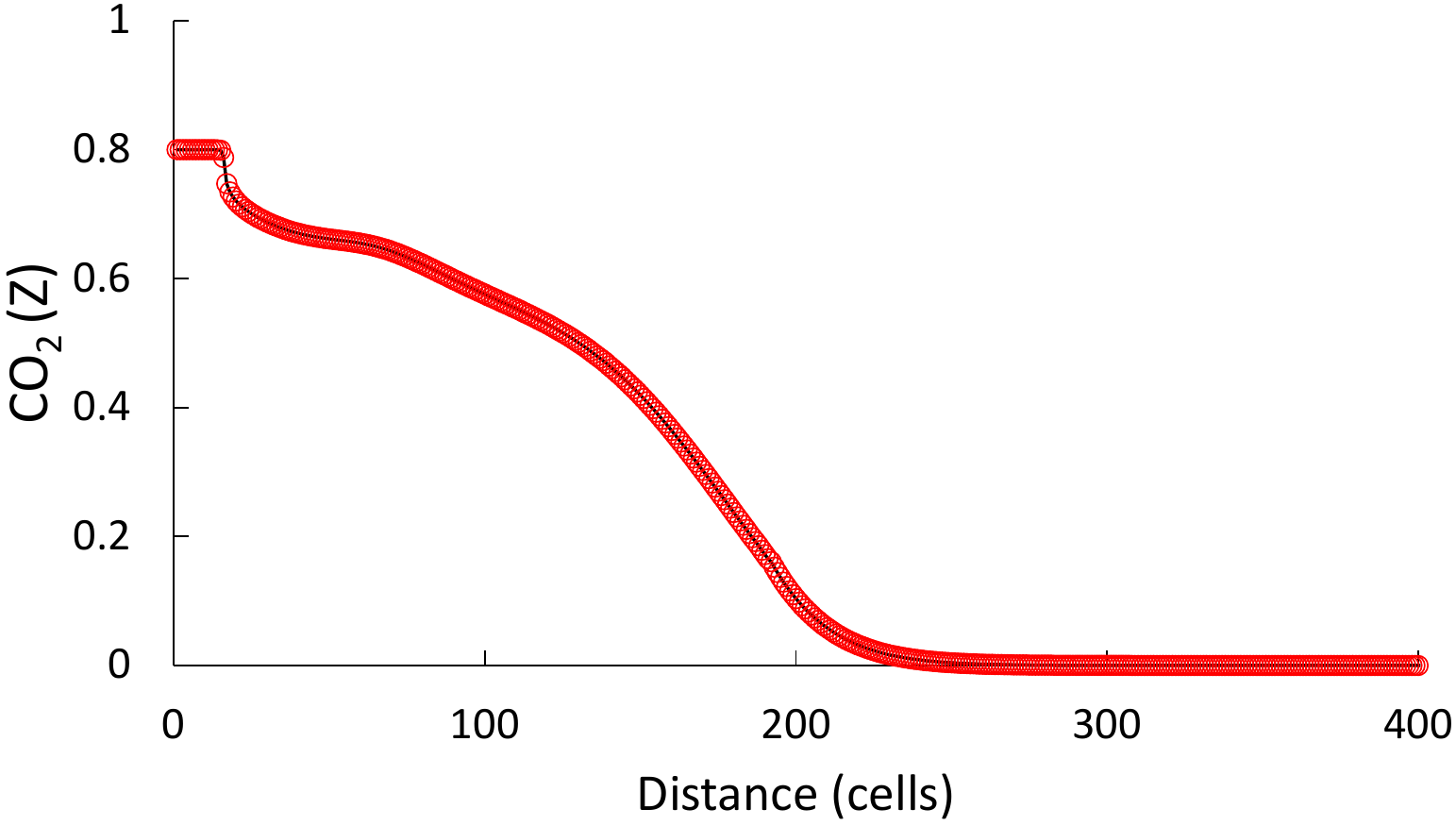}}
\caption{Gas saturation and overall composition profiles of Case 6b: four-component fluid.}
\label{fig:sg_z_6b}
\end{figure}

The initial compositions are $\left \{ \textrm{C}_1 (10 \%), \textrm{CO}_2 (0 \%), \textrm{C}_4 (30 \%), \textrm{C}_{10} (60 \%) \right \}$, at an initial pressure of 870.2 psi and a temperature of 360 $\textrm{K}$. Production pressure is 725.2 psi. Injection pressure is 1305.3 psi, and injection gas mixture is $\left \{ 20 \%, 80 \%, 0 \%, 0 \% \right \}$. Gas saturation and overall composition profiles are plotted in \textbf{Fig.~\ref{fig:sg_z_6b}}. 

\begin{figure}[!htb]
\centering
\subfloat[Case 6a]{
\includegraphics[scale=0.46]{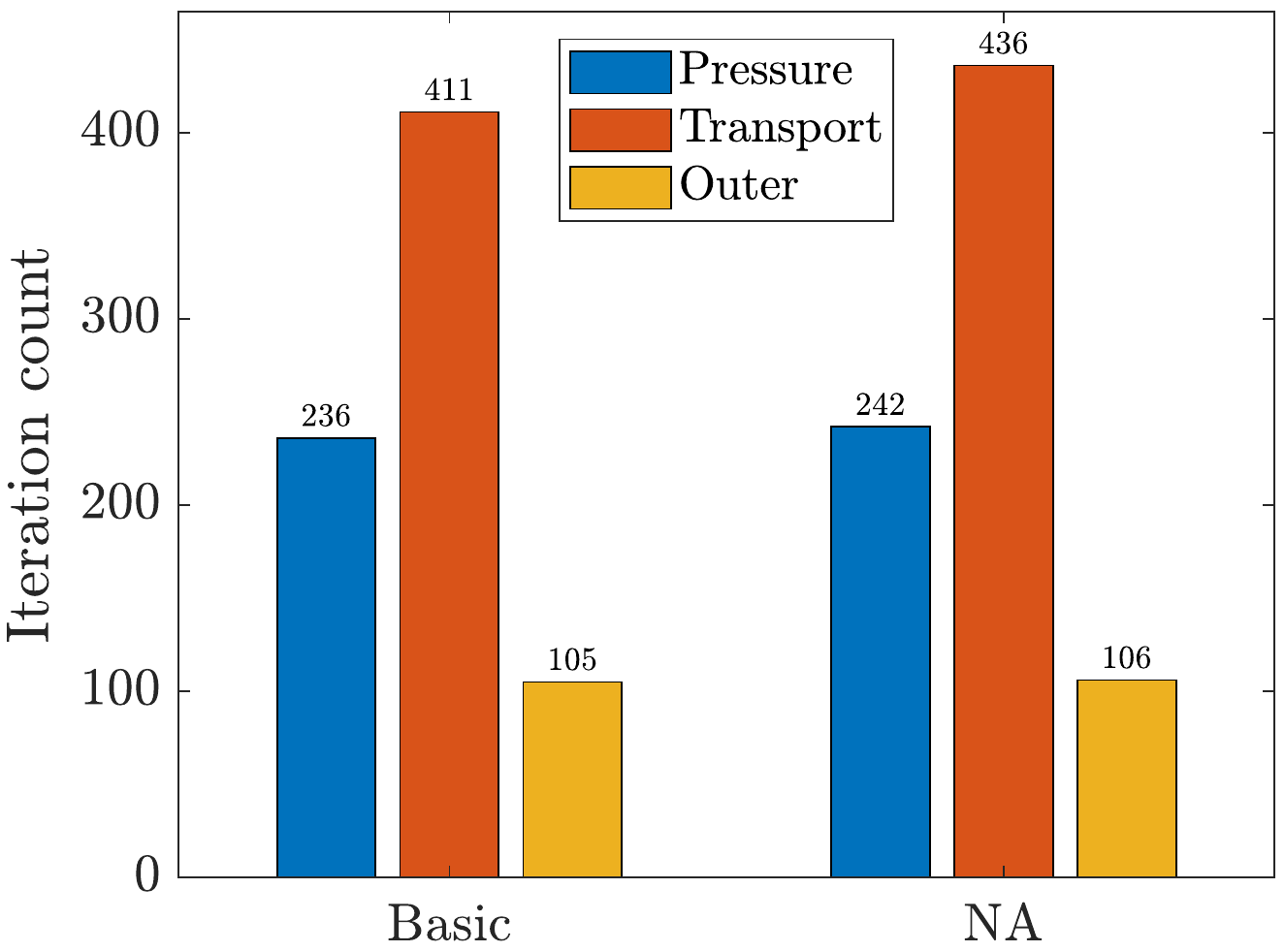}}
\
\subfloat[Case 6b]{
\includegraphics[scale=0.46]{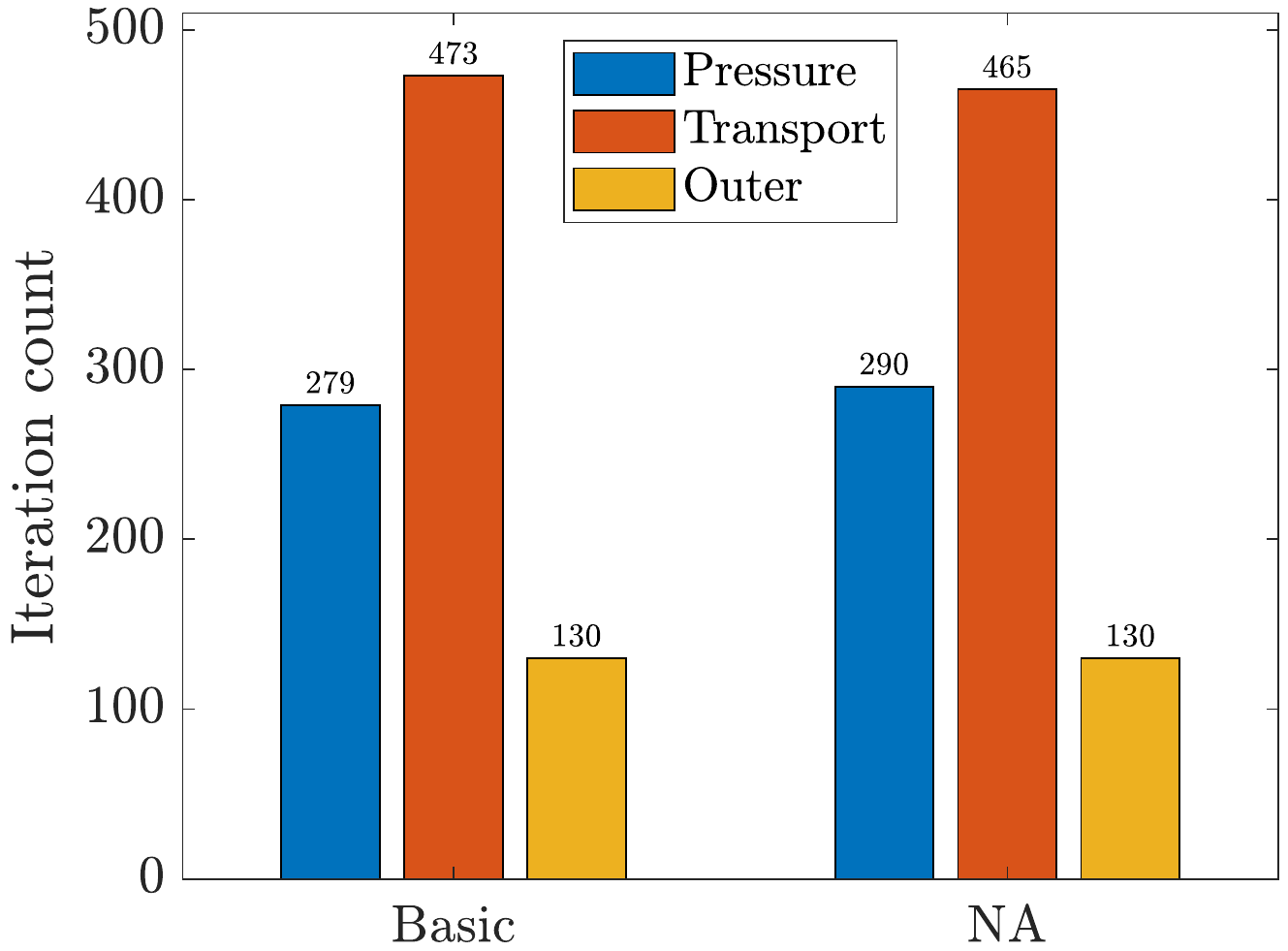}}
\caption{Iteration performance of SFI with and without the nonlinear acceleration for Case 6: one-dimensional gas injection.}
\label{fig:out_iter_6}
\end{figure}

\begin{figure}[!htb]
\centering
\subfloat[Case 6a]{
\includegraphics[scale=0.46]{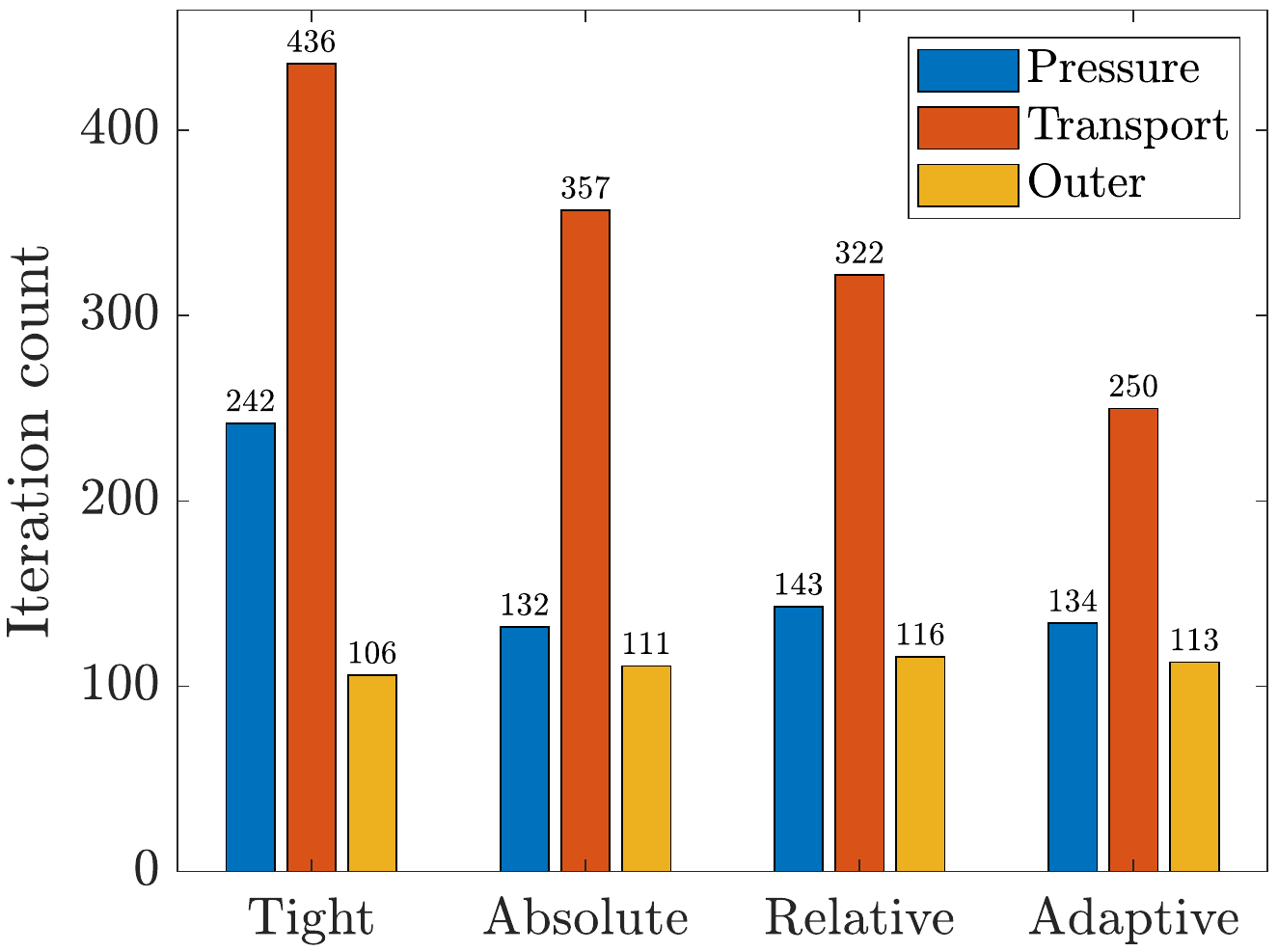}}
\
\subfloat[Case 6b]{
\includegraphics[scale=0.46]{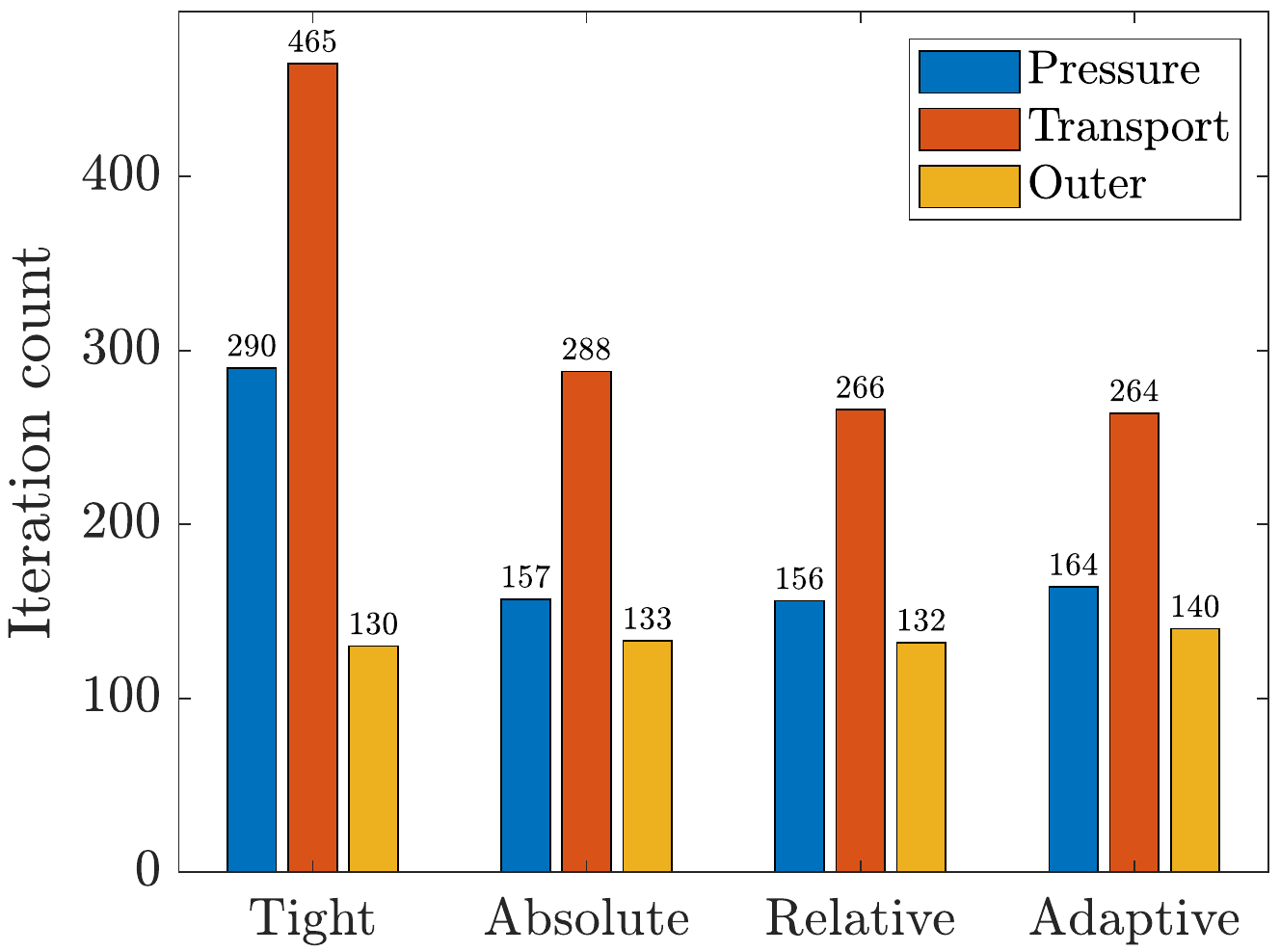}}
\caption{Iteration performance of Case 6: one-dimensional gas injection.}
\label{fig:iter_6}
\end{figure}

We compared the outer-loop performance of SFI with and without the nonlinear acceleration. The results are summarized in \textbf{Fig.~\ref{fig:out_iter_6}}. In this simulation case, the basic SFI method shows a good convergence behavior, indicating weak coupling effects between the sub-problems. The iteration performance of the inner solvers for Case 6 is summarized in \textbf{Fig.~\ref{fig:iter_6}}. As can be seen, the inexact methods effectively resolve the over-solving issue, and thus significantly improve the overall efficiency.

\subsection{Case 7: gravity-driven lock exchange}

Again we tested a lock-exchange problem on a homogeneous $60 \times 60$ square model. A~uniform cell size 3.048 m is specified, and porosity is 0.2. Oil ($\textrm{C}_{10}$) initially occupies the left half of the domain, while gas ($\textrm{C}_1$) fills the right half. Initial pressure is 2001.5 psi, with a temperature of~373 $\textrm{K}$. Total simulation time is 20 days. 

\begin{figure}[!htb]
\centering
\subfloat[Basic]{
\includegraphics[scale=0.47]{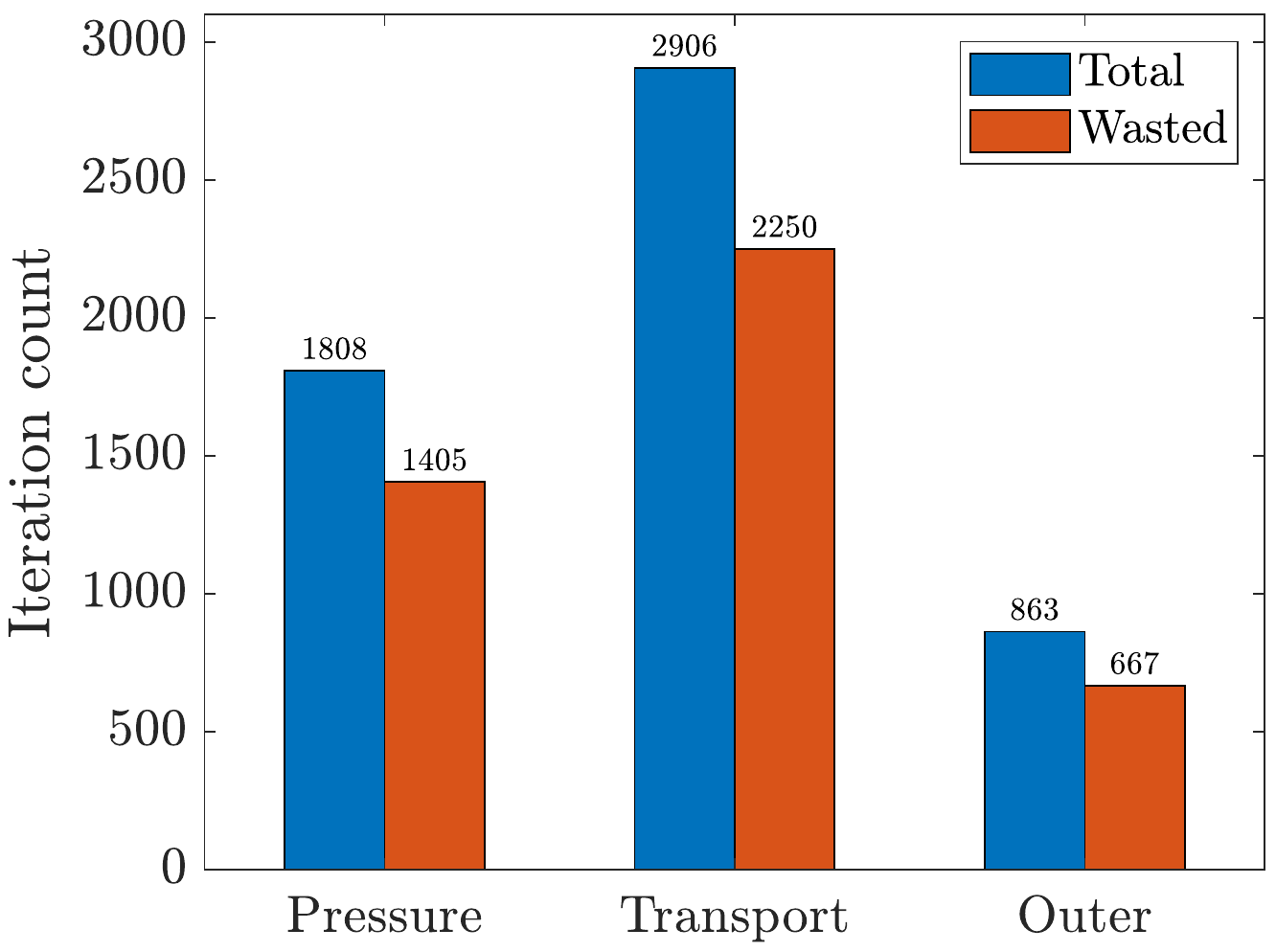}}
\
\subfloat[NA]{
\includegraphics[scale=0.47]{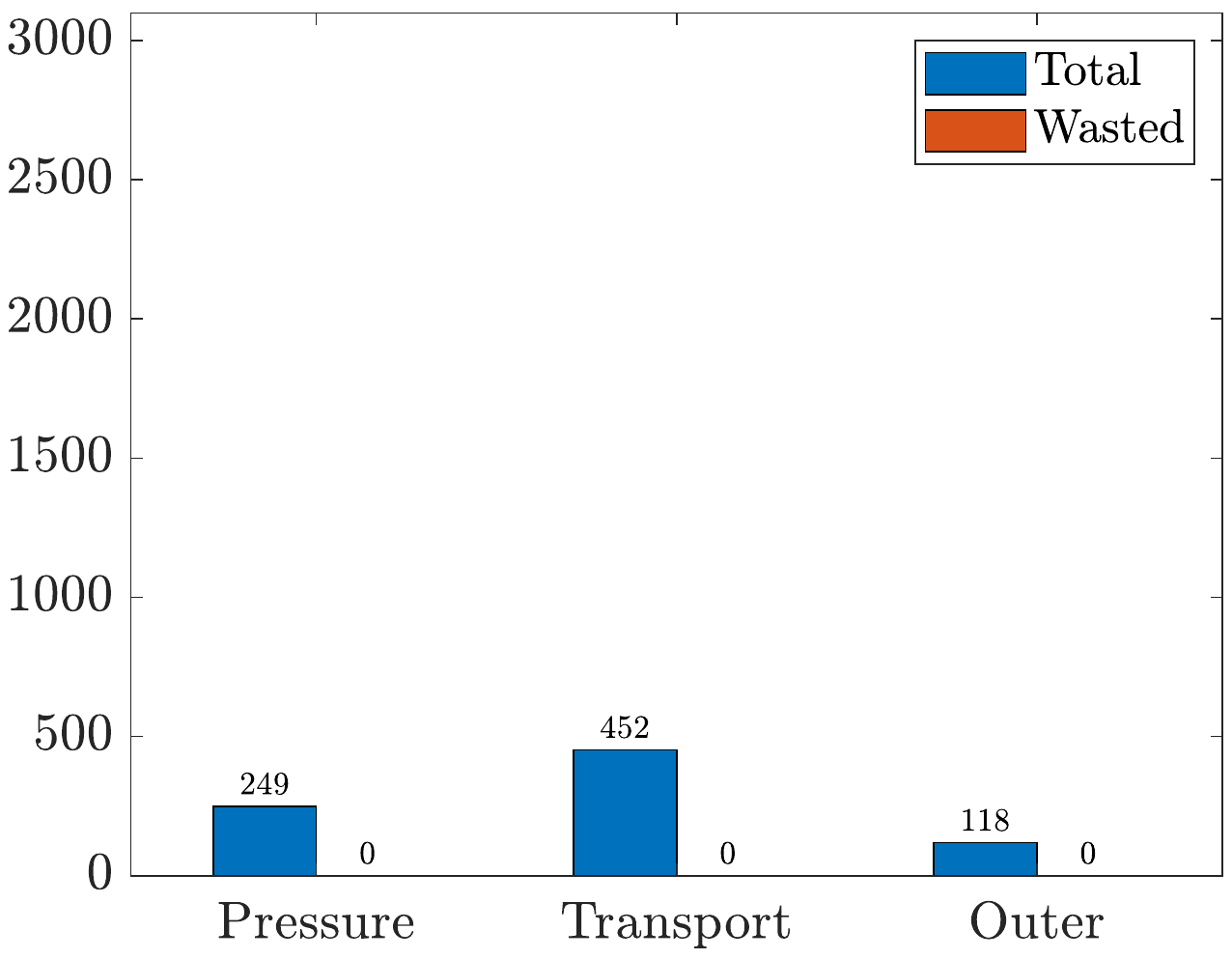}}
\caption{Iteration performance of SFI with and without the nonlinear acceleration for Case 7: gravity-driven lock exchange.}
\label{fig:out_iter_7}
\end{figure}

Although the model is homogeneous and has no wells, strong coupling exists between the flow and transport. We compared the outer-loop performance of SFI with and without the nonlinear acceleration. \textbf{Fig.~\ref{fig:out_iter_7}} reports the total number of iterations. We also report the number of wasted iterations that correspond to the iterations spent on unconverged timesteps. It can be seen that the basic SFI process is quite inefficient, suffering from severe restrictions on allowable time-step size. By comparison, the NA method leads to superior convergence performance. Here CPU timing results are not further provided because the QN update is cheap and thus the additional cost from NA is assumed to be negligible. As a result, the overall simulation time will be proportional to the numbers of outer iterations.

\begin{figure}[!htb]
\centering
\includegraphics[scale=0.58]{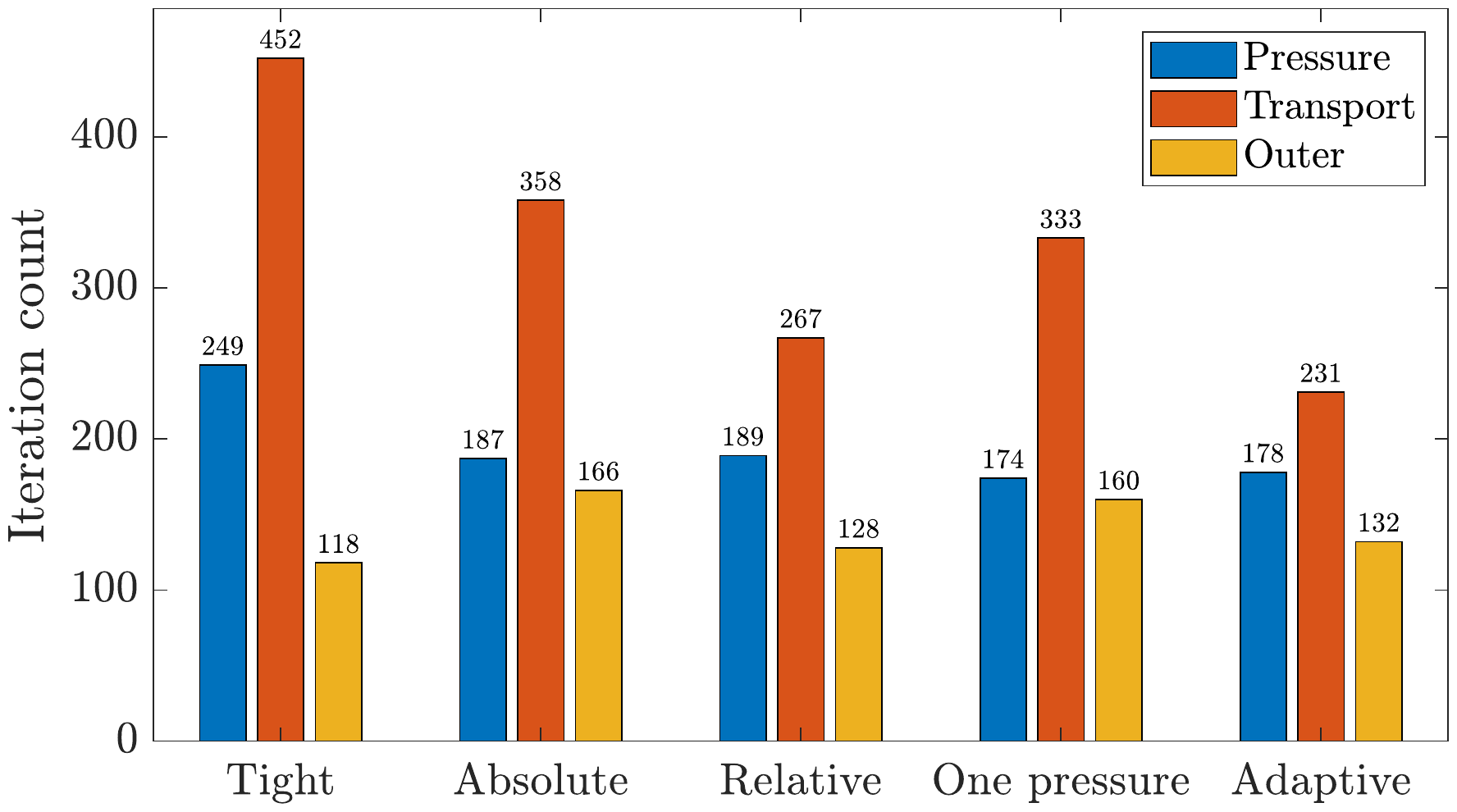}
\caption{Iteration performance of Case 7: gravity-driven lock exchange.}
\label{fig:iter_7}
\end{figure}

The iteration performance of the inner solvers for Case 7 is summarized in \textbf{Fig.~\ref{fig:iter_7}}. The adaptive strategy provides the optimal efficiency. In addition, we also ran a test with only one iteration for pressure. We can see that the transport iteration count largely increases, indicating that the quality of the pressure solution has a big impact on the transport solver.


\subsection{Case 8: viscous and gravitational forces, SPE 10 model}

\begin{figure}[!htb]
\centering
\subfloat[Gas saturation]{
\includegraphics[scale=0.36]{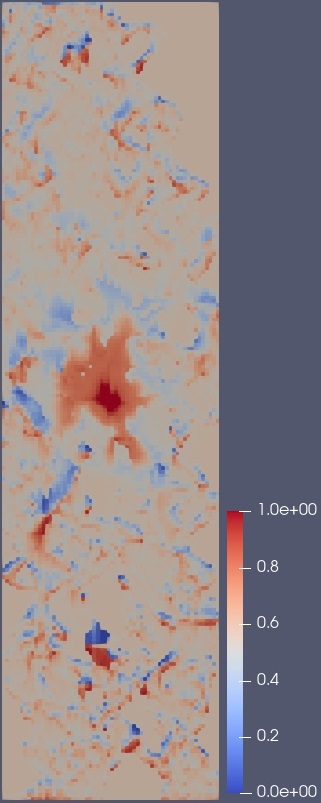}}
\qquad \ \ 
\subfloat[Overall composition of $\textrm{CO}_2$]{
\includegraphics[scale=0.36]{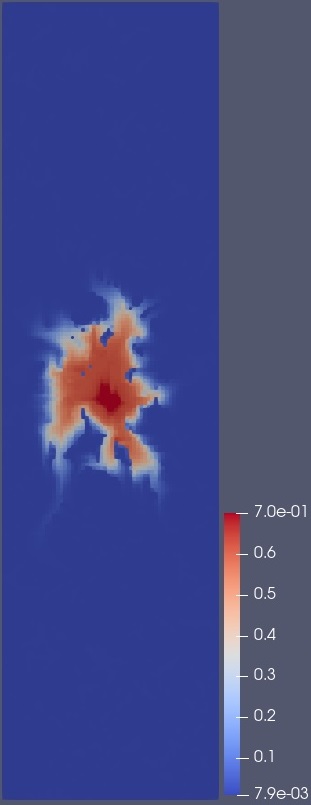}}
\caption{Gas saturation and overall composition profiles of Case 8: viscous and gravitational forces, SPE 10 model.}
\label{fig:re_comp_8}
\end{figure}

We further tested a case on the bottom layer of the SPE 10 model. The domain is~vertically placed, and thus gravity is in effect. A uniform cell size 9.84~ft is specified, and porosity is~0.2. The four-component fluid system is used with initial compositions $\left \{ \textrm{C}_1 (50 \%), \textrm{CO}_2 (1 \%), \textrm{C}_4 (29 \%), \textrm{C}_{10} (20 \%) \right \}$, at an initial pressure of 1450~psi and a~temperature of 373~$\textrm{K}$. Injection pressure is 1740~psi, and injection gas mixture as~$\left \{ 28 \%, 70 \%, 1 \%, 1 \% \right \}$. Total simulation time is 2~days. Gas saturation and overall composition profiles are plotted in \textbf{Fig.~\ref{fig:re_comp_8}}. 

We ran another case (Case 8b) with total simulation time changed to 30 days. Initial compositions become $\left \{ 25 \%, 1 \%, 9 \%, 65 \% \right \}$. The heterogeneous model is more challenging for the outer loop of SFI to converge, due to a large variation in the CFL numbers throughout the domain.

\begin{figure}[!htb]
\centering
\subfloat[Basic]{
\includegraphics[scale=0.47]{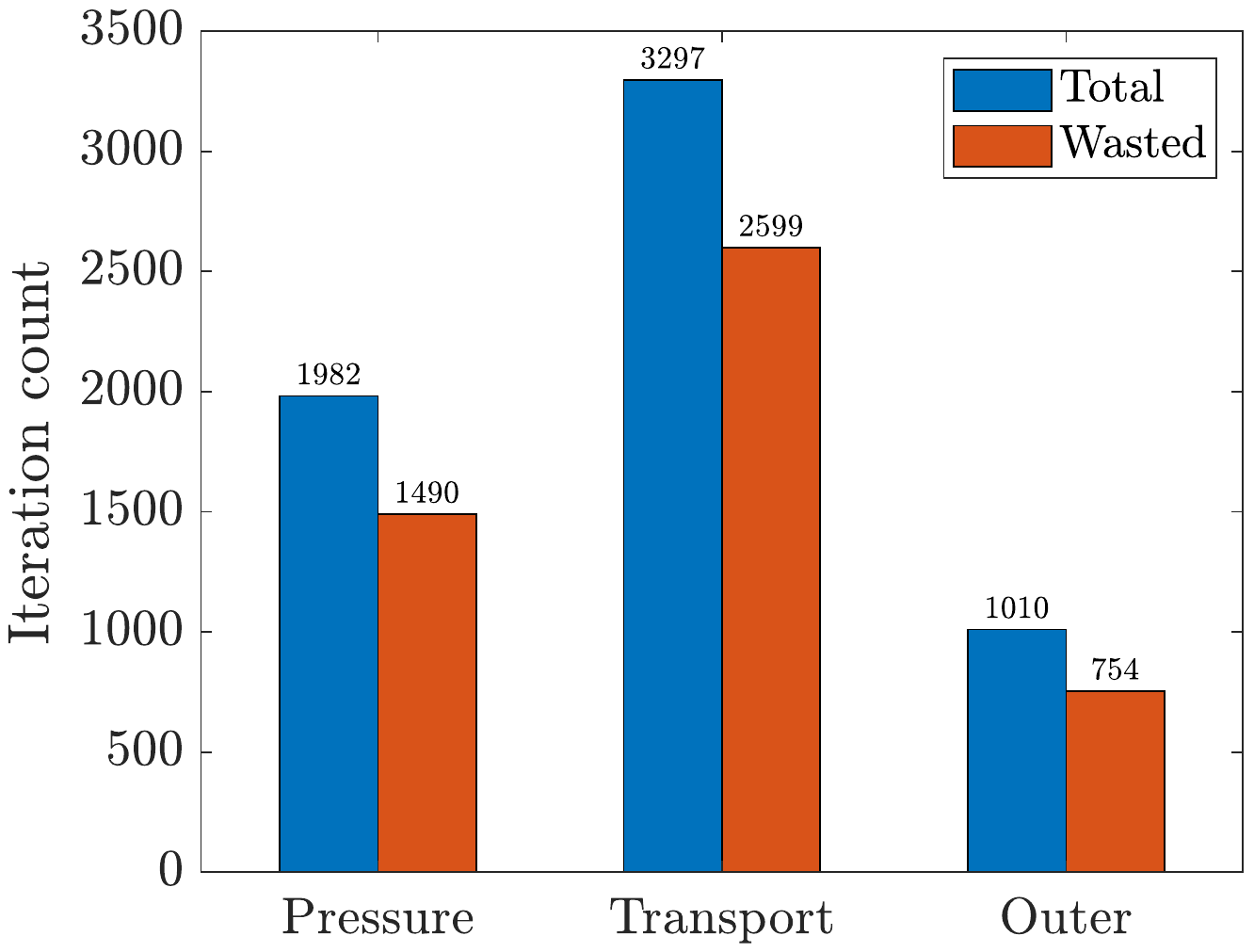}}
\
\subfloat[NA]{
\includegraphics[scale=0.47]{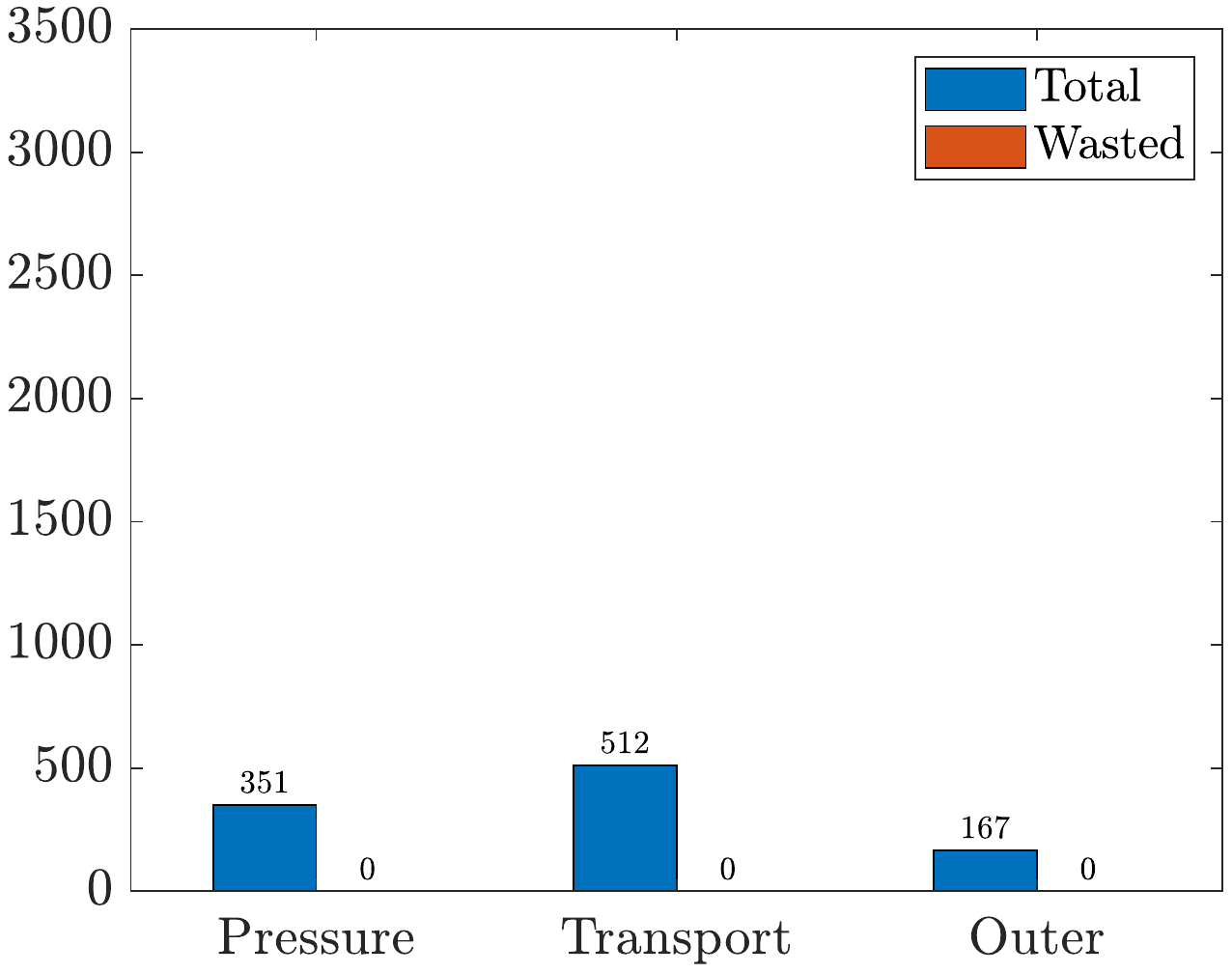}}
\caption{Iteration performance of SFI with and without the nonlinear acceleration for Case 8a.}
\label{fig:out_iter_8a}
\end{figure}

\begin{figure}[!htb]
\centering
\subfloat[Basic]{
\includegraphics[scale=0.47]{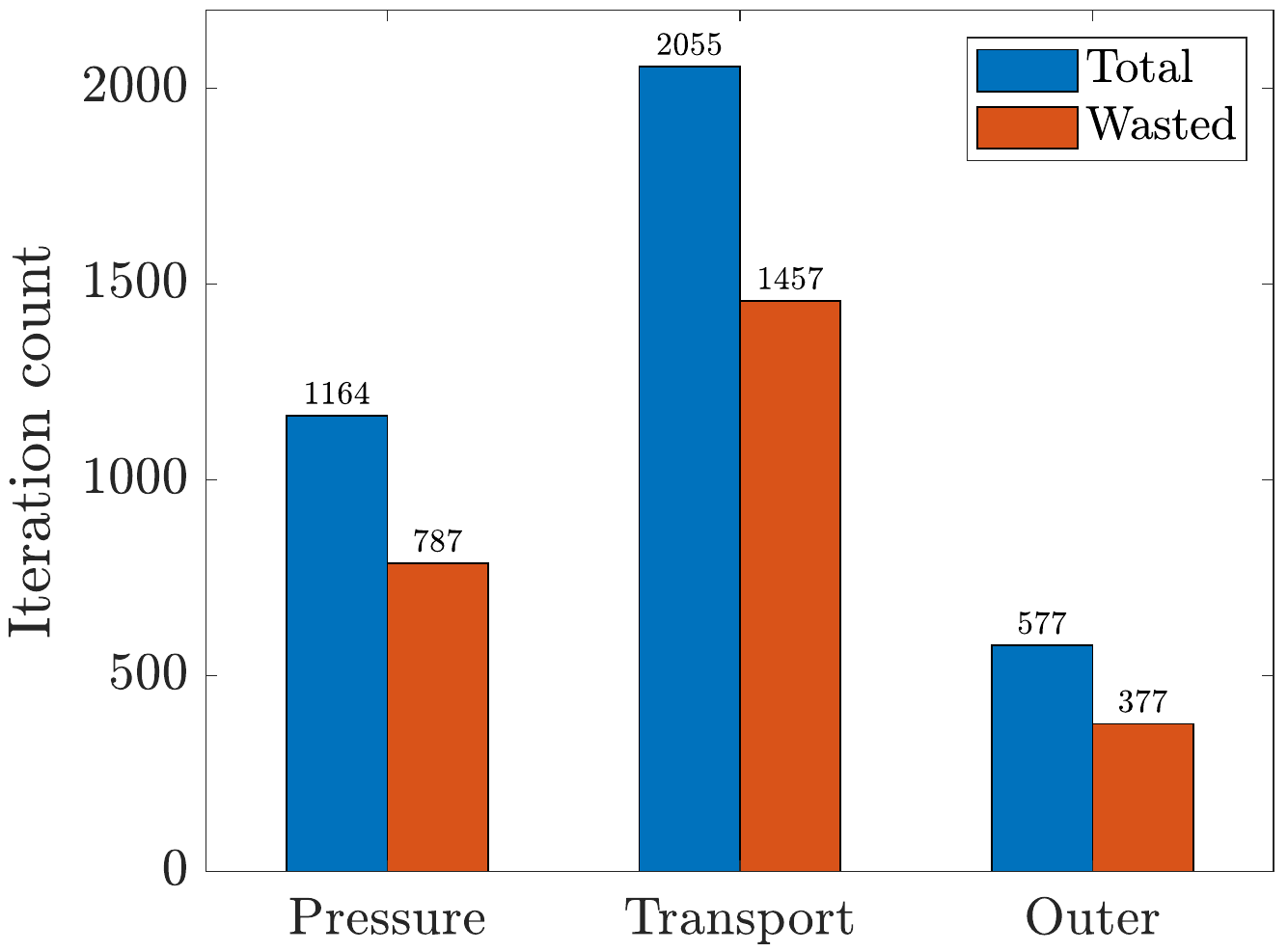}}
\
\subfloat[NA]{
\includegraphics[scale=0.47]{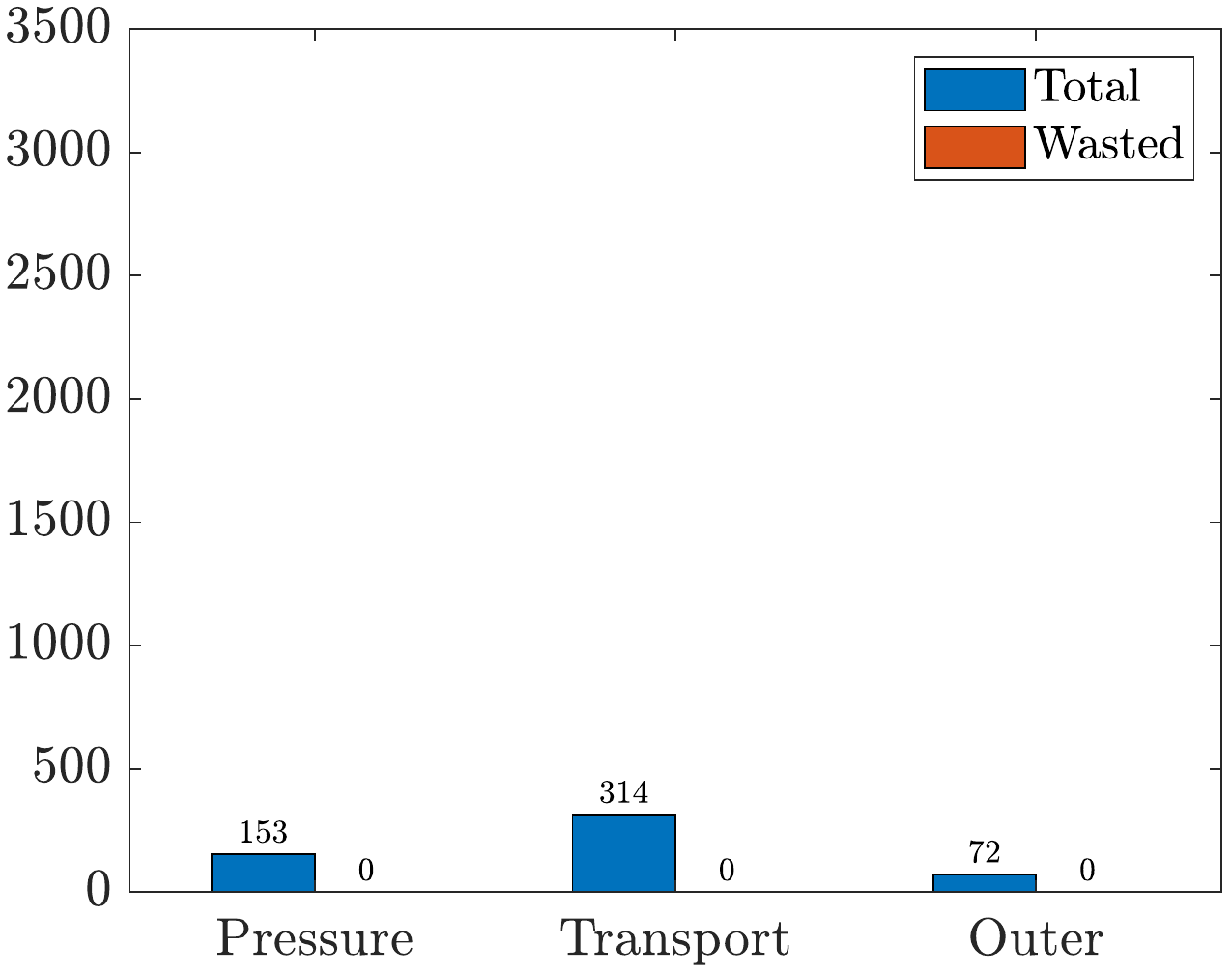}}
\caption{Iteration performance of SFI with and without the nonlinear acceleration for Case 8b.}
\label{fig:out_iter_8b}
\end{figure}

\begin{figure}[!htb]
\centering
\subfloat[Case 8a]{
\includegraphics[scale=0.46]{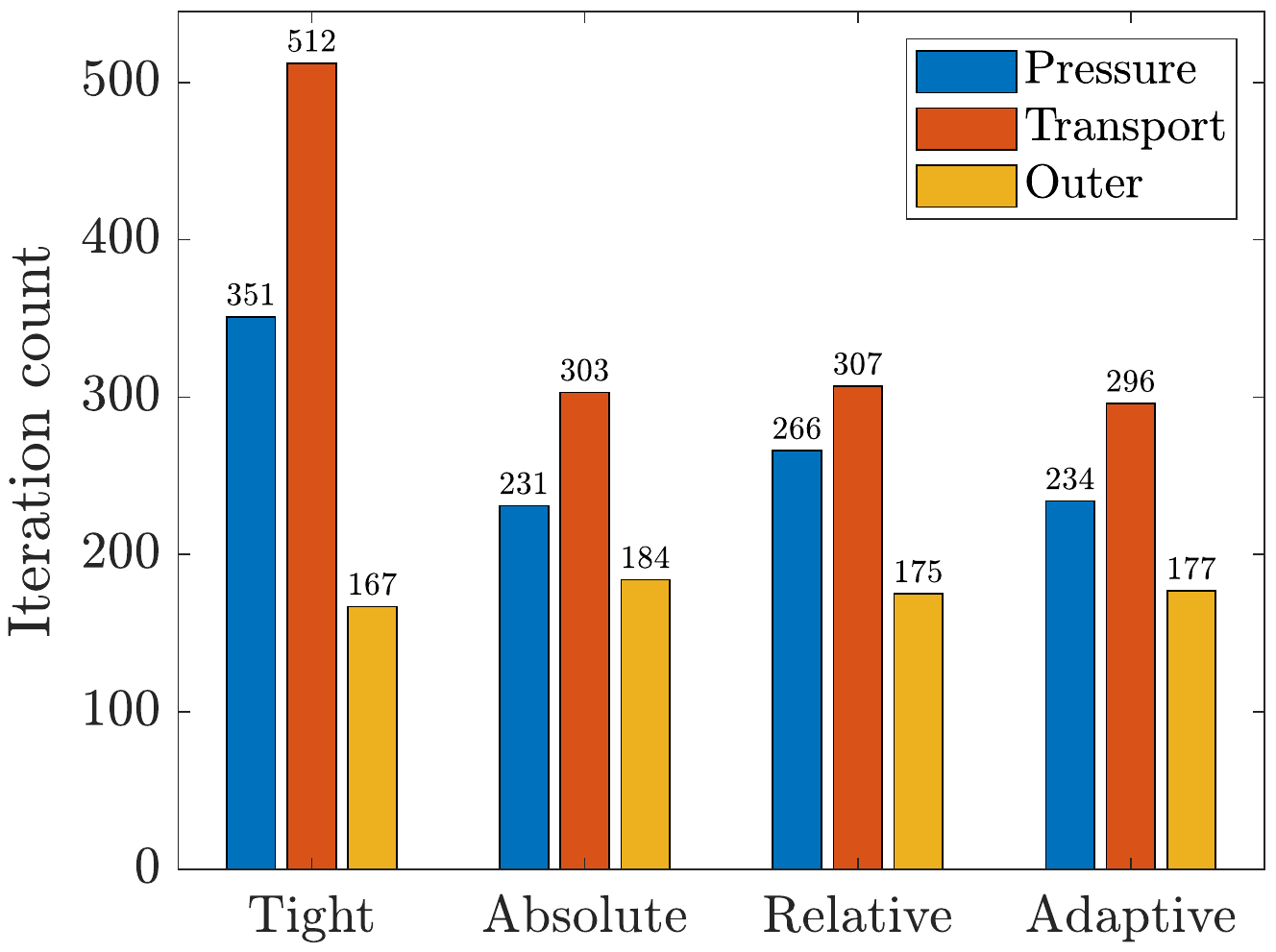}}
\
\subfloat[Case 8b]{
\includegraphics[scale=0.46]{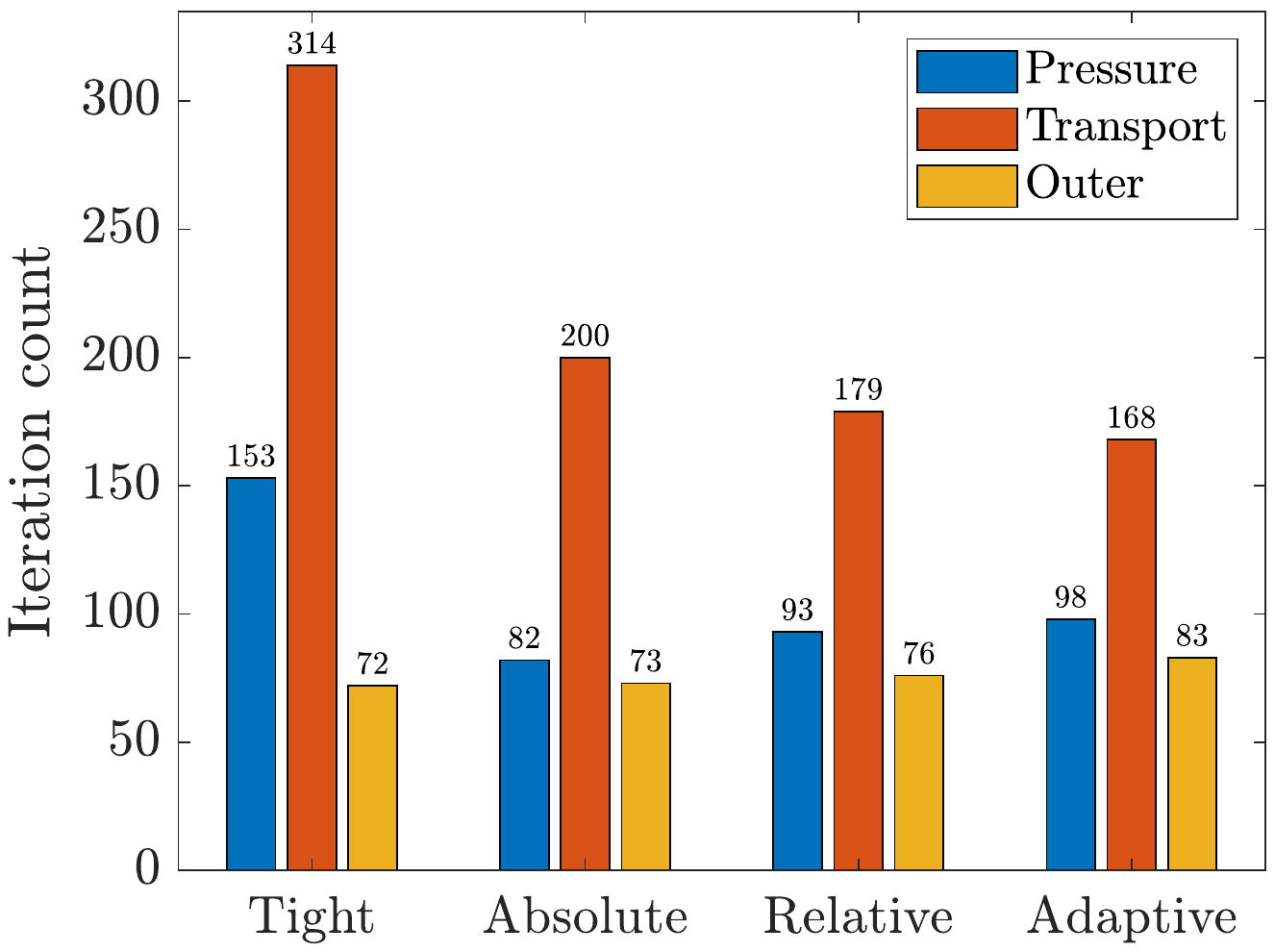}}
\caption{Iteration performance of Case 8: viscous and gravitational forces, SPE 10 model.}
\label{fig:iter_8}
\end{figure}

We first compared the performance of SFI with and without NA. The results are summarized in \textbf{Fig.~\ref{fig:out_iter_8a}} and \textbf{\ref{fig:out_iter_8b}}. The basic SFI induces huge numbers of outer iterations and wasted computations. In contrast, SFI-NA does not require any time-step cuts. The iteration performance of the inner solvers for Case 8 is summarized in \textbf{Fig.~\ref{fig:iter_8}}. The inexact methods show consistent improvements upon the basic method.


\section{Summary}

In the standard SFI method, the sub-problems are usually solved to high precision at every outer iteration. This may result in wasted computations that contribute little or no progress towards the coupled solution. We find that there is no need for one sub-problem to strive for perfection (`over-solving'), while the coupled (outer) residual remains high due to the other sub-problem. The objective of this work is to minimize the cost of inner solvers while not degrading the convergence rate of SFI.

We first extended a nonlinear-acceleration (NA) framework to multi-component compositional models, for ensuring robust outer-loop convergence. We then developed inexact-type methods that alleviate over-solving. An adaptive strategy was proposed for providing relative tolerances based on the convergence rates of coupled problems.

The new SFI solver was tested using several complex cases.~The problems involve multi-phase and EoS-based compositional fluid systems.~The results demonstrate that the basic SFI method is quite inefficient. Away from a coupled solution, additional accuracy achieved in inner solvers is wasted, contributing to little reduction of the overall residual. By comparison, the inexact method adaptively provides the relative tolerances adequate for the sub-problems.~We show that the new solver consistently improves the efficiency of the sub-problems, and nearly 50\% reduction in the total simulation cost is obtained.

\section*{Acknowledgements}

This work was supported by the Chevron/Schlumberger INTERSECT Research \& Prototyping project. The authors thank Chevron for permission to publish the paper. Petroleum Research Institute for Reservoir Simulation (SUPRI-B) at Stanford University is gratefully acknowledged for providing the AD-GPRS platform.

We would like to thank Hamdi Tchelepi for inspiring discussions on the inexact methods. We also thank Jiawei Li, Huanquan Pan at Stanford University, and Olav M{\o}yner at SINTEF Digital, for constructive discussions on AD-GPRS, compositional model, sequential formulation, etc.

\section*{References}

Acs, G., Doleschall, S. and Farkas, E., 1985. General purpose compositional model. Society of Petroleum Engineers Journal, 25(04), pp.543-553.

Appleyard, J. R., Cheshire, I. M., 1982. The cascade method for accelerated convergence in implicit simulators. In European Petroleum Conference, 259–290.

Aziz, K., Settari, A., 1979. Petroleum Reservoir Simulation. Chapman \& Hall.

Baker, L.E., 1988, January. Three-phase relative permeability correlations. In SPE Enhanced Oil Recovery Symposium. Society of Petroleum Engineers.

Birken, P., 2015. Termination criteria for inexact fixed-point schemes. Numerical Linear Algebra with Applications, 22(4), pp.702-716.

Brenier, Y. and Jaffré, J., 1991. Upstream differencing for multiphase flow in reservoir simulation. SIAM journal on numerical analysis, 28(3), pp.685-696.

Cao, H., 2002. Development of techniques for general purpose simulators (Doctoral dissertation, Stanford University).

Coats, K.H., 1980. An equation of state compositional model. Society of Petroleum Engineers Journal, 20(05), pp.363-376.

Coats, K.H., 2000. A note on IMPES and some IMPES-based simulation models. SPE Journal, 5(03), pp.245-251.

Collins, D.A., Nghiem, L.X., Li, Y.K. and Grabonstotter, J.E., 1992. An efficient approach to adaptive-implicit compositional simulation with an equation of state. SPE reservoir engineering, 7(02), pp.259-264.

Christie, M.A. and Blunt, M.J., 2001, January. Tenth SPE comparative solution project: A comparison of upscaling techniques. In SPE reservoir simulation symposium. Society of Petroleum Engineers.

Dembo, R.S., Eisenstat, S.C. and Steihaug, T., 1982. Inexact newton methods. SIAM Journal on Numerical analysis, 19(2), pp.400-408.

Dawson, C.N., Klíe, H., Wheeler, M.F. and Woodward, C.S., 1997. A parallel, implicit, cell‐centered method for two‐phase flow with a preconditioned Newton–Krylov solver. Computational Geosciences, 1(3), pp.215-249.

Degroote, J., Souto-Iglesias, A., Van Paepegem, W., Annerel, S., Bruggeman, P. and Vierendeels, J., 2010. Partitioned simulation of the interaction between an elastic structure and free surface flow. Computer methods in applied mechanics and engineering, 199(33), pp.2085-2098.

Eisenstat, S.C. and Walker, H.F., 1996. Choosing the forcing terms in an inexact Newton method. SIAM Journal on Scientific Computing, 17(1), pp.16-32.

Garipov, T.T., Tomin, P., Rin, R., Voskov, D.V. and Tchelepi, H.A., 2018. Unified thermo-compositional-mechanical framework for reservoir simulation. Computational Geosciences, 22(4), pp.1039-1057.

Hajibeygi, H. and Tchelepi, H.A., 2014. Compositional multiscale finite-volume formulation. SPE Journal, 19(02), pp.316-326.

Jenny, P., Lee, S.H. and Tchelepi, H.A., 2006. Adaptive fully implicit multi-scale finite-volume method for multi-phase flow and transport in heterogeneous porous media. Journal of Computational Physics, 217(2), pp.627-641.

Jenny, P., Tchelepi, H.A. and Lee, S.H., 2009. Unconditionally convergent nonlinear solver for hyperbolic conservation laws with S-shaped flux functions. Journal of Computational Physics, 228(20), pp.7497-7512.

Jiang, J. and Younis, R.M., 2017. Efficient C1-continuous phase-potential upwind (C1-PPU) schemes for coupled multiphase flow and transport with gravity. Advances in Water Resources, 108, pp.184-204.

Jiang, J. and Tchelepi, H.A., 2019. Nonlinear acceleration of sequential fully implicit (SFI) method for coupled flow and transport in porous media. Computer Methods in Applied Mechanics and Engineering, 352, pp.246-275.

Klie, H.M., Krylov-secant methods for solving large-scale systems of coupled nonlinear parabolic equations, Ph.D. thesis, Dept. of Computational and Applied Mathematics, Rice University, Houston, TX (September 1996).

Klie, H. and Wheeler, M.F., 2005, January. Krylov-secant methods for accelerating the solution of fully implicit formulations. In SPE Reservoir Simulation Symposium. Society of Petroleum Engineers.

Kuttler, U. and Wall, W.A., 2008. Fixed-point fluid-structure interaction solvers with dynamic relaxation. Computational Mechanics, 43(1), pp.61-72.

Kozlova, A., Li, Z., Natvig, J.R., Watanabe, S., Zhou, Y., Bratvedt, K. and Lee, S.H., 2016. A real-field multiscale black-oil reservoir simulator. SPE Journal.

Lee, S.H., Wolfsteiner, C. and Tchelepi, H.A., 2008. Multiscale finite-volume formulation for multiphase flow in porous media: black oil formulation of compressible, three-phase flow with gravity. Computational Geosciences, 12(3), pp.351-366.

Lee, S.H., Efendiev, Y. and Tchelepi, H.A., 2015. Hybrid upwind discretization of nonlinear two-phase flow with gravity. Advances in Water Resources, 82, pp.27-38.

Lee, S.H. and Efendiev, Y., 2016. C1-Continuous relative permeability and hybrid upwind discretization of three phase flow in porous media. Advances in Water Resources, 96, pp.209-224.

Lie, K.A., Krogstad, S., Ligaarden, I.S., Natvig, J.R., Nilsen, H.M. and Skaflestad, B., 2012. Open-source MATLAB implementation of consistent discretisations on complex grids. Computational Geosciences, 16(2), pp.297-322.

Lie, K.A., M{\o}yner, O., Natvig, J.R., Kozlova, A., Bratvedt, K., Watanabe, S. and Li, Z., 2017. Successful application of multiscale methods in a real reservoir simulator environment. Computational Geosciences, 21(5-6), pp.981-998.

Li, B. and Tchelepi, H.A., 2015. Nonlinear analysis of multiphase transport in porous media in the presence of viscous, buoyancy, and capillary forces. Journal of Computational Physics, 297, pp.104-131.

Li, J., Wong, Z.Y., Tomin, P. and Tchelepi, H., 2019, March. Sequential Implicit Newton Method for Coupled Multi-Segment Wells. In SPE Reservoir Simulation Conference. Society of Petroleum Engineers.

M{\o}yner, O. and Lie, K.A., 2016. A multiscale restriction-smoothed basis method for compressible black-oil models. SPE Journal, 21(06), pp.2-079.

M{\o}yner, O. and Tchelepi, H.A., 2018. A Mass-Conservative Sequential Implicit Multiscale Method for Isothermal Equation-of-State Compositional Problems. SPE Journal, 23(06), pp.2-376.

Moncorgé, A., Tchelepi, H.A. and Jenny, P., 2018. Sequential fully implicit formulation for compositional simulation using natural variables. Journal of Computational Physics, 371, pp.690-711.

Moncorgé, A., M{\o}yner, O., Tchelepi, H.A. and Jenny, P., 2019. Consistent upwinding for sequential fully implicit multiscale compositional simulation. Computational Geosciences, pp.1-18.

Rin, R., Tomin, P., Garipov, T., Voskov, D. and Tchelepi, H., 2017, February. General implicit coupling framework for multi-physics problems. In SPE Reservoir Simulation Conference. Society of Petroleum Engineers.

Schlumberger: ECLIPSE 2013.2 Technical Description (2013)

Senecal, J.P. and Ji, W., 2017. Approaches for mitigating over-solving in multiphysics simulations. International Journal for Numerical Methods in Engineering, 112(6), pp.503-528.

Trangenstein, J.A. and Bell, J.B., 1989a. Mathematical structure of the black-oil model for petroleum reservoir simulation. SIAM Journal on Applied Mathematics, 49(3), pp.749-783.

Trangenstein, J.A. and Bell, J.B., 1989b. Mathematical structure of compositional reservoir simulation. SIAM journal on scientific and statistical computing, 10(5), pp.817-845.

Voskov, D.V. and Tchelepi, H.A., 2012. Comparison of nonlinear formulations for two-phase multi-component EoS based simulation. Journal of Petroleum Science and Engineering, 82, pp.101-111.

Watts, J.W., 1986. A compositional formulation of the pressure and saturation equations. SPE (Society of Petroleum Engineers) Reserv. Eng.;(United States), 1(3).

Walker, H.F. and Ni, P., 2011. Anderson acceleration for fixed-point iterations. SIAM Journal on Numerical Analysis, 49(4), pp.1715-1735.

Watanabe, S., Li, Z., Bratvedt, K., Lee, S.H. and Natvig, J., 2016, August. A Stable Multi-phase Nonlinear Transport Solver with Hybrid Upwind Discretization in Multiscale Reservoir Simulator. In ECMOR XV-15th European Conference on the Mathematics of Oil Recovery (pp. cp-494). European Association of Geoscientists \& Engineers.

Younis, R. and Aziz, K., 2007, January. Parallel automatically differentiable data-types for next-generation simulator development. In SPE Reservoir Simulation Symposium. Society of Petroleum Engineers.

Younis, R., Tchelepi, H.A. and Aziz, K., 2010. Adaptively Localized Continuation-Newton Method--Nonlinear Solvers That Converge All the Time. SPE Journal, 15(02), pp.526-544.

Zhou, Y., Tchelepi, H.A. and Mallison, B.T., 2011, January. Automatic differentiation framework for compositional simulation on unstructured grids with multi-point discretization schemes. In SPE Reservoir Simulation Symposium. Society of Petroleum Engineers.


\end{document}